\documentclass[11pt]{amsart}

\usepackage[T1]{fontenc}
\usepackage{lmodern}
\usepackage[margin=1in]{geometry}
\usepackage{amsmath,amssymb,mathtools}
\usepackage{enumitem}
\setlist{itemsep=2pt,topsep=4pt,parsep=0pt}
\usepackage{booktabs}
\usepackage{graphicx}
\usepackage{tikz}
\usepackage{microtype}
\usepackage[hidelinks]{hyperref}

\numberwithin{equation}{section}
\allowdisplaybreaks[1]

\newtheorem{theorem}{Theorem}[section]
\newtheorem{proposition}[theorem]{Proposition}
\newtheorem{lemma}[theorem]{Lemma}
\newtheorem{corollary}[theorem]{Corollary}

\newtheorem{definition}[theorem]{Definition}

\theoremstyle{remark}
\newtheorem{remark}[theorem]{Remark}
\newtheorem{example}[theorem]{Example}

\newcommand{\C}{\mathbb C}

\newcommand{\Z}{\mathbb Z}
\newcommand{\R}{\mathbb R}
\newcommand{\E}{\mathbb E}

\newcommand{\calD}{\mathcal D}
\newcommand{\calA}{\mathcal A}
\newcommand{\calF}{\mathcal F}
\newcommand{\calM}{\mathcal M}
\newcommand{\dd}{\mathrm d}

\newcommand{\ket}[1]{\lvert #1\rangle}
\newcommand{\ad}{\operatorname{ad}}

\newcommand{\supp}{\operatorname{supp}}
\newcommand{\fall}[2]{#1^{\downarrow #2}}

\title[A Difference Operator Method for Partitions with Respect to the Gelfand Measure]{A Difference Operator Method for Partitions with Respect to the Gelfand Measure}
\author[Guo-Niu Han]{Guo-Niu Han}
\address{Universit\'e de Strasbourg, CNRS, IRMA UMR 7501,
  F-67000 Strasbourg, France}
\email{guoniu.han@unistra.fr}
\author[Huan Xiong]{Huan Xiong}
\address{Institute of Advanced Study of Mathematics,
  Harbin Institute of Technology, Heilongjiang 150001, China}
\email{huan.xiong.math@gmail.com}
\date{November 30, 2026}

\begin{document}

\begin{abstract}
Stanley proved Han's polynomiality conjecture for Plancherel averages.
We develop a difference-operator method for the corresponding
single-tableau sums. Pointwise nilpotence and explicit highest differences
give generating functions \(P(x)e^{x+x^2/2}\), with \(P(x)\) polynomial,
for every fixed inner diagram. We determine all coefficients of the
factorial content and hook polynomials and prove integrality of their
derivatives. From the same difference-operator identities, we derive
the limit shape and joint Gaussian fluctuations for fixed inner diagrams,
recovering the classical Gelfand laws when the inner diagram is empty.
For fixed-prefix Gelfand measures, we obtain uniform moment
error bounds and the first two content-variance corrections.
\end{abstract}

\maketitle

\section{Introduction}\label{sec:introduction}

For a partition \(\lambda\) of \(n\), let \(f^\lambda\) count its
standard tableaux. The Plancherel measure assigns mass
\((f^\lambda)^2/n!\) to \(\lambda\). Han conjectured that the averages
of even hook-length power sums under this measure are polynomials in
the diagram size \cite{HanConjectures}. Stanley proved Han's conjecture
and established the more general polynomiality theorem for every
symmetric function of squared hook lengths \cite{StanleyHooks}.
His theorem also treats functions symmetric separately in squared hook
lengths, contents and shifted parts
\cite[Theorems~3.3, 4.3 and~4.4]{StanleyHooks}.
Explicit evaluations include Okada's hook identity, proved by Panova
\cite{Panova}, and the content identity of Fujii, Kanno, Moriyama and
Okada \cite{FKMO}. Difference-operator proofs and skew extensions were
developed by Han and Xiong \cite{HanXiong}.

We study the corresponding structure for the single-tableau weight
\(f^\lambda\), where the polynomial occurs in the exponential
generating variable \(x\).
\[
 \sum_{n\ge0}\frac{x^n}{n!}\sum_{\lambda\vdash n}f^\lambda g(\lambda)
       =P_g(x)e^{x+x^2/2}.
\]
Here \(g\) is a complex-valued function on partitions, and \(P_g(x)\)
is a polynomial for the classes of observables studied below.
The pointwise statement behind this identity is our main object.
Young's lattice orders partitions by inclusion of their Young diagrams
\cite[Section~2]{StanleyDifferential}.
The rank of a partition is its number of boxes.
For a partition \(\lambda\), define the signed up--down difference operator by
\begin{equation}\label{eq:D-intro}
 Dg(\lambda)=\sum_{\lambda^+}\bigl(g(\lambda^+)-g(\lambda)\bigr)
       +\sum_{\lambda^-}\bigl(g(\lambda)-g(\lambda^-)\bigr).
\end{equation}
Here \(\lambda^+\) and \(\lambda^-\) range over the partitions obtained
from \(\lambda\) by adding and removing one box, respectively.
For the observables below we give an explicit bound \(R\) with
\(D^{R+1}g=0\) on the entire lattice. A highest-symbol formula computes
\(D^Rg\) by evaluating the terms of largest allowed weight in \(g\).
The polynomial factor is then determined by finitely many local
differences. We treat squared hook lengths, contents and shifted parts
jointly and allow any fixed inner partition.

Using the Gelfand weights and involution normalization
\cite[Definition~2 and Section~3.1]{Meliot}, for \(x\ge0\) we define
the Poissonized Gelfand expectation by
\begin{equation}\label{eq:poissonized-motivation}
 \mathbb E_xg=e^{-x-x^2/2}
          \sum_\lambda f^\lambda g(\lambda)\frac{x^{|\lambda|}}{|\lambda|!}.
\end{equation}
For general \(g\) the expression is read coefficientwise; its analytic
interpretation is used when the series is absolutely convergent.
The weights in this expectation sum to one by
\eqref{eq:involution-egf}. This normalization is also the exponential
specialization of the Schur measure with one free boundary
\cite[Section~2.1, (2.3)]{BBNV}. Its rank weights are
\(I_nx^n/(n!I(x))\). We say \(g\) has polynomial growth if
\(|g(\lambda)|\le C(1+|\lambda|)^A\) for some constants \(C,A\ge0\)
and all \(\lambda\). For such \(g\), tableau branching gives
\begin{equation}\label{eq:poissonized-differentiation-intro}
 \frac{d}{dx}\mathbb E_xg=\mathbb E_x(Dg),
\end{equation}
as proved in Lemma~\ref{clt:poisson-derivative}. If \(D^{R+1}g=0\),
this gives \(P_g(x)=\sum_{j=0}^R(D^jg)(\varnothing)x^j/j!\).
Theorem~\ref{thm:skew-exponential-polynomiality} gives the corresponding
formula for a fixed inner diagram.

We make this structure explicit in two ways. First, we determine all
coefficients of the factorial content and factorial hook polynomials by
finite binomial sums, and prove integrality after differentiation.
Second, we quantify the first inner-shape correction to joint Gaussian
fluctuation moments at a common outer rank. Stanley's relative-dimension
expansion and the classical joint Gaussian limits identify its leading
covariance. The operator identities provide uniform remainder estimates
and the next content-variance term. We determine the correction's exact
order, give polynomial tests detecting it, and distinguish changes of
mean from corrections that remain after exact centering.

\subsection{Diagrams, corner observables and tableau measures}\label{sec:diagram-notation}

A partition \(\lambda=(\lambda_1,\lambda_2,\ldots)\) is a weakly
decreasing sequence of nonnegative integers with finitely many positive
parts. Let \(\mathbb Y\) be the set of partitions and \(\mathbb Y_n\) its subset of size \(n\).
We write \(|\lambda|=\sum_i\lambda_i\), \(\ell(\lambda)\) for its
number of positive parts, \(\lambda\vdash n\) for \(|\lambda|=n\), and
\(\varnothing\) for the empty partition. For a nonnegative integer \(r\),
\(\fall{a}{r}=\prod_{j=0}^{r-1}(a-j)\) is the falling factorial, with
\(\fall{a}{0}=1\). We use \(\binom ar=\fall{a}{r}/r!\) for these \(r\)
and set the binomial coefficient to zero when the lower argument is
negative or is not an integer. In particular, \(\binom ar=0\) for
integers \(0\le a<r\). Throughout, \([x^j]F(x)\) denotes the coefficient
of \(x^j\) in a power or Laurent series. We write \(\mathbf1_A\) for the
indicator of a condition \(A\).
We identify \(\lambda\) with
the boxes \((i,j)\) such that \(1\le j\le\lambda_i\), in rows numbered
from the top. The content of \((i,j)\) is \(c(i,j)=j-i\). Its hook
length is \(h_\lambda(i,j)=\lambda_i-j+\lambda'_j-i+1\), where
\(\lambda'=(\lambda'_1,\lambda'_2,\ldots)\) is the conjugate partition
with \(\lambda'_j=\#\{i:\lambda_i\ge j\}\).
The hook at a box \(\square=(i,j)\) consists of that box, the boxes
to its right in row \(i\), and the boxes below it in column \(j\);
its number of boxes is \(h_\lambda(\square)\). Put
\(H_\lambda=\prod_{\square\in\lambda}h_\lambda(\square)\), with
\(H_\varnothing=1\).

For \(\mu\subseteq\lambda\), the skew shape \(\lambda/\mu\) consists
of the boxes of \(\lambda\) outside \(\mu\), in their original
coordinates. We write \(m=|\mu|\) for the inner size,
\(n=|\lambda/\mu|\) for the number of added boxes, and \(N=|\lambda|=n+m\)
for the outer size. A standard tableau of shape \(\lambda/\mu\) fills
these boxes by \(1,\ldots,n\) in increasing order along rows and down
columns. Its number \(f^{\lambda/\mu}\) also counts chains of successive
box additions from \(\mu\) to \(\lambda\). We use \(f^{\mu/\mu}=1\),
\(f^{\lambda/\mu}=0\) if \(\mu\not\subseteq\lambda\), and
\(f^\lambda=f^{\lambda/\varnothing}=|\lambda|!/H_\lambda\)
\cite{FRT,StanleyEC2}. Powers of degree zero are interpreted as \(1\),
including at content \(0\). For integers \(j,r\ge0\), put
\[
 C_j(\lambda)=\sum_{\square\in\lambda}c(\square)^j,\qquad
 C_j(\lambda/\mu)=C_j(\lambda)-C_j(\mu),\qquad
 \mathcal H_r(\lambda)=\sum_{\square\in\lambda}h_\lambda(\square)^{2r}.
\]
Both \(C_0(\lambda)\) and \(\mathcal H_0(\lambda)\) equal
\(|\lambda|\). The difference \(\mathcal H_r(\lambda)-\mathcal H_r(\mu)\)
includes the changes in the hook lengths of boxes already in \(\mu\).

\begin{definition}[Corner observables]\label{def:corner-observables-intro}
Let \(X(\lambda)\) and \(Y(\lambda)\) be the contents of boxes that can
respectively be added to or removed from \(\lambda\) while leaving a
Young diagram. We call these the addable and removable corner contents, respectively.
For \(x\in X(\lambda)\) and \(y\in Y(\lambda)\), write
\(\lambda^{x+}\) and \(\lambda^{y-}\) for the partitions obtained
by adding or deleting the corresponding box. For a polynomial
\(\psi\in\mathbb C[z]\) and integer \(d\ge0\), define
\begin{equation}\label{eq:corner-powers-intro}
\begin{aligned}
 Q_\psi(\lambda)
   &=\sum_{a\in X(\lambda)}\psi(a)-\sum_{b\in Y(\lambda)}\psi(b),\\
 q_d(\lambda)=Q_{z^d}(\lambda)
   &=\sum_{a\in X(\lambda)}a^d-\sum_{b\in Y(\lambda)}b^d .
\end{aligned}
\end{equation}
In particular \(q_0=1\), \(q_1=0\) and
\(q_2(\lambda)=2|\lambda|\) (see \cite[(6.2)--(6.3)]{HanXiong}).

\end{definition}

\begin{figure}[htbp]
\centering
\begin{minipage}[t]{0.32\textwidth}
\centering\vspace{0pt}
\begin{tikzpicture}[x=0.59cm,y=-0.59cm,
  every node/.style={font=\small,inner sep=0pt}]
 \path[use as bounding box] (-0.55,-0.55) rectangle (6.2,4.15);
 \foreach \x/\y in {1/0,2/0,3/0,4/0,1/1,1/2}
   \fill[blue!14] (\x,\y) rectangle ++(1,1);
 \fill[blue!28] (1,0) rectangle (2,1);
 \foreach \r/\len in {1/5,2/3,3/2}
   \foreach \j in {1,...,\len}
     \draw[black!65,line width=0.4pt]
       (\j-1,\r-1) rectangle (\j,\r);
 \draw[blue!65!black,line width=0.9pt]
   (1,0)--(5,0)--(5,1)--(2,1)--(2,3)--(1,3)--cycle;
 \draw[black,line width=1pt] (1,0) rectangle (2,1);
 \foreach \x/\y/\h in
   {0.5/0.5/7,2.5/0.5/4,3.5/0.5/2,4.5/0.5/1,
    0.5/1.5/4,1.5/1.5/3,2.5/1.5/1,0.5/2.5/2,1.5/2.5/1}
   \node at (\x,\y) {\(\h\)};
 \node[font=\small\bfseries] at (1.5,0.5) {\(\boldsymbol6\)};
 \foreach \j in {1,...,5}
   \node[font=\scriptsize,anchor=south] at (\j-0.5,-0.13) {\(\j\)};
 \foreach \i in {1,...,3}
   \node[font=\scriptsize,anchor=east] at (-0.13,\i-0.5) {\(\i\)};
\end{tikzpicture}
\par\small (a) Hook lengths
\end{minipage}\hfill
\begin{minipage}[t]{0.32\textwidth}
\centering\vspace{0pt}
\begin{tikzpicture}[x=0.59cm,y=-0.59cm,
  every node/.style={font=\small,inner sep=0pt}]
 \path[use as bounding box] (-0.55,-0.55) rectangle (6.2,4.15);
 \foreach \r/\len in {1/5,2/3,3/2}
   \foreach \j in {1,...,\len}
     \draw[black!65,line width=0.4pt]
       (\j-1,\r-1) rectangle (\j,\r);
 \fill[blue!14] (1,0) rectangle (2,1);
 \draw[black,line width=1pt] (1,0) rectangle (2,1);
 \foreach \x/\y/\c in
   {0.5/0.5/0,1.5/0.5/1,2.5/0.5/2,3.5/0.5/3,4.5/0.5/4,
    0.5/1.5/-1,1.5/1.5/0,2.5/1.5/1,0.5/2.5/-2,1.5/2.5/-1}
   \node at (\x,\y) {\(\c\)};
\end{tikzpicture}
\par\small (b) Contents
\end{minipage}\hfill
\begin{minipage}[t]{0.32\textwidth}
\centering\vspace{0pt}
\begin{tikzpicture}[x=0.59cm,y=-0.59cm,
  every node/.style={font=\small,inner sep=0pt}]
 \path[use as bounding box] (-0.55,-0.55) rectangle (6.2,4.15);
 \foreach \x/\y in {4/0,2/1,1/2}
   \fill[black!16] (\x,\y) rectangle ++(1,1);
 \foreach \r/\len in {1/5,2/3,3/2}
   \foreach \j in {1,...,\len}
     \draw[black!50,line width=0.4pt]
       (\j-1,\r-1) rectangle (\j,\r);
 \foreach \x/\y/\c in {5/0/5,3/1/2,2/2/0,0/3/-3} {
   \fill[white] (\x,\y) rectangle ++(1,1);
   \draw[blue!65!black,densely dashed,line width=0.8pt]
     (\x,\y) rectangle ++(1,1);
   \node[text=blue!65!black] at (\x+0.5,\y+0.5) {\(\c\)};
 }
 \foreach \x/\y/\c in {4/0/4,2/1/1,1/2/-1} {
   \draw[black,line width=0.8pt] (\x,\y) rectangle ++(1,1);
   \node at (\x+0.5,\y+0.5) {\(\c\)};
 }
\end{tikzpicture}
\par\small (c) Corner contents
\end{minipage}
\caption{The Young diagram of \(\lambda=(5,3,2)\), with rows numbered
from the top. In (a), the entries are hook lengths, and the shaded
boxes form the hook of the marked box \(\square=(1,2)\), of length
\(h_\lambda(\square)=6\). In (b), the entries are contents, with
\(c(\square)=1\) in the marked box. In (c), dashed boxes are addable
and shaded boxes are removable; their entries give the two corner
sets in Definition~\ref{def:corner-observables-intro}.}
\label{fig:partition-hook-content-corners}
\end{figure}
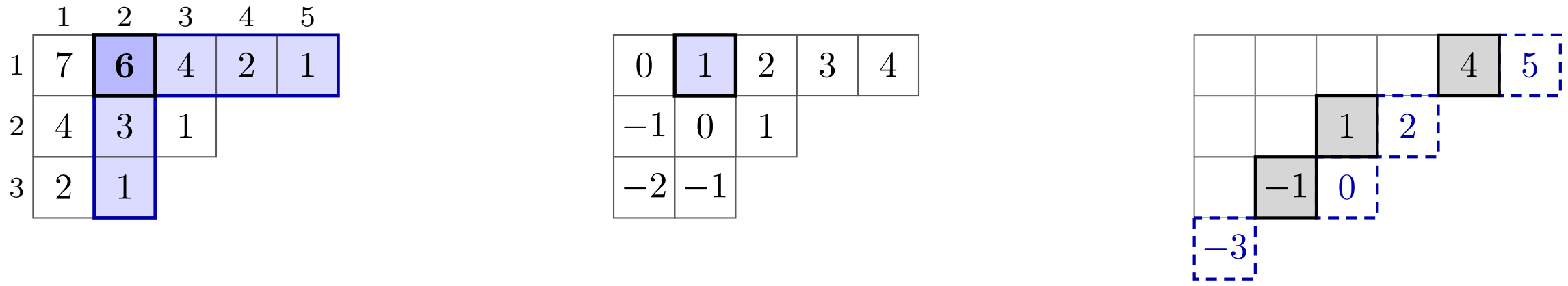

Figure~\ref{fig:partition-hook-content-corners} illustrates these
definitions for \(\lambda=(5,3,2)\), whose Young diagram has ten boxes.
For this partition,
\[
 X(\lambda)=\{-3,0,2,5\},\qquad Y(\lambda)=\{-1,1,4\},
\]
so \(q_0=1\), \(q_1=0\), and \(q_2=20=2|\lambda|\).

Lemma~\ref{lem:corner-increments} records the local corner changes.

Fix an inner partition \(\mu\). For \(\lambda\supseteq\mu\) with
\(|\lambda/\mu|=n\), the single-tableau skew measure is
\[
 M_{\mu,n}(\lambda)=\frac{f^{\lambda/\mu}}{Z_\mu(n)},\qquad
 Z_\mu(n)=\sum_{\substack{\rho\supseteq\mu\\|\rho/\mu|=n}}f^{\rho/\mu}.
\]
It gives mass zero to other partitions. Equivalently, choose a standard
tableau uniformly over all tableaux of size \(N\) and condition on its
first \(m\) entries forming a fixed tableau of shape \(\mu\).
The conditional law of its outer shape is \(M_{\mu,n}\), as shown in
Lemma~\ref{lem:fixed-tableau-prefix}. The empty-base specialization
is the Gelfand measure
\begin{equation}\label{eq:Young-measure}
 M_n(\lambda)=M_{\varnothing,n}(\lambda)=\frac{f^\lambda}{I_n},
 \qquad I_n=Z_\varnothing(n),\quad \lambda\vdash n .
\end{equation}
By the Robinson--Schensted correspondence
\cite{Robinson,Schensted,StanleyEC2}, \(I_n\) counts the involutions
in the symmetric group \(S_n\). Decomposing an involution into fixed points and
transpositions gives \(I_0=I_1=1\), \(I_{n+1}=I_n+nI_{n-1}\)
for \(n\ge1\), and
\begin{equation}\label{eq:involution-egf}
 I(x)=\sum_{n\ge0}I_n\frac{x^n}{n!}=e^{x+x^2/2}.
\end{equation}

\subsection{Exponential polynomiality and difference degrees}

Our difference-operator method gives pointwise identities for functions
on diagrams and their mixed products. For example,
\(\mathbb E_xC_1=0\) by diagram conjugation, although \(C_1\) is not the zero
function on diagrams. The commutator identities give explicit nilpotence
bounds, exact highest-symbol identities and local coefficient formulas.
In the skew setting, they are used for inserted correlations and their
uniform remainder estimates.

A nonzero function \(g\) has \(D\)-degree \(R\) if \(D^Rg\ne0\)
and \(D^{R+1}g=0\), with \(D^0g=g\).
Such functions and the zero function are called \(D\)-polynomials.
We write \(\deg_Dg\) for the \(D\)-degree, with \(\deg_D0=-\infty\),
and write \(\mathbf1\) for the constant function equal to one. Put
\[
 S_{\mu,g}(n)=
 \sum_{\substack{\lambda\supseteq\mu\\|\lambda/\mu|=n}}
                   f^{\lambda/\mu}g(\lambda).
\]
\begin{theorem}[Exponential polynomiality for a fixed skew base]
\label{thm:skew-exponential-polynomiality}
Let \(g\) be a function on all partitions, let \(R\) be a nonnegative
integer with \(D^{R+1}g=0\), and fix an inner partition \(\mu\). Then
\begin{equation}\label{eq:skew-polynomiality-intro}
 \sum_{n\ge0}S_{\mu,g}(n)\frac{x^n}{n!}
                   =P_{\mu,g}(x)e^{x+x^2/2},
\end{equation}
where
\begin{equation}\label{eq:skew-polynomial-factor-intro}
 P_{\mu,g}(x)=
 \sum_{\eta\subseteq\mu}\frac{f^{\mu/\eta}x^{|\mu/\eta|}}{|\mu/\eta|!}
       \sum_{j=0}^R\frac{(D^jg)(\eta)}{j!}x^j
\end{equation}
has degree at most \(R+|\mu|\). In particular,
\[
 \sum_\lambda f^\lambda g(\lambda)\frac{x^{|\lambda|}}{|\lambda|!}
  =e^{x+x^2/2}\sum_{j=0}^R\frac{(D^jg)(\varnothing)}{j!}x^j .
\]
\end{theorem}
All generating-function identities are understood coefficientwise.
Section~\ref{sec:skew-transfer} proves the stronger form in
Lemma~\ref{lem:skew-transfer}, including the sums weighted by \(D^kg\).
The result reduces polynomiality to pointwise nilpotence.
We next establish the requisite bounds and evaluate their highest
nonzero differences.

For an indexing partition \(\nu\), define on all partitions
\begin{equation}\label{eq:skew-content-extension-intro}
 F_{\nu,\mu}(\lambda)
    =\prod_{i=1}^{\ell(\nu)}
            \bigl(C_{2\nu_i}(\lambda)-C_{2\nu_i}(\mu)\bigr).
\end{equation}
On \(\lambda\supseteq\mu\) this is the product of the even content
power sums of \(\lambda/\mu\). The global extension allows \(D\)
to delete a box of \(\mu\). All empty products equal one.

A polynomial in \(q_j\), \(j\ge2\), has \emph{corner weight}
at most \(w\) if each of its monomials \(\prod_jq_j^{e_j}\)
has \(\sum_jje_j\le w\). Constants have weight zero. Write
\(A_{[w]}\) for the sum of monomials of weight exactly \(w\), and put
\[
 \varepsilon(A_{[w]})=A_{[w]}(a_2,a_3,\ldots),\qquad
 a_{2r}=\binom{2r}{r},\quad a_{2r+1}=0\quad(r\ge0).
\]

\begin{theorem}[Pointwise polynomiality and highest differences]
\label{thm:degree-intro}
Let \(w\) be a nonnegative integer and \(A\) a corner polynomial of
weight at most \(w\). Then
\[
 D^{w+1}A=0,\qquad
 D^wA=w!\varepsilon(A_{[w]})\mathbf1,\qquad
 [x^{m+w}]P_{\mu,A}=\frac{f^\mu}{m!}\varepsilon(A_{[w]}),
 \quad m=|\mu|.
\]
If \(\varepsilon(A_{[w]})\ne0\), both degree bounds are exact.
The highest-symbol evaluation is multiplicative in the following sense.
If \(B\) has weight at most \(v\), then
\(\varepsilon((AB)_{[w+v]})=\varepsilon(A_{[w]})\varepsilon(B_{[v]})\).
In particular, let \(\psi_1,\ldots,\psi_s\) be nonzero polynomials of degrees
\(d_1,\ldots,d_s\), let \(G=\prod_{i=1}^sQ_{\psi_i}\), and put
\(w=\sum_{i=1}^s d_i\). Then
\[
 D^{w+1}G=0.
\]
For a fixed inner partition \(\mu\) and an indexing partition \(\nu\),
put \(R=2|\nu|+2\ell(\nu)\). Then \(F_{\nu,\mu}\) has exact
\(D\)-degree \(R\), with
\begin{equation}\label{eq:skew-content-degree-intro}
 D^RF_{\nu,\mu}
   =R!\prod_{i=1}^{\ell(\nu)}
           \frac{(2\nu_i)!}{(\nu_i+1)!^2}\,\mathbf1,
 \qquad D^{R+1}F_{\nu,\mu}=0 .
\end{equation}
\end{theorem}

The operator form is Theorem~\ref{thm:corner-symbol}.
Its proof uses the lowering of polynomial coefficient degree under
commutation with the Heisenberg generators, the lifted particle shifts
used in Section~\ref{sec:commutator-calculus}
\cite[Sections~2.2.4--2.2.6]{OkounkovPandharipande}
\cite[Section~1.3]{BHY}.
The discussion after Lemma~\ref{lem:highest-central-commutator} also
derives the highest commutator directly from the Lascoux--Thibon
vertex-operator relation \cite[(16)]{LascouxThibon}.
If \(a_i\) is the leading coefficient of \(\psi_i\), then
\begin{equation}\label{eq:corner-product-top-intro}
 D^wG=w!\left(\prod_{i=1}^s a_i\binom{d_i}{d_i/2}\right)\mathbf1,
\end{equation}
where the binomial coefficient is zero when \(d_i\) is odd.
In fact \(\deg_DG\le2\sum_i\lfloor d_i/2\rfloor\);
when all \(d_i\) are even, the degree is exactly \(w\).
Corollary~\ref{cor:corner-parity-degree} proves these refinements.

\subsection{Stanley's polynomiality and its Gelfand analogue}

Let \(\Lambda_{\mathbb Q}=\mathbb Q[p_1,p_2,\ldots]\) be the ring
of symmetric functions, with power sums \(p_j(z)=\sum_i z_i^j\) and
\(p_\nu=\prod_i p_{\nu_i}\), with \(p_\varnothing=1\).
Evaluation on a finite list sets all remaining variables to zero
\cite[Chapter~I]{MacdonaldSF}. For \(|\lambda|=N\), define
\[
 \boldsymbol h_\lambda=(h_\lambda(\square)^2:\square\in\lambda),
 \quad \boldsymbol c_\lambda=(c(\square):\square\in\lambda),
 \quad \boldsymbol a_\lambda=(\lambda_i+N-i:1\le i\le N).
\]
Here \(\lambda_i=0\) for \(i>\ell(\lambda)\); all lists are empty
for \(\lambda=\varnothing\).
The entries of \(\boldsymbol a_\lambda\) are Stanley's shifted parts.
An element of \(\Lambda_{\mathbb Q}^{\otimes3}\), the tensor product
of three copies of this ring, is a finite linear combination of
\(p_\nu(\boldsymbol h)p_\kappa(\boldsymbol c)p_\rho(\boldsymbol a)\).
It is symmetric separately in the three lists. Each indexing partition
may be empty.

\begin{lemma}[Stanley's Plancherel polynomiality; \cite{StanleyHooks}]
\label{lem:stanley-polynomiality}
For every \(F\in\Lambda_{\mathbb Q}^{\otimes3}\),
\[
 \frac1{n!}\sum_{\lambda\vdash n}(f^\lambda)^2
       F(\boldsymbol h_\lambda;\boldsymbol c_\lambda;
                                   \boldsymbol a_\lambda)\in\mathbb Q[n].
\]
\end{lemma}
This is \cite[Theorems~3.3, 4.3 and~4.4]{StanleyHooks}; taking the
second and third factors constant gives squared-hook polynomiality.
For the single-tableau weight we obtain the following pointwise bound.

The generating-function assertion in the following theorem is also a
consequence of known character formulas, as explained immediately after
the statement. The pointwise vanishing
\(D^{\mathfrak w(F)+1}g_F=0\) and the finite local formula for its
polynomial factor are the contributions of the difference-operator
calculus. These hold on the whole lattice before an inner diagram is
chosen.

\begin{theorem}[Gelfand analogue of Stanley's polynomiality]
\label{thm:gelfand-stanley-polynomiality}
Let \(F\in\Lambda_{\mathbb Q}^{\otimes3}\) and define, on all partitions,
\[
 g_F(\lambda)=F(\boldsymbol h_\lambda;\boldsymbol c_\lambda;
                                           \boldsymbol a_\lambda).
\]
For the unique power-sum expansion
\[
 F=\sum_{\nu,\kappa,\rho}a_{\nu,\kappa,\rho}
               p_\nu(\boldsymbol h)p_\kappa(\boldsymbol c)p_\rho(\boldsymbol a),
\]
put
\[
 \mathfrak w(F)=\max_{a_{\nu,\kappa,\rho}\ne0}
 \bigl(2|\nu|+2\ell(\nu)+|\kappa|+2\ell(\kappa)
                         +2|\rho|+2\ell(\rho)\bigr),\qquad
 \mathfrak w(0)=0.
\]
Then \(D^{\mathfrak w(F)+1}g_F=0\).
For every fixed inner partition \(\mu\), there is a polynomial
\(P_{\mu,F}\in\mathbb Q[x]\) of degree at most \(|\mu|+\mathfrak w(F)\)
such that
\begin{equation}\label{eq:gelfand-stanley-egf}
 \sum_{n\ge0}\frac{x^n}{n!}
 \sum_{\substack{\lambda\supseteq\mu\\|\lambda/\mu|=n}}
        f^{\lambda/\mu}g_F(\lambda)=P_{\mu,F}(x)e^{x+x^2/2}.
\end{equation}
Its polynomial factor is given by
\eqref{eq:skew-polynomial-factor-intro} with \(g=g_F\) and
\(R=\mathfrak w(F)\). In particular,
\begin{equation}\label{eq:gelfand-stanley-ordinary}
 \sum_{n\ge0}\frac{x^n}{n!}\sum_{\lambda\vdash n}f^\lambda g_F(\lambda)
 =e^{x+x^2/2}\sum_{j=0}^{\mathfrak w(F)}
                       \frac{(D^jg_F)(\varnothing)}{j!}x^j.
\end{equation}
\end{theorem}

Each of the three lists has \(N\) entries, with hook lengths at most
\(N\), absolute contents at most \(N\), and shifted parts at most
\(2N\). The finite power-sum expansion therefore gives
\(|g_F(\lambda)|\le C_F(1+N)^{\mathfrak w(F)}\).
Thus these observables have polynomial growth and their Poissonized
sums are entire in \(x\).

Section~\ref{sec:gelfand-stanley-proof} proves the pointwise assertion
and the local formula by expressing
\(g_F\) in the corner powers with weight at most \(\mathfrak w(F)\).
The hooks and shifted parts are those of the outer diagram, with shifted
parts computed at size \(N\). Example~\ref{ex:shifted-parts-mixed-sharp}
gives a mixed product whose highest difference attains the bound.

\paragraph{A character-theoretic derivation of the generating function.}
For a partition \(\rho\), let \(m_j(\rho)\) denote the number of parts
equal to \(j\), and set \(z_\rho=\prod_{j\ge1}j^{m_j(\rho)}m_j(\rho)!\).
The irreducible representation of \(S_N\) indexed by \(\lambda\vdash N\)
has dimension \(f^\lambda\) \cite[Chapter~I, Section~7]{MacdonaldSF}.
Write \(\chi^\lambda(\nu)\) for its character evaluated at cycle type
\(\nu\vdash N\).
If \(N\ge|\rho|\), define the normalized central character by
\begin{equation}\label{eq:central-character}
 \Sigma_\rho(\lambda)=\fall{N}{|\rho|}
       \frac{\chi^\lambda(\rho\cup1^{N-|\rho|})}{f^\lambda}.
\end{equation}
Here \(\rho\cup1^a\) appends \(a\) parts equal to one. We set
\(\Sigma_\rho(\lambda)=0\) when \(N<|\rho|\).
The \(\Sigma_\rho\) form a linear basis of the algebra of corner
polynomials, and those with \(|\rho|+\ell(\rho)\le w\) form a basis of
its weight-\(\le w\) subspace \cite[Section~2]{Meliot}. Let \(R_\rho\)
count the square roots of a permutation whose cycle type consists of the
parts of \(\rho\) of size at least two, in the symmetric group on their
total number of letters. By M\'eliot's finite-rank formula
\cite[Section~3.1]{Meliot},
\[
 \mathbb E_x\Sigma_\rho
   =R_\rho x^{|\rho|}
          \frac{I^{(m_1(\rho))}(x)}{I(x)}.
\]
Indeed, in size \(N\) the first moment is
\(\fall{N}{|\rho|} R_\rho I_{N-|\rho|+m_1(\rho)}/I_N\) for
\(N\ge|\rho|\), and summation gives the displayed identity. Its degree
is at most \(|\rho|+m_1(\rho)\), since \(I'=(1+x)I\).

For \(|\mu|=m\), set \(c_\mu=f^\mu/m!\). The shifted Schur function
\(s^*_\mu\), evaluated on a diagram \(\lambda\), satisfies
\[
 s_\mu^*(\lambda)
   =\fall{N}{m}\frac{f^{\lambda/\mu}}{f^\lambda}
   =\sum_{\rho\vdash m}\frac{\chi^\mu(\rho)}{z_\rho}\Sigma_\rho(\lambda)
   \quad(N\ge m).
\]
It vanishes in sizes below \(m\). The first formula is
\cite[Theorem~8.1]{OkounkovOlshanskiShifted}, and the second follows from
\cite[Theorem~3.1, equation~(12)]{StanleySkew}. Thus \(s_\mu^*\) has corner
weight at most \(2m\). Applying the character expectation to
\(g_Fs_\mu^*\) gives a polynomial of degree at most
\(\mathfrak w(F)+2m\). By the definition of \(\mathbb E_x\),
this expectation is the normalized series in
\eqref{eq:gelfand-stanley-egf} multiplied by \(x^m\), so division gives
the stated bound on the generating-function factor. The difference
operator additionally yields the pointwise vanishing and computes
this polynomial from local values of \(D^jg_F\). For example,
\((Dq_2)(\lambda)=4|Y(\lambda)|+2\), \(D^2q_2=4\), and \(D^3q_2=0\).

\subsection{Hook polynomials}

Hook-length identities and polynomial partition averages have been
studied in \cite{HanConjectures,StanleyHooks,OlshanskiAverages,Panova,HanXiong}.
We now determine the highest differences of hook power sums
and their fixed-base differences.

\begin{theorem}[Hook polynomials and fixed-base differences]
\label{thm:hook-polynomiality-intro}
Fix \(\mu\) and nonnegative integers \(r_1,\ldots,r_s\). Set
\[
 F^{\mathrm h}_{\boldsymbol r,\mu}(\lambda)
   =\prod_{i=1}^s\bigl(\mathcal H_{r_i}(\lambda)-\mathcal H_{r_i}(\mu)\bigr),
 \qquad R=\sum_{i=1}^s(2r_i+2).
\]
For \(r\ge0\) put
\begin{equation}\label{eq:hook-leading-constant}
 \beta_r^{\mathrm h}=\frac{1}{2(r+1)^2}
                   \binom{2r}{r}\binom{2r+2}{r+1}.
\end{equation}
Then \(F^{\mathrm h}_{\boldsymbol r,\mu}\) has exact \(D\)-degree \(R\), and
\begin{equation}\label{eq:hook-highest-difference}
 D^R F^{\mathrm h}_{\boldsymbol r,\mu}
   =R!\prod_{i=1}^s\beta_{r_i}^{\mathrm h}\,\mathbf1,
 \qquad D^{R+1}F^{\mathrm h}_{\boldsymbol r,\mu}=0.
\end{equation}
More generally, let \(F=\sum_\nu a_\nu p_\nu\) be a symmetric
function in one alphabet and put
\(g_F(\lambda)=F(h_\lambda(\square)^2:\square\in\lambda)\).
Its \(D\)-degree is at most
\begin{equation}\label{eq:hook-symmetric-degree}
 W(F)=\max_{a_\nu\ne0}\bigl(2|\nu|+2\ell(\nu)\bigr),\qquad W(0)=0.
\end{equation}
Consequently Theorem~\ref{thm:skew-exponential-polynomiality} gives
\(\sum_{n\ge0}S_{\mu,g_F}(n)x^n/n!=I(x)P_{\mu,g_F}(x)\), with
\(\deg P_{\mu,g_F}\le|\mu|+W(F)\).
\end{theorem}

\subsection{Factorial content and hook polynomials}

For both families we give every coefficient of the normalized
single-tableau generating function.

\begin{theorem}[The factorial content polynomial]
\label{thm:factorial-content-intro}
For \(k\ge0\), define
\[
 F_k^{\mathrm c}(\lambda)=\sum_{\square\in\lambda}
                \prod_{r=0}^{k-1}\bigl(c(\square)^2-r^2\bigr).
\]
There is a polynomial \(Q_k\in\mathbb Q[x]\) such that
\begin{equation}\label{eq:factorial-content-egf}
 \sum_\lambda\frac{x^{|\lambda|}}{H_\lambda}F_k^{\mathrm c}(\lambda)
                         =Q_k(x)e^{x+x^2/2}.
\end{equation}
With the generalized binomial coefficients above, one has
\begin{equation}\label{eq:factorial-content-complete}
 Q_k(x)=(2k)!\sum_{d=k}^{2k+1}
         \frac{4^{d-k}}{(d+1)!}\binom{d/2}{d-k}\,x^{d+1}.
\end{equation}
In particular,
\begin{equation}\label{eq:factorial-content-derivative}
 [x^d]Q'_k(x)=\frac{(2k)!}{d!}\,4^{d-k}\binom{d/2}{d-k}
                 \qquad(k\le d\le 2k+1).
\end{equation}
The coefficients of \(x^{k+1},\ldots,x^{2k+2}\) in \(Q_k\) are
strictly positive rational numbers, and all its other coefficients
vanish. The coefficients of \(x^k,\ldots,x^{2k+1}\) in \(Q'_k\) are
strictly positive integers, and all its other coefficients vanish.
In particular,
\begin{align}\label{eq:factorial-content-extreme-coefficients}
 [x^{2k+1}]Q'_k&=\frac1{2k+1}\binom{2k+2}{k+1},\notag\\
 [x^{2k}]Q'_k&=2^{2k},\qquad [x^k]Q'_k=\frac{(2k)!}{k!}.
\end{align}
\end{theorem}

Fujii, Kanno, Moriyama and Okada computed the Plancherel average of
this same statistic \cite[Theorem~A.1]{FKMO}. Han and Xiong gave a
difference-operator proof and, for \(k=1\), a skew marked-content
formula \cite[Theorems~1.7--1.9]{HanXiong}. In our notation the ordinary
Plancherel formula reads
\[
 n!\sum_{\lambda\vdash n}\frac{F_k^{\mathrm c}(\lambda)}{H_\lambda^2}
       =\frac{(2k)!}{(k+1)!}\binom{n}{k+1}.
\]
Its Poisson transform with mean \(x^2\) consists of the single term
\((2k)!\,x^{2k+2}/(k+1)!^2\), exactly the highest term in
\eqref{eq:factorial-content-complete}. The single-tableau weight in
\eqref{eq:factorial-content-egf} gives the additional degrees
\(k+1,\ldots,2k+1\). Section~\ref{sec:factorial-content-proof}
extracts their coefficients from the content exponential and proves
integrality after differentiation.

Replacing the contents by hook lengths and shifting the factors gives a
second family of polynomials.

Here are the finite arrays used in the complete hook formula. For
integers \(p\ge0\) and \(h\), put
\[
 t_{p,h}=\begin{cases}
 4^h\displaystyle\binom{(p-1)/2}{h},&0\le2h\le p,\\
 0,&\text{otherwise},
 \end{cases}
 \qquad \operatorname{Cat}_i=\frac1{i+1}\binom{2i}{i}\quad(i\ge0).
\]
For \(n\ge0\) and integers \(0\le2h\le n\), define
\begin{align*}
 \mathfrak a_{n,h}&=\sum_{p=0}^n\binom np
               \sum_{i=0}^h t_{p,i}t_{n-p,h-i},\\
 \mathfrak e_{n,h}&=\sum_{i=0}^h\binom n{2i}
                          \operatorname{Cat}_i t_{n-2i,h-i}.
\end{align*}
Both arrays are zero outside this range. For \(n\ge2\) and even
\(d\ge0\), put
\begin{equation}\label{eq:hook-short-arrays-intro}
 \mathfrak b_{n,d}=\mathfrak a_{n,n-d/2}
       +\frac12\mathfrak a_{n,n-d/2+1}
       +2\mathfrak e_{n-2,n-d/2}.
\end{equation}
For odd \(d\), set \(\mathfrak b_{n,d}=0\).
These are integer arrays; in particular \(\mathfrak a_{n,h}\) is even
for \(n\ge1\). Section~\ref{sec:hook-all-coefficients} proves these
facts and derives the formula below.

\begin{theorem}[The factorial hook polynomial]\label{thm:factorial-hook-intro}
For \(k\ge0\) put
\(\mathcal J_k(\lambda)=\sum_{\square\in\lambda}
            \prod_{j=1}^k(h_\lambda(\square)^2-j^2)\).
There is a polynomial \(P_k\in\mathbb Q[x]\) of degree \(2k+2\),
divisible by \(x^{k+1}\), such that
\begin{equation}\label{eq:factorial-hook-egf}
 \sum_\lambda\frac{x^{|\lambda|}}{H_\lambda}\mathcal J_k(\lambda)
                =P_k(x)e^{x+x^2/2}.
\end{equation}
For \(k\ge1\), all coefficients are given by
\begin{equation}\label{eq:hook-all-binomial-intro}
 P'_k(x)=\frac{(k!)^2}{2}
 \sum_{n=k+1}^{2k+2}\frac{x^{n-1}}{(n-1)!}
 \sum_{j=\lfloor(n-1)/2\rfloor}^{k}
          (-1)^{k-j}\binom{2j}{j}\mathfrak b_{n,2j+2}.
\end{equation}
For \(k\ge1\), its derivative satisfies
\begin{align}\label{eq:factorial-hook-three-coefficients}
 [x^{2k+1}]P'_k(x)
   &=\frac{1}{k+1}\binom{2k}{k}\binom{2k+2}{k+1},\notag\\
 [x^{2k}]P'_k(x)&=\binom{4k+1}{2k},\qquad
 [x^k]P'_k(x)=\frac{2^k(2k+1)!}{(k+1)!}.
\end{align}
In particular \(P'_k\in x^k\mathbb Z[x]\).
The equivalent local formula is
\begin{equation}\label{eq:factorial-hook-all-coefficients}
 [x^j]P'_k(x)=\frac{D^{j+1}\mathcal J_k(\varnothing)}{j!},
       \qquad j\ge0.
\end{equation}
For \(k=0\), \(P_0=x+x^2\).
\end{theorem}

The Okada--Panova Plancherel average of \(\mathcal J_k\) in size \(n\)
is \(\beta_k^{\mathrm h}\fall{n}{k+1}\)
\cite[Theorem~1]{Panova}. Its Poisson transform with mean \(x^2\) is
\(\beta_k^{\mathrm h}x^{2k+2}\), the highest term of \(P_k\).
Section~\ref{sec:hook-polynomials} obtains all coefficients and proves
integrality by taking a constant term in two auxiliary variables.
Proposition~\ref{prop:odd-hook-counterexample} gives a counterexample
to pointwise polynomiality for the first odd hook power.

\subsection{Skew limit shape and Gaussian fluctuations}

M\'eliot established the limit shape and Gaussian boundary fluctuations
for the ordinary Gelfand measure \cite{Meliot,MeliotKerov}.
Raposo's framework for global and multilevel fluctuations of random
Young tableaux includes this measure
\cite[Theorems~2.19 and~2.22, Example~2.24]{Raposo}.
The classical fixed-prefix enumeration gives the normalization of the
conditioned law \cite{McKayMorseWilf,StanleySkew}. We prove the
fixed-inner-shape limit laws directly from the difference operator \(D\)
in Section~\ref{sec:limit-shape}. The highest differences and exact skew
transfer yield the content moment laws and the uniform limit shape;
polynomial operator conjugation and coefficient extraction give the
joint Gaussian fluctuations. Setting \(\mu=\varnothing\) recovers the
classical ordinary Gelfand results cited above.
The same calculus gives an \(O(N^{-1/2})\) error for each fixed mixed
moment and the inserted identities needed for the inner-shape corrections.
In particular, it yields the two variance terms
\(-2m(m-1)\) and \(6\kappa_3(\mu)/\sqrt N\), with remainder
\(O_\mu(N^{-1})\), in \eqref{eq:response-variance-refined-content};
here \(\kappa_3(\mu)=C_2(\mu)-\binom m2\).
Corollary~\ref{cor:adv-centered-inner-correction} describes the effect
of separate exact expectation centering.

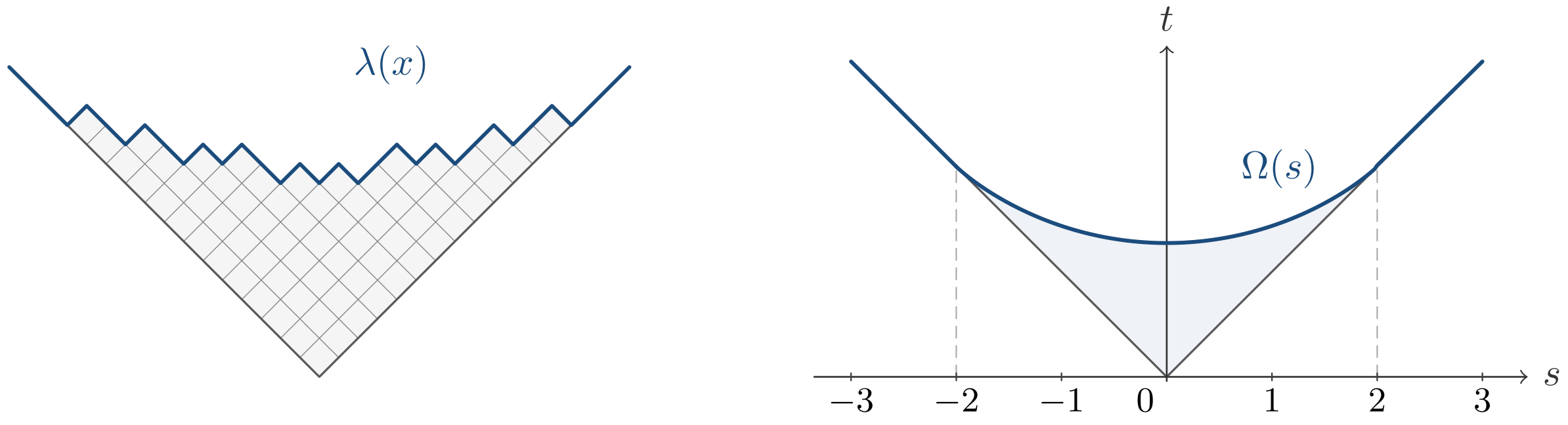
\begin{figure}[!ht]
\centering
\resizebox{0.98\linewidth}{!}{%
\begin{tikzpicture}[line cap=round,line join=round,font=\small]
  \definecolor{profileblue}{RGB}{27,76,125}
  \path[use as bounding box] (0,-0.83) rectangle (14.5,3.35);
  \begin{scope}[shift={(3.15,0)}]
    \begin{scope}[x=0.175cm,y=0.175cm]
      \foreach \row/\len in {1/13,2/11,3/9,4/8,5/6,6/5,7/4,8/4,9/3,10/2,11/2,12/1,13/1} {
        \foreach \col in {1,...,\len} {
          \path[draw=black!45,fill=black!4,line width=0.22pt]
            ({\col-\row},{\col+\row-2}) --
            ({\col-\row+1},{\col+\row-1}) --
            ({\col-\row},{\col+\row}) --
            ({\col-\row-1},{\col+\row-1}) -- cycle;
        }
      }
      \draw[black!65,line width=0.55pt] (-16,16)--(0,0)--(16,16);
      \draw[profileblue,line width=0.95pt] (-16,16) -- (-13,13) -- (-12,14) -- (-11,13) -- (-10,12) -- (-9,13) -- (-8,12) -- (-7,11) -- (-6,12) -- (-5,11) -- (-4,12) -- (-3,11) -- (-2,10) -- (-1,11) -- (0,10) -- (1,11) -- (2,10) -- (4,12) -- (5,11) -- (6,12) -- (7,11) -- (9,13) -- (10,12) -- (12,14) -- (13,13) -- (16,16);
    \end{scope}
    \node[text=profileblue] at (0.65,2.82) {\(\lambda(x)\)};
    \node at (0,-0.66) {(a) An unscaled partition boundary};
  \end{scope}
  \begin{scope}[shift={(10.8,0)},x=0.95cm,y=0.95cm]
    \path[fill=profileblue!7]
       (-2,2) plot[domain=-2:2,samples=181]
       (\x,{2/pi*(\x*asin(\x/2)*pi/180+sqrt(max(0,4-\x*\x)))})
       -- (0,0) -- cycle;
    \draw[black!35,densely dashed,line width=0.35pt] (-2,0)--(-2,2) (2,0)--(2,2);
    \draw[black!65,line width=0.55pt] (-3,3)--(0,0)--(3,3);
    \draw[->,black!80,line width=0.45pt] (-3.35,0)--(3.43,0) node[right] {\(s\)};
    \draw[->,black!80,line width=0.45pt] (0,-0.04)--(0,3.15) node[above] {\(t\)};
    \foreach \t in {-3,-2,-1,1,2,3} {
      \draw[black!75,line width=0.4pt] (\t,0.035)--(\t,-0.035);
      \node[below,inner sep=3pt,font=\footnotesize] at (\t,0) {\(\t\)};
    }
    \node[below left,inner sep=3pt,font=\footnotesize] at (0,0) {\(0\)};
    \draw[profileblue,line width=1pt] (-3,3)--(-2,2)
      plot[domain=-2:2,samples=181]
       (\x,{2/pi*(\x*asin(\x/2)*pi/180+sqrt(max(0,4-\x*\x)))})
      -- (3,3);
    \node[text=profileblue] at (1.07,1.98) {\(\Omega(s)\)};
    \node at (0,-0.695) {(b) The limiting profile};
  \end{scope}
\end{tikzpicture}%
}
\caption{An unscaled partition boundary (left) and the normalized limiting
profile (right). The blue curves are \(y=\lambda(x)\) and
\(t=\Omega(s)\); the black V-shaped lines are \(y=|x|\) and \(t=|s|\),
respectively. The panels use different scales.}
\label{fig:russian-profiles}
\end{figure}

We use Russian coordinates
\(x=u-v\), \(y=u+v\), where the box in row \(i\), column \(j\),
occupies \([j-1,j]\times[i-1,i]\) in the \((u,v)\)-plane.
Each rotated box has area \(2\), and its horizontal center is its content. The upper boundary \(\lambda(x)\) has slopes
\(\pm1\), equals \(|x|\) outside a compact interval, and has slope
jumps \(+2\) at \(X(\lambda)\) and \(-2\) at \(Y(\lambda)\).
For \(N>0\), in normalized coordinates \(s=x/\sqrt N\), \(t=y/\sqrt N\), its graph is
\(t=\bar\lambda_N(s):=N^{-1/2}\lambda(\sqrt N\,s)\).
The Logan--Shepp--Vershik--Kerov shape for the Plancherel measure
\((f^\lambda)^2/N!\), \(\lambda\vdash N\), is \cite{LS,VK,VK3}
\[
 \Omega(s)=
 \begin{cases}
 \displaystyle\frac{2}{\pi}
       \left(s\arcsin(s/2)+\sqrt{4-s^2}\right),&|s|\le2,\\
 |s|,&|s|\ge2.
 \end{cases}
\]
Figure~\ref{fig:russian-profiles} pairs an unscaled finite example of
\(\lambda(x)\) with the normalized limiting profile \(\Omega(s)\).
The black V-shaped lines are \(y=|x|\) on the left and \(t=|s|\) on the
right; \(\Omega(s)=|s|\) for \(|s|\ge2\).
Theorem~\ref{thm:skew-shape-intro} gives this same uniform limit for the
outer and skew profiles under \(M_{\mu,n}\), for fixed \(\mu\) and
\(n\to\infty\).

For \(\lambda\supseteq\mu\) and \(n=|\lambda/\mu|=N-|\mu|>0\),
the skew normalization uses \(s=x/\sqrt n\):
\[
 \begin{aligned}
 \omega_{\lambda/\mu}(x)&=\lambda(x)-\mu(x)+|x|,\\
 \overline\omega_{\lambda/\mu,n}(s)&=n^{-1/2}\omega_{\lambda/\mu}(\sqrt n\,s),\\
 \Delta_{\lambda/\mu,n}^Y(s)&=
 \frac{\sqrt n}{2}\bigl(\overline\omega_{\lambda/\mu,n}(s)-\Omega(s)\bigr).
 \end{aligned}
\]
Here \(\mu(x)\) is the ordinary boundary of the inner partition.
The skew profile keeps the original content coordinates of the remaining
boxes. For \(\mu=\varnothing\), it reduces to \(\lambda(x)\) and \(n=N\).
The factor \(\sqrt n/2\) specifies our fluctuation normalization.
The analogous outer fluctuation at size \(N\) uses
\(\sqrt N(\bar\lambda_N-\Omega)/2\).

With \(a_d\) as above, set
\begin{equation}\label{eq:intro-clt-coefficients}
 b_{dr}=d\binom{d-1}{(d-1-r)/2}\quad(d\ge2,\ r\ge1),
 \qquad \gamma_j=\frac{a_{j+2}}{(j+1)(j+2)}\quad(j\ge0).
\end{equation}
These indices are integers. Define the modified Chebyshev polynomials by
\(u_0(x)=1\), \(u_1(x)=x\), and
\(u_{k+1}(x)=xu_k(x)-u_{k-1}(x)\) for \(k\ge1\). Equivalently,
\(u_k(2\cos\theta)=\sin((k+1)\theta)/\sin\theta\) for \(0<\theta<\pi\).
We write \(\mathcal N(a,v)\) for the real Gaussian law with mean \(a\)
and variance \(v\), and \(\operatorname{Cov}(U,V)=
\mathbb E[UV]-\mathbb E[U]\mathbb E[V]\) for covariance.

For \(1\le p<\infty\), the norm of a statistic \(F\) is
\[
 \|F\|_{L^p(M_{\mu,n})}
  =\left(\sum_\lambda M_{\mu,n}(\lambda)|F(\lambda)|^p\right)^{1/p}.
\]
Convergence in this norm uses the measure at the indicated rank.
For profile functions we write
\(\|f\|_\infty=\sup_{x\in\mathbb R}|f(x)|\).

\begin{theorem}[Skew moment laws and limit shape]
\label{thm:skew-shape-intro}
Fix a partition \(\mu\), put \(N=n+|\mu|\), and sample \(\lambda\)
from \(M_{\mu,n}\). For every \(j\ge0\) and finite \(p\ge1\),
\[
 N^{-j/2-1}C_j(\lambda/\mu)\longrightarrow\gamma_j
                 \quad\text{in }L^p(M_{\mu,n}).
\]
The genuine skew profile has the uniform limit shape
\[
 \sup_{s\in\mathbb R}
    |\overline\omega_{\lambda/\mu,n}(s)-\Omega(s)|
                         \longrightarrow0
                 \quad\text{in probability}.
\]
The full outer profile \(N^{-1/2}\lambda(\sqrt N\,s)\) has the same
uniform limit. The content moment limit is unchanged if \(N\) is
replaced by \(n\).
\end{theorem}

\begin{theorem}[Joint Gaussian fluctuations for a fixed skew base]
\label{thm:skew-asymptotics-intro}\label{thm:fixed-skew-asymptotics}
Fix a partition \(\mu\), put \(N=n+|\mu|\), and sample \(\lambda\)
from \(M_{\mu,n}\).
The following convergences hold jointly for each fixed finite family
of indices, in distribution and in all mixed moments.
\begin{enumerate}[label=\textup{(\roman*)},leftmargin=*]
\item For \(j\ge0\),
\[
 \frac{C_j(\lambda/\mu)-\gamma_jN^{j/2+1}}{N^{(j+1)/2}}
       \longrightarrow\frac{V_{j+2}}{(j+1)(j+2)},
\]
where \(V_d\) is defined in \textup{(ii)}.
\item For \(d\ge2\),
\[
 V^\mu_{d,N}:=\frac{q_d(\lambda)-a_dN^{d/2}}{N^{(d-1)/2}}
       \longrightarrow V_d=\sum_{r\ge2}b_{dr}\Xi_r,
\]
where the \(\Xi_r\) are independent and
\(\Xi_r\sim\mathcal N(\mathbf1_{\{r\text{ odd}\}},2r)\).
Only finitely many \(r\) have \(b_{dr}\ne0\) for a fixed \(d\).
Exact expectation centering gives the limits
\(\sum_{r\ge2}b_{dr}(\Xi_r-\mathbf1_{\{r\text{ odd}\}})\), and
\[
 \lim_{n\to\infty}
 \operatorname{Cov}_{M_{\mu,n}}(V^\mu_{d,N},V^\mu_{e,N})
                   =\sum_{r\ge2}2r\,b_{dr}b_{er}.
\]
\item For \(k\ge1\),
\[
 \int_{\mathbb R}u_k(s)\Delta^Y_{\lambda/\mu,n}(s)\,\mathrm ds
                     \longrightarrow\frac{\Xi_{k+1}}{k+1}.
\]
The limits at distinct indices are independent, with means
\(\mathbf1_{\{k+1\text{ odd}\}}/(k+1)\) and variances \(2/(k+1)\).
The zeroth skew mode is identically zero.
The corresponding Gaussian random distribution is specified in
Corollary~\ref{thm:polynomial-gaussian-field}.
\end{enumerate}
The content and corner limits are unchanged if \(N\) is replaced by
\(n\). The inner shape and all index families are fixed throughout.
\end{theorem}

Section~\ref{sec:limit-shape} proves these assertions from the exact
skew \(D\)-transfer, highest differences and polynomial operator moments.
Section~\ref{sec:clt} develops their quantitative mixed-moment and
inner-shape refinements.
All three classes of statistics are affine functions of a common finite
corner vector, so their joint limits use the same family \((\Xi_r)_{r\ge2}\).
Corollary~\ref{cor:fixed-mixed-moment-rate} gives an
\(O(n^{-1/2})\) error for each fixed mixed moment, including exact
expectation centering.
The convergence to the field is stated for fixed polynomial tests.
Subscripts on \(O(\cdot)\) indicate the fixed data on which an error
constant may depend.

\subsection{The first detectable inner-shape correction}

Stanley's exact relative-dimension formula
\cite[Theorem~3.1, equation~(12)]{StanleySkew} writes the density of the
fixed-prefix law, up to normalization, as a finite combination of
character ratios. Dousse--F\'eray studied the asymptotics of these
relative dimensions for fixed and growing inner diagrams
\cite[Theorem~1 and Corollary~3]{DousseFeray}.
Together with the joint Gaussian limits of normalized characters,
Stanley's expansion determines the leading covariance
\(\operatorname{Cov}(A_\mu,P(Y))\) in
Theorem~\ref{thm:inner-response-intro}, including the leading variance
correction \(-2m(m-1)\). The operator calculus gives uniform remainder
estimates for polynomial fluctuation moments at the same outer size and
the next variance term \(6\kappa_3(\mu)/\sqrt N\), with separate exact
centering. A class-algebra identity below explains the shape independence
of the leading variance correction.

Write \(\mathbb E_{\mu,N}\) and \(\mathbb E_{0,N}\) for expectation
under \(M_{\mu,N-m}\) and \(M_N=M_{\varnothing,N}\), respectively.
The associated variances are denoted by \(\operatorname{Var}_{\mu,N}\)
and \(\operatorname{Var}_{0,N}\). Under both laws use the same observable
\(Y_{d,N}=(q_d-a_dN^{d/2})/N^{(d-1)/2}\), with \(d\ge2\).
This is \(V^\mu_{d,N}\) under the fixed-prefix law, and its Gaussian
limit \(Y_d=\sum_{r\ge2}b_{dr}\Xi_r\) is \(V_d\) in
Theorem~\ref{thm:skew-asymptotics-intro}.
For a finite family of such indices, write \(Y_N=(Y_{d,N})_d\)
and \(Y=(Y_d)_d\). We call the first nonzero coefficient of a
moment difference, as the polynomial test varies, the first response.

\begin{theorem}[First-response classification at a common outer size]
\label{thm:inner-response-intro}
For every fixed \(\mu\) let \(\kappa_3(\mu)=C_2(\mu)-\binom m2\).
When \(m\ge2\), define the order parameter \(s\) and the Gaussian
polynomial \(A_\mu\) by
\[
 (s,A_\mu)=
 \begin{cases}
 (2,C_1(\mu)\Xi_2),&C_1(\mu)\ne0,\\
 (3,\kappa_3(\mu)\Xi_3),
       &C_1(\mu)=0,\ \kappa_3(\mu)\ne0,\\
 (4,C_3(\mu)\Xi_4-\frac{m(m-1)}4(\Xi_2^2-2)),
       &C_1(\mu)=0,\ \kappa_3(\mu)=0.
 \end{cases}
\]
For every fixed real polynomial \(P\) in finitely many variables
\(x_d\) indexed by corner degrees \(d\ge2\),
\[
 \mathbb E_{\mu,N}P(Y_N)-\mathbb E_{0,N}P(Y_N)
 =N^{-s/2}\operatorname{Cov}(A_\mu,P(Y))
                       +O_{\mu,P}(N^{-(s+1)/2}).
\]
The first coefficient is detected by \(P=x_3,x_4,x_3^2\) in the
three respective cases. If each \(Y_N\) is further centered by its
expectation under the corresponding law, the displayed leading
coefficient vanishes for every polynomial test in the first two cases.
In the third case, the coefficient for the new test vector is
\(\mathbb E[-m(m-1)(\Xi_2^2-4)P(Y-\mathbb EY)/4]\), and the
order-\(N^{-2}\) correction remains detectable. For \(\mu=\varnothing\)
or \((1)\), the laws at the same outer size agree exactly.

For every fixed \(\mu\), the variance of the outer content sum
\(C_1=C_1(\lambda)\) satisfies
\[
 \operatorname{Var}_{\mu,N}(C_1)-\operatorname{Var}_{0,N}(C_1)
 =-2m(m-1)+\frac{6\kappa_3(\mu)}{\sqrt N}+O_\mu(N^{-1}).
\]
Both variances use their respective exact expectations. The skew
sum \(C_1(\lambda/\mu)\) differs from \(C_1(\lambda)\) by a constant,
so either sum gives the same variance under \(M_{\mu,N-m}\).
\end{theorem}

Theorem~\ref{thm:adv-first-nonzero-inner},
Corollaries~\ref{cor:adv-centered-inner-correction}
and~\ref{cor:response-three-levels} prove the classification with uniform
coefficient-norm bounds on each fixed polynomial class.
Theorem~\ref{thm:response-universal-variance} gives the last assertion.
For \(Y_{3,N}=6C_1/N\), its two correction coefficients are
\(-72m(m-1)\) at order \(N^{-2}\) and \(216\kappa_3(\mu)\)
at order \(N^{-5/2}\), with remainder \(O_\mu(N^{-3})\).
The base \((3,1,1)\) gives the first coefficient \(-1440\).

\begin{remark}[Classical Gelfand limits]\label{thm:known-asymptotics-intro}
Setting \(\mu=\varnothing\) gives \(\lambda/\mu=\lambda\), \(N=n\)
and \(M_{\mu,n}=M_n\). Thus Theorems~\ref{thm:skew-shape-intro}
and~\ref{thm:skew-asymptotics-intro} recover the classical Gelfand limits of
M\'eliot \cite{Meliot,MeliotKerov}. In particular, for \(r\ge0\),
\[
 n^{-r-1}C_{2r}(\lambda)\longrightarrow\frac{(2r)!}{(r+1)!^2}
                         \quad\text{in }L^2(M_n).
\]
\end{remark}

\subsection{Method and organization}

Sections~\ref{sec:preliminaries}--\ref{sec:product} derive the skew transfer
and the pointwise difference bounds. The operator \(D\) acts on functions
of partitions. In its particle representation, each commutation with \(D\)
lowers the degree of a polynomial coefficient until the last commutator
is scalar. Sections~\ref{sec:content}--\ref{sec:hook-polynomials} extract
the exact factorial content and hook coefficients from finite operator
identities. Section~\ref{sec:limit-shape} proves the fixed-inner-shape
limit shape and Gaussian fluctuations from these difference-operator
identities; its empty-inner-shape specialization recovers the classical
Gelfand limits. Section~\ref{sec:clt} gives quantitative mixed-moment
bounds and the correction formulas at a common outer size.
Section~\ref{sec:conclusion} summarizes the results and discusses
extensions to general \(\beta\)-Plancherel measures.

\section{Boundary sequences and tableau enumeration}\label{sec:preliminaries}

\subsection{Corners and contents}

We use the content and corner conventions of
Definition~\ref{def:corner-observables-intro}. The addable and removable
corners alternate along the boundary, starting and ending with an addable
corner.
In particular,
\begin{equation}\label{eq:X-Y-cardinality}
  |X(\lambda)|=|Y(\lambda)|+1.
\end{equation}

\begin{example}\label{ex:corners-211}
Let \(\lambda=(2,1,1)\).  The contents of its boxes are
\[
\begingroup
\setlength{\unitlength}{6mm}
\begin{picture}(2,3)
\linethickness{0.4pt}
\put(0,0){\line(1,0){1}}
\put(0,1){\line(1,0){1}}
\put(0,2){\line(1,0){2}}
\put(0,3){\line(1,0){2}}
\put(0,0){\line(0,1){3}}
\put(1,0){\line(0,1){3}}
\put(2,2){\line(0,1){1}}
\put(0,2){\makebox(1,1){\(0\)}}
\put(1,2){\makebox(1,1){\(1\)}}
\put(0,1){\makebox(1,1){\(-1\)}}
\put(0,0){\makebox(1,1){\(-2\)}}
\end{picture}
\endgroup
\]
One checks that
\[
   X(\lambda)=\{-3,0,2\},
   \qquad
   Y(\lambda)=\{-2,1\}.
\]
Here \(|X(\lambda)|=3\) and \(|Y(\lambda)|=2\).
\end{example}

\subsection{The \texorpdfstring{\(0,1\)}{0,1} boundary sequence}
For a partition \(\lambda\), define its boundary sequence
\((w_i(\lambda))_{i\in\mathbb Z}\) by
\[
 w_i(\lambda)=
 \begin{cases}
  0,&i=\lambda_j-j\text{ for some }j\ge1,\\
  1,&\text{otherwise}.
 \end{cases}
\]
Thus \(w_i=0\) for \(i\ll0\) and \(w_i=1\) for \(i\gg0\).
The corner contents satisfy
\begin{equation}\label{eq:inner-outer-01-10}
 x\in X(\lambda)\Longleftrightarrow(w_{x-1},w_x)=(0,1),\qquad
 y\in Y(\lambda)\Longleftrightarrow(w_{y-1},w_y)=(1,0).
\end{equation}
We shall use the following classical local identity for signed corner
sums \cite[Lemma~6.3, (6.6)]{HanXiong}.
\begin{lemma}[Corner increments; \cite{HanXiong}]\label{lem:corner-increments}
For a polynomial \(\psi\), \(x\in X(\lambda)\), and \(y\in Y(\lambda)\),
\[
\begin{aligned}
 Q_\psi(\lambda^{x+})-Q_\psi(\lambda)
   &=\psi(x+1)-2\psi(x)+\psi(x-1),\\
 Q_\psi(\lambda^{y-})-Q_\psi(\lambda)
   &=-\psi(y+1)+2\psi(y)-\psi(y-1).
\end{aligned}
\]
\end{lemma}

\begin{lemma}[Classical tableau enumeration; \cite{FRT,Robinson,Schensted,StanleyEC2}]
\label{lem:classical-tableau-enumeration}\label{lem:tableau-classical-inputs}
For \(\lambda\vdash n\), the hook-length formula is
\(f^\lambda=n!/H_\lambda\) \cite{FRT}.
The Robinson--Schensted correspondence is a bijection
\[
 S_n\longleftrightarrow
 \coprod_{\lambda\vdash n}
 \{(P,Q):P,Q\text{ are standard tableaux of shape }\lambda\},
\]
under which inversion interchanges \(P\) and \(Q\)
\cite{Robinson,Schensted,StanleyEC2}.
Consequently involutions \(\sigma^2=\mathrm{id}\) correspond to
pairs \(P=Q\). Thus
\(\sum_{\lambda\vdash n} f^\lambda=I_n\) and
\(\sum_{\lambda\vdash n}(f^\lambda)^2=n!\).
\end{lemma}
Decomposition of an involution into fixed points and transpositions
gives \eqref{eq:involution-egf} and its stated recurrence.

\subsection{Tableau branching and involution coefficients}

Write \(\rho\gtrdot\lambda\) if \(\rho\) is obtained from \(\lambda\)
by adding one box, and \(\rho\lessdot\lambda\) if it is obtained by
deleting one box. Write \(I\) for the identity operator on the space
under consideration, and \([A,B]=AB-BA\) for the commutator. For functions on partitions set
\begin{equation}\label{eq:adjacency-operators}
 (\mathsf U g)(\lambda)=\sum_{\rho\gtrdot\lambda}g(\rho),\qquad
 (\mathsf L g)(\lambda)=\sum_{\rho\lessdot\lambda}g(\rho).
\end{equation}
These formulas act on functions. On partition states, their transpose
matrices sum over additions and deletions, respectively. Counting
common upper and lower covers of two partitions gives
\begin{equation}\label{eq:skew-AB}
 [\mathsf U,\mathsf L]=I,\qquad D=\mathsf U-\mathsf L-I.
\end{equation}
Equation~\eqref{eq:skew-AB} is the Young-lattice instance of Stanley's
differential-poset relation \cite[Theorem~2.2]{StanleyDifferential},
with the action on functions used here.

\begin{lemma}[Tableau branching; \cite{StanleyEC2,StanleyDifferential}]\label{lem:tableau-branching}
The deletion identity holds for nonempty \(\lambda\), and the addition
identity holds for every \(\lambda\).
\begin{equation}\label{eq:branching-down}
 f^\lambda=\sum_{\rho\lessdot\lambda}f^\rho,
\end{equation}
\begin{equation}\label{eq:branching-up}
 \sum_{\rho\gtrdot\lambda}f^\rho=(|\lambda|+1)f^\lambda.
\end{equation}
\end{lemma}

\begin{lemma}[Involution coefficients; \cite{MoserWyman}]\label{lem:involution-asymptotic}
The involution numbers defined by \eqref{eq:involution-egf} satisfy
the Moser--Wyman asymptotic \cite[(2.5)]{MoserWyman}
\[
 I_n\sim 2^{-1/2}e^{-1/4}(n/e)^{n/2}e^{\sqrt n}.
\]
For every fixed integer \(k\ge0\),
\begin{equation}\label{eq:In-ratio}
 \frac{I_{n-k}}{I_n}\sim n^{-k/2}\qquad(n\to\infty).
\end{equation}
\end{lemma}

\section{The difference operator}\label{sec:D-operator}

We write \(\delta_{ij}=\mathbf1_{\{i=j\}}\) for the Kronecker delta.

\begin{definition}[Matrix conventions]
\label{def:operator-conventions}\label{def:partition-matrix-element}
For a countable set \(\mathcal S\), the space
\(\mathbb C^{(\mathcal S)}=\bigoplus_{s\in\mathcal S}\mathbb Ce_s\)
consists of finite vectors, and \(\mathbb C^{\mathcal S}\) consists of
all functions \(\mathcal S\to\mathbb C\). A matrix is row-finite or
column-finite when each row
or column, respectively, has finite support. Its direct actions are
\begin{equation}\label{eq:function-action-same-matrix}
 (Tf)(s)=\sum_tT_{st}f(t),
\end{equation}
\begin{equation}\label{eq:ket-matrix-coefficients}
 Te_t=\sum_sT_{st}e_s.
\end{equation}
The first requires row finiteness and the second column finiteness.
For partitions, \(\ket\lambda=e_\lambda\) denotes the basis vector
and \(\langle\mu|\) selects the \(e_\mu\)-coordinate. The row \(\mu\),
column \(\lambda\) entry is
\begin{equation}\label{eq:partition-matrix-element}
 \langle\mu|T|\lambda\rangle=T_{\mu\lambda}.
\end{equation}
The algebraic dual action is the transpose action
\begin{equation}\label{eq:function-action-transpose}
 (T^{\mathsf T}f)(s)=\sum_tT_{ts}f(t).
\end{equation}
Products act from right to left. We use
\(\operatorname{ad}_D^0(T)=T\) and
\(\operatorname{ad}_D^{r+1}(T)=[D,\operatorname{ad}_D^r(T)]\).
The adjoint action is locally nilpotent on a family if each member is
killed by some positive iterate; the iterate may depend on that member.
\end{definition}

\begin{lemma}[Row-finite algebra]\label{lem:row-finite-algebra}
Row-finite matrices are closed under products and commutators, acting
on all functions by ordinary matrix multiplication. Column-finite
matrices have the corresponding property on finite vectors. In
particular, \([D,AB]=[D,A]B+A[D,B]\).
\end{lemma}
\begin{proof}
The support of a row of \(AB\) lies in the union of finitely many
finite row supports of \(B\). This justifies regrouping all sums.
Transpose for column finiteness; expand the commutator for the last assertion.
\end{proof}

Multiplication by \(g\) sends \(f(\lambda)\) to \(g(\lambda)f(\lambda)\),
so its matrix is diagonal. The difference operator \(D\) is row-finite
and satisfies \(D\mathbf1=0\). Thus the \(j\)th difference of \(g\)
is obtained by applying the \(j\)th commutator of its multiplication
operator to \(\mathbf1\).

\subsection{Skew Gelfand transfer}\label{sec:skew-transfer}
For a function \(g\) on partitions, recall that
\[
 S_{\mu,g}(n)=\sum_{\substack{\lambda\supseteq\mu\\|\lambda/\mu|=n}}
                 f^{\lambda/\mu}g(\lambda).
\]
Put \(I_t=0\) for \(t<0\). The transfer formula uses \(g\) on the whole
lattice: deletion can take a diagram outside
\(\{\lambda:\lambda\supseteq\mu\}\).

The unweighted enumeration is due to Sagan--Stanley
\cite[Corollaries~3.5 and~6.4]{SaganStanleySkew} and Stanley
\cite[Theorem~2.2]{StanleySkew}. McKay--Morse--Wilf
\cite[Theorem~1]{McKayMorseWilf} gave the corresponding asymptotics
for a prescribed tableau prefix. The following weighted identity
follows from the up--down relation on Young's lattice
\cite{StanleyDifferential}.

\begin{lemma}[Skew transfer; differential-poset identity]\label{lem:skew-transfer}
Suppose \(D^{R+1}g=0\) on every partition, with \(R\) a nonnegative
integer.
Fix \(\mu\). For each \(\eta\subseteq\mu\) put \(t_\eta=|\mu/\eta|\).
For integers \(0\le k\le R\), define
\begin{equation}\label{eq:skew-transfer-polynomial}
 P_{\mu,g,k}(x)=
 \sum_{\eta\subseteq\mu}\frac{f^{\mu/\eta}x^{t_\eta}}{t_\eta!}
 \sum_{j=0}^{R-k}\frac{D^{j+k}g(\eta)}{j!}x^j .
\end{equation}
For \(k=0\), this is the polynomial \(P_{\mu,g}\) in
Theorem~\ref{thm:skew-exponential-polynomiality}. Its degree is at most \(|\mu|+R-k\), and
\begin{equation}\label{eq:skew-transfer-egf}
 \sum_{n\ge0}S_{\mu,D^kg}(n)\frac{x^n}{n!}
 =I(x)P_{\mu,g,k}(x).
\end{equation}
Equivalently, for \(n\ge0\),
\begin{equation}\label{eq:skew-transfer-coeff}
 S_{\mu,D^kg}(n)
 =\sum_{\eta\subseteq\mu}
   \sum_{\substack{0\le j\le R-k\\t_\eta+j\le n}}
   \frac{n!\,f^{\mu/\eta}D^{j+k}g(\eta) I_{n-t_\eta-j}}
        {t_\eta!\,j!\,(n-t_\eta-j)!}.
\end{equation}
For integers \(k>R\), set \(P_{\mu,g,k}=0\). Both weighted sums
vanish because \(D^kg=0\).
\end{lemma}

For \(g=F_{\nu,\mu}\), Theorem~\ref{thm:degree-intro} gives the exact
truncation order. The exponential factor \(I(x)\) reflects the choice
of the local operator \(D\). In the related skew transfer of Han--Xiong
\cite[Theorems~1.3 and~1.4]{HanXiong}, the operator only sums over box
additions and the observable contains a hook-product denominator.

\begin{proof}[Proof of Lemma~\ref{lem:skew-transfer}]

Put \(a=|\mu|\). By \eqref{eq:skew-AB},
\([D,\mathsf L]=I\).

In particular,
\[
 D \mathsf L^r=\mathsf L^rD+r\mathsf L^{r-1}\quad(r\ge1),
 \qquad e^{x\mathsf L}D=(D-xI)e^{x\mathsf L}.
\]
All exponentials in this proof are formal power series. Their
coefficients are well defined because every vertex of Young's lattice
has finitely many neighbours. Put
\[
 F(x)=e^{x+x^2/2}e^{x\mathsf L}e^{xD}g.
\]
Differentiating and using the preceding commutation identity yields
\[
 F'(x)=(1+x)F(x)+\mathsf L F(x)+(D-xI)F(x)
       =(I+\mathsf L+D)F(x)=\mathsf U F(x),
 \qquad F(0)=g.
\]
The coefficient recursion for this formal differential equation is
unique, so
\begin{equation}\label{eq:skew-factorization}
 e^{x\mathsf U}g=e^{x+x^2/2}e^{x\mathsf L}e^{xD}g.
\end{equation}
A chain is saturated when successive shapes differ by one box.
Counting such chains gives, for any function \(h\),
\[
 (\mathsf U^ng)(\mu)=S_{\mu,g}(n),
 \qquad
 (\mathsf L^rh)(\mu)=
  \sum_{\substack{\eta\subseteq\mu\\|\mu/\eta|=r}}
       f^{\mu/\eta}h(\eta).
\]
In particular, \((\mathsf L^rh)(\mu)=0\) for \(r>a\). Apply
\eqref{eq:skew-factorization} to \(D^kg\) and evaluate at \(\mu\).
Since \(D^{R+1}g=0\), the \(D\)-exponential truncates after degree \(R-k\);
the \(\mathsf L\)-exponential at \(\mu\) truncates after degree \(a\). This gives
\eqref{eq:skew-transfer-polynomial}--\eqref{eq:skew-transfer-egf} and the
degree bound. Coefficient extraction using \eqref{eq:involution-egf}
gives \eqref{eq:skew-transfer-coeff}.
\end{proof}

\begin{remark}[Exponential polynomiality]\label{rem:skew-polynomiality}
The polynomial in Lemma~\ref{lem:skew-transfer} is the exponential
transform \(e^{-x-x^2/2}\sum_n S_{\mu,D^kg}(n)x^n/n!\).
In contrast with the Plancherel-weighted skew identities for the
Han--Xiong operator \cite[Theorems~1.3--1.4]{HanXiong}, the untransformed sum need not be a polynomial in
\(n\). Already \(\mu=\varnothing\) and \(g=1\) give
\(S_{\varnothing,1}(n)=I_n\), with \(I_{2r}\ge(2r)!/(2^r r!)\), so this sequence is not a polynomial in \(n\). The difference comes from using the present
operator \(\mathsf U-\mathsf L-I\) and the weights \(f^{\lambda/\mu}\).
\end{remark}

\begin{corollary}[The skew normalizing constant; Sagan--Stanley and Stanley]\label{cor:skew-normalizer}
The skew-tableau enumeration formula of Sagan--Stanley
\cite[Corollary~3.5]{SaganStanleySkew} and Stanley
\cite[Theorem~2.2]{StanleySkew} is
\begin{equation}\label{eq:skew-normalizer}
 Z_\mu(n):=\sum_{|\lambda/\mu|=n}f^{\lambda/\mu}
 =\sum_{\substack{\eta\subseteq\mu\\|\mu/\eta|\le n}}
   \binom{n}{|\mu/\eta|}f^{\mu/\eta}I_{n-|\mu/\eta|}.
\end{equation}
Thus \(M_{\mu,n}(\lambda)=f^{\lambda/\mu}/Z_\mu(n)\) defines a probability
measure on partitions containing \(\mu\) with \(|\lambda/\mu|=n\).
For every \(D\)-polynomial \(g\), its expectation under \(M_{\mu,n}\) is given
exactly by \eqref{eq:skew-transfer-coeff} with \(k=0\), divided by
\eqref{eq:skew-normalizer}.
\end{corollary}
\begin{proof}
Apply Lemma~\ref{lem:skew-transfer} to the constant function \(\mathbf1\),
for which \(D\mathbf1=0\).
Every term in the defining sum for \(Z_\mu(n)\) is nonnegative, and there
exists a saturated chain of length \(n\) starting at \(\mu\), so \(Z_\mu(n)>0\).
\end{proof}

\begin{example}[Two fixed inner shapes]\label{ex:small-skew-normalizers}
For \(\mu=(1)\), every standard tableau starts with its unique box;
hence \(f^{\lambda/(1)}=f^\lambda\) and
\[
 Z_{(1)}(n)=I_{n+1}=I_n+nI_{n-1},\qquad
 \sum_{n\ge0}Z_{(1)}(n)\frac{x^n}{n!}=I(x)(1+x).
\]
For \(\mu=(2)\), each relevant tableau count for its three
subpartitions is one. Thus
\[
 Z_{(2)}(n)=I_n+nI_{n-1}+\binom n2I_{n-2},\qquad
 \sum_{n\ge0}Z_{(2)}(n)\frac{x^n}{n!}
     =I(x)(1+x+x^2/2).
\]
At \(n=2\) the outer shapes are
\((4),(3,1),(2,2),(2,1,1)\), with skew counts \(1,2,1,1\),
whose sum is \(Z_{(2)}(2)=5\).
\end{example}

\begin{lemma}[A fixed tableau prefix; \cite{McKayMorseWilf}]\label{lem:fixed-tableau-prefix}
Let \(\mu\) have size \(m\), let \(T\) be a fixed standard tableau of
shape \(\mu\), and choose uniformly among all standard tableaux on
\(N=n+m\) boxes, over all outer shapes. For a fixed outer shape \(\lambda\),
the number of completions of the prescribed entries \(1,\ldots,m\) in
\(T\) is \(f^{\lambda/\mu}\). Thus the event has probability \(Z_\mu(n)/I_N\),
and the conditional law of \(\lambda\) is \(M_{\mu,n}\). As \(n\to\infty\)
with \(\mu\) fixed,
\[
 \frac{Z_\mu(n)}{I_N}\longrightarrow\frac{f^\mu}{m!}
 \qquad\text{by~\cite[Theorem~1]{McKayMorseWilf}.}
\]
The initial-shape event is the disjoint union over its \(f^\mu\)
possible tableaux \(T\). Its probability tends to \((f^\mu)^2/m!\),
and it gives the same conditional law.
\end{lemma}

\begin{example}[A specified prefix and its shape]
\label{ex:prefix-shape-distinction}
For \(\mu=(2,1)\), the transfer polynomial is
\(1+2x+x^2+x^3/3\), and
\[
 Z_{(2,1)}(n)=I_n+2nI_{n-1}+n(n-1)I_{n-2}
                     +2\binom n3I_{n-3}.
\]
Here \(f^{(2,1)}=2\). Each specified prefix has limiting probability
\(2/3!=1/3\), while their union, specifying only the shape, has
limiting probability \(2/3\). Both conditionings give \(M_{(2,1),n}\).
\end{example}

The same up--down operator occurs in the hook-weight branching identity
of Carde, Loubert, Potechin and Sanborn \cite[Lemma~1.2]{CLPS}:
for \(W_\vartheta(\lambda)=\prod_{\square\in\lambda}
(1+\vartheta^{h_\lambda(\square)})/(1-\vartheta^{h_\lambda(\square)})\),
one has \(DW_\vartheta=\frac{2\vartheta}{1-\vartheta}W_\vartheta\).
We study its action on polynomial corner observables and their exact
highest differences.
The global identity
\(\sum_\lambda x^{|\lambda|}W_\vartheta(\lambda)/H_\lambda
=\exp((1+\vartheta)x/(1-\vartheta)+x^2/2)\)
is the exponential specialization of Kawanaka's Schur-function identity
\cite[Theorem~1.1]{Kawanaka}. These identities may be read formally in
the hook parameter \(\vartheta\), with \(W_\vartheta(\varnothing)=1\).

\subsection{Leading expectations for a fixed inner partition}

\begin{corollary}[Skew leading coefficient]\label{cor:expectation-asymptotic}
Fix a partition \(\mu\). Suppose \(D^{d+1}g=0\) on all partitions,
where \(d\) is a nonnegative integer. As \(n\to\infty\),
\[
 \E_{M_{\mu,n}}g=O(n^{d/2}).
\]
If \(D^dg(\varnothing)\ne0\), then
\begin{equation}\label{eq:asymp-from-D}
 \E_{M_{\mu,n}}g\sim
       \frac{D^dg(\varnothing)}{d!}n^{d/2}.
\end{equation}
\end{corollary}
\begin{proof}
Put \(m=|\mu|\) and \(c_\mu=f^\mu/m!>0\).
By \eqref{eq:skew-normalizer} and \eqref{eq:In-ratio}, the unique
highest-order normalizing term is \(\eta=\varnothing\), and
\[
 Z_\mu(n)/I_n\sim c_\mu n^{m/2}.
\]
In \eqref{eq:skew-transfer-coeff} with \(k=0\), the term indexed by
\((\eta,j)\), after division by \(I_n\), is
\(O(n^{(m-|\eta|+j)/2})\). There are finitely many terms. The unique
candidate of order \(n^{(m+d)/2}\) is
\((\eta,j)=(\varnothing,d)\), with coefficient
\(c_\mu D^dg(\varnothing)/d!\). Division by the normalizer proves
the bound and, when this coefficient is nonzero, the asymptotic formula.
\end{proof}

\subsection{The empty-base transfer as an operator identity}

Setting \(\mu=\varnothing\) in the skew transfer gives the coefficient
functional used in the Gaussian calculation.

\begin{proposition}[Exponential transfer]\label{prop:transfer}
For \(n\ge0\), let
\[
   A_g(n)=\sum_{|\lambda|=n} f^\lambda g(\lambda),
\]
and set \(A_g(-1)=0\).
If \(B_g(n)=\sum_{|\lambda|=n}f^\lambda Dg(\lambda)\), then, for \(n\ge0\),
\begin{equation}\label{eq:transfer-recurrence}
   A_g(n+1)-A_g(n)-nA_g(n-1)=B_g(n).
\end{equation}
Consequently, if \(D^r g=0\) for an integer \(r\ge1\), then
\begin{equation}\label{eq:D-polynomial-egf}
  \sum_{n\ge0}A_g(n)\frac{x^n}{n!}
  =\left(\sum_{j=0}^{r-1}\frac{D^j g(\varnothing)}{j!}x^j\right)
    \exp\left(x+\frac{x^2}{2}\right).
\end{equation}
\end{proposition}
\begin{proof}
Setting \(\mu=\varnothing\) and \(k=0\) in
Lemma~\ref{lem:skew-transfer} gives \eqref{eq:D-polynomial-egf}.
For the recurrence, the one-box corner identity gives
\(D=\mathsf U-\mathsf L-I\). The tableau branching identities therefore
give \(B_g(n)=A_g(n+1)-A_g(n)-nA_g(n-1)\), with
\(A_g(-1)=0\). This proves \eqref{eq:transfer-recurrence}.
\end{proof}

\section{Finite-band commutator calculus}\label{sec:commutator-calculus}

We use the classical infinite-wedge representation and its normal-ordering
central extension (see \cite[Chapter~14]{Kac},
\cite[Appendix~A]{OkounkovWedge}, and
\cite[Sections~2.1--2.2]{OkounkovPandharipande}).
The classical commutation relations are stated below with references;
the subsequent calculations identify our boundary operators and their
degree bounds.

Let \(v_i=\mathbf1_{\{i\ge0\}}\) be the boundary sequence of the
empty partition, called the vacuum. Let \(\calM\) be the set of
\(0,1\) sequences \(w\) on \(\Z\) for which \(w-v\) has finite support,
where \(\supp(u)=\{i\in\Z:u_i\ne0\}\).
A site \(i\) is occupied when \(w_i=0\) and empty when \(w_i=1\).
The charge
\[
  \operatorname{ch}(w)=\sum_{i\in\Z}(v_i-w_i)
\]
counts the net excess of occupied sites relative to the vacuum.
If \(i_1>i_2>\cdots\) are the occupied sites of \(w\), then
\(i_j=-j+\operatorname{ch}(w)\) for large \(j\), by counting
occupied sites above a common index below all changes. Thus
\(\operatorname{ch}(w)=0\) exactly for the boundary sequences
of partitions, recovered by \(\lambda_j=i_j+j\).
We first work on \(\calM\) and then restrict to this charge-zero sector.

For \(x\in\Z\), set
\[
  \epsilon_x(w)=w_x-w_{x-1}.
\]
Then \(\epsilon_x(w)=1\) precisely when \((w_{x-1},w_x)=(0,1)\), and
\(\epsilon_x(w)=-1\) precisely when \((w_{x-1},w_x)=(1,0)\).
On the charge-zero sector these are, respectively, addable and removable
corner contents; on the other sectors they specify admissible particle
moves. Let \(s_xw\)
denote the sequence obtained from \(w\) by interchanging \(w_{x-1}\) and
\(w_x\). Thus \(s_xw\) is obtained by moving a particle one step to the right
if \(\epsilon_x(w)=1\), and by moving a particle one step to the left if
\(\epsilon_x(w)=-1\).

Apply Definition~\ref{def:operator-conventions} with \(\mathcal S=\calM\).
We write \(V=\bigoplus_{w\in\calM}\C\ket w\) for the space spanned by
boundary states and \(\calF=\C^{\calM}\) for its algebraic dual, the
space of all functions on \(\calM\). With the direct matrix action
\eqref{eq:function-action-same-matrix}, multiplication by \(h\) is
diagonal on both spaces.
The difference operator on functions is
\[
  (Df)(w)=\sum_{x\in\Z}\epsilon_x(w)\bigl(f(s_xw)-f(w)\bigr).
\]
Only finitely many corners contribute. Each swap preserves charge,
so \(D\) and diagonal multiplication restrict to each charge sector.
The two eventual values of \(w_x\) are \(0\) on the left and \(1\)
on the right. Telescoping the jumps gives
\(\sum_x\epsilon_x(w)=1\), hence
\[
  D=\mathsf J-I,
  \qquad
  (\mathsf Jf)(w)=\sum_{x\in\Z}\epsilon_x(w)f(s_xw),
\]
and therefore
\[
  \ad_D=\ad_{\mathsf J}.
\]

\subsection{Signed corner sums as diagonal particle operators}

In boundary coordinates the signed corner sum of
Definition~\ref{def:corner-observables-intro} is
\(Q_\psi(w)=\sum_x(w_x-w_{x-1})\psi(x)\).

\begin{lemma}[Discrete summation by parts]\label{lem:summation-by-parts}
Let
\[
  h_\psi(i)=\psi(i+1)-\psi(i).
\]
Then
\[
  Q_\psi(w)=\psi(0)+\sum_{i\in\Z}(v_i-w_i)h_\psi(i).
\]
The sum on the right is finite.
\end{lemma}

\begin{proof}
At the vacuum, \(v_x-v_{x-1}=\delta_{x0}\), so \(Q_\psi(v)=\psi(0)\).
Subtract the vacuum contribution and set \(u=w-v\). Since \(u\)
has finite support, shifting indices gives
\[
 Q_\psi(w)-\psi(0)=\sum_x(u_x-u_{x-1})\psi(x)
                  =\sum_i u_i\bigl(\psi(i)-\psi(i+1)\bigr).
\]
\end{proof}

Write \(\mathsf M_\psi\) for multiplication by \(Q_\psi\). A particle move
changes the occupation numbers in the lemma while keeping the site
weights \(h_\psi(i)\) fixed.

\subsection{A closed algebra of signed particle-move operators}

The signed particle-move operators defined below are closed under
commutators up to a scalar term.

On \(V\), define particle deletion and insertion as follows. For a
sequence \(w\), put
\[
  \rho_i(w)=\#\{j>i: w_j=0\}.
\]
This number is finite because \(w_j=1\) for all sufficiently large \(j\).

Define an operator \(a_i\), which deletes a particle at \(i\) according to
\[
  a_i\ket w=
  \begin{cases}
    (-1)^{\rho_i(w)}\ket{w^{i:0\to1}}, & w_i=0,\\
    0, & w_i=1.
  \end{cases}
\]
Here \(w^{i:0\to1}\) denotes the sequence obtained by changing
\(w_i=0\) into \(w_i=1\). Define \(w^{i:1\to0}\) by the reverse change.
The operator \(a_i^*\) inserts a particle at \(i\) according to
\[
  a_i^*\ket w=
  \begin{cases}
    (-1)^{\rho_i(w)}\ket{w^{i:1\to0}}, & w_i=1,\\
    0, & w_i=0.
  \end{cases}
\]
Deletion lowers the charge by one and insertion raises it by one.
Their signs give the canonical anticommutation relations (CAR) below.

\begin{lemma}[Canonical anticommutation relations; \cite{Kac,OkounkovWedge}]\label{lem:CAR}
For \(i,j\in\mathbb Z\),
\[
 a_ia_j+a_ja_i=0,\qquad
 a_i^*a_j^*+a_j^*a_i^*=0,\qquad
 a_ia_j^*+a_j^*a_i=\delta_{ij}I.
\]
\end{lemma}

For \(r,s\in\Z\), define the charge-preserving operator
\[
  F_{rs}=a_r^*a_s.
\]
If \(r\ne s\), then \(F_{rs}\) moves a particle from \(s\) to \(r\), with the
CAR sign, provided there is a particle at \(s\) and a hole at \(r\); otherwise it
gives zero. In the adjacent moves used below, the two CAR signs multiply to
\(+1\). If \(r=s\), then
\[
  F_{rr}\ket w=(1-w_r)\ket w,
\]
the particle-indicator at \(r\). Let
\[
  n_r^0=1-v_r
\]
be the particle-indicator of the vacuum at \(r\). Subtracting the
vacuum value on the diagonal gives the normal-ordered operators
\[
  E_{rs}=\begin{cases}
    F_{rs}, & r\ne s,\\
    F_{rr}-n_r^0I, & r=s.
  \end{cases}
\]

\begin{lemma}[Normal-ordered matrix-unit relations; \cite{Kac,OkounkovPandharipande}]\label{lem:matrix-unit}
For all \(r,s,t,u\in\Z\),
\[
  [E_{rs},E_{tu}]
  =\delta_{st}E_{ru}-\delta_{ru}E_{ts}
  +\delta_{st}\delta_{ru}(n_r^0-n_s^0)I.
\]
In particular, modulo scalar operators, the \(E_{rs}\)'s satisfy the ordinary
matrix-unit commutation relations.
\end{lemma}

\subsection{Finite-band polynomial particle operators}\label{sec:finite-band-operators}

On the one-particle space with basis \((e_i)_{i\in\Z}\), let
\(\mathsf N\) denote position and \(S\) the right shift,
\[
  \mathsf{N}e_i=ie_i,
  \qquad
  Se_i=e_{i+1}.
\]
Set \(\nabla=S-S^{-1}\) and \(T_0=S+S^{-1}\).
A one-particle operator \(T=(t_{rs})_{r,s\in\Z}\) is called finite-band if
there exists \(R_T\ge0\) such that \(t_{rs}=0\) whenever \(|r-s|>R_T\). For a
finite-band operator \(T\), define its lifted operator by
\[
  \Theta(T)=\sum_{r,s\in\Z} t_{rs}E_{rs}.
\]

The lift acts by a finite sum on each \(\ket w\). Indeed, put
\(K=\supp(w-v)\). A nonzero off-diagonal move has a particle at
\(s\), a hole at \(r\), and \(|r-s|\le R_T\). If both endpoints
avoid \(K\), then \(s<0\le r\); otherwise both lie within
distance \(R_T\) of \(K\). Thus all endpoints lie in the finite
\(R_T\)-neighbourhood of \(K\cup\{0\}\). Diagonal terms are
\(E_{ii}\ket w=(v_i-w_i)\ket w\) and vanish outside \(K\). The same argument with source and target
interchanged proves row finiteness. Thus every finite-band lift acts
both on \(V\) and on \(\calF\); products are justified by
Lemma~\ref{lem:row-finite-algebra}.

For a nonnegative integer \(p\), let \(\calD_p\) be the space of
operators admitting a representation
\[
  T=\sum_{|k|\le R}\phi_k(\mathsf{N})S^k,
  \qquad
  \phi_k(z)\in\C[z],\quad \deg \phi_k\le p,
\]
for some finite nonnegative integer \(R=R(T)\). This bandwidth bound may depend
on \(T\), but is independent of the source index \(s\); no common bandwidth
bound is imposed on \(\calD_p\).
In \(\phi_k(\mathsf N)S^k\) the shift acts first. With \(s\) the source
and \(r=s+k\) the target, the coefficient is \(\phi_k(s+k)\). Thus, with
\(\phi_k=0\) for \(|k|>R\),
\[
  t_{rs}=\phi_{r-s}(r),\qquad
  \Theta(T)=\sum_{|k|\le R}\sum_{s\in\Z}\phi_k(s+k)E_{s+k,s}.
\]
The terms with \(k\ne0\) move a particle from \(s\) to \(s+k\);
the term with \(k=0\) multiplies \(\ket w\) by
\(\sum_s\phi_0(s)(v_s-w_s)\).

Set \(\calD_{-1}=\{0\}\). Define
\[
  \calA_p=\Theta(\calD_p)+\C I\qquad (p\ge0),
\]
and also
\[
  \calA_{-1}=\C I,
  \qquad
  \calA_{-2}=\{0\}.
\]

The lifted commutator equals the lift of the one-particle commutator
up to a scalar. This is the finite-band form of the classical
projective representation of matrices supported on a bounded strip
around the diagonal \cite[Section~1, Lemma~1.3]{BlochOkounkov}.
The scalar for \(\nabla\) comes from crossing the vacuum boundary
between sites \(-1\) and \(0\).

\begin{lemma}[Lifted finite-band commutators; \cite{BlochOkounkov}]\label{lem:lifted-commutator}
Let \(U=(u_{ij})\) and \(T=(t_{ij})\) be finite-band one-particle operators. Then
\begin{equation}\label{eq:finite-band-commutator}
 [\Theta(U),\Theta(T)]=\Theta([U,T])+c(U,T)I,\qquad
 c(U,T)=\sum_{i,j\in\Z}u_{ij}t_{ji}(n_i^0-n_j^0).
\end{equation}
The scalar sum is finite. In particular, for \(\nabla=S-S^{-1}\) and
\(T=\sum_k f_k(\mathsf N)S^k\) one has
\begin{equation}\label{eq:boundary-cocycle}
 c(\nabla,T)=-f_{-1}(-1)-f_1(0).
\end{equation}
\end{lemma}

Repeated commutation lowers the coefficient degree.

\begin{lemma}[Closed filtered commutator lemma]\label{lem:closed-filtration}
For every nonnegative integer \(p\), one has
\[
  [\Theta(\nabla),\calA_p]\subseteq\calA_{p-1}.
\]
Consequently, if \(X\in\calA_p\), then
\[
  \ad_{\Theta(\nabla)}^{p+2}(X)=0.
\]
\end{lemma}

\begin{proof}
It suffices to commute \(\nabla\) with a single band. Using
\(S^{\pm1}f(\mathsf N)=f(\mathsf N\mp1)S^{\pm1}\), we obtain
\begin{equation}\label{eq:one-particle-difference}
\begin{aligned}
 [\nabla,f(\mathsf N)S^k]
   ={}&-\bigl(f(\mathsf N+1)-f(\mathsf N)\bigr)S^{k-1}\\
     &-\bigl(f(\mathsf N)-f(\mathsf N-1)\bigr)S^{k+1}.
\end{aligned}
\end{equation}
Both differences lower the degree by one and vanish on constants.
Thus \([\nabla,\calD_p]\subseteq\calD_{p-1}\). By
Lemma~\ref{lem:lifted-commutator} lifting adds at most a scalar, which
belongs to \(\calA_{p-1}\) even for \(p=0\). After \(p+1\)
commutators only a scalar can remain, and the next commutator kills it.
\end{proof}

\subsection{Adjoint local nilpotence of signed corner multiplication}

\begin{lemma}[Lifts on partition functions]\label{lem:partition-lifts}
Let \(h_\psi(z)=\psi(z+1)-\psi(z)\). On the
function space with direct matrix action,
\[
 D=-\Theta(\nabla)-I,\qquad \mathsf M_\psi=\psi(0)I+\Theta(h_\psi(\mathsf N)).
\]
\end{lemma}
\begin{proof}
Put \(\mathsf J_{\mathrm{ket}}=\Theta(S-S^{-1})\). An adjacent particle move has sign \(+1\),
so \(\mathsf J_{\mathrm{ket}}\ket w=\sum_x\epsilon_x(w)\ket{s_xw}\).
The reverse move has the opposite corner sign, hence \(\mathsf J_{\mathrm{ket}}^{\mathsf T}=-\mathsf J_{\mathrm{ket}}\).
The signed sum \(\mathsf J\) in \(D=\mathsf J-I\) acts on functions by
\(\mathsf J=\mathsf J_{\mathrm{ket}}^{\mathsf T}=-\Theta(\nabla)\). This proves the first identity with the
matrix conventions of Definition~\ref{def:partition-matrix-element}.

For the second, \(E_{ii}\ket w=(v_i-w_i)\ket w\).
The diagonal entry of \(\Theta(h_\psi(\mathsf N))\) is therefore
\(Q_\psi(w)-\psi(0)\) by Lemma~\ref{lem:summation-by-parts}.
All these lifts are column-finite by the support argument above. They
are also row-finite, since reversing a bounded move from a given output
allows only finitely many inputs. Lemma~\ref{lem:row-finite-algebra}
then identifies compositions with matrix products on both spaces.
\end{proof}

\begin{proposition}[Adjoint local nilpotence]\label{prop:local-nilpotence}
Let \(\psi(z)\in\C[z]\), and let \(\mathsf M_\psi\) be multiplication by
\(Q_\psi(\lambda)\). If \(\psi\ne0\) and \(d=\deg\psi\), then
\[
  \ad_D^{d+1}(\mathsf M_\psi)=0.
\]
If \(\psi=0\), the same conclusion is trivial, with any positive exponent.
\end{proposition}

\begin{proof}
For constant \(\psi\), Lemma~\ref{lem:partition-lifts} gives
\(\mathsf M_\psi=\psi(0)I\). For \(d\ge1\), the polynomial \(h_\psi\)
has degree \(d-1\), so \(\mathsf M_\psi\in\calA_{d-1}\).
Commutation with \(D\) is commutation with \(-\Theta(\nabla)\).
Lemma~\ref{lem:closed-filtration} now gives the assertion.
\end{proof}

\begin{example}[Low-degree corner powers]\label{rem:higher-power-examples}
For \(\psi(z)=z^2\), \(Q_\psi=2|\lambda|\). Let \(\mathsf N_0\) denote
multiplication by \(|\lambda|\) on partition functions. Then
\([\mathsf U,\mathsf N_0]=\mathsf U\) and
\([\mathsf L,\mathsf N_0]=-\mathsf L\). The relations in \eqref{eq:skew-AB}
give
\[
 [D,2\mathsf N_0]=2(\mathsf U+\mathsf L),\qquad
 \ad_D^2(2\mathsf N_0)=4I,\qquad \ad_D^3(2\mathsf N_0)=0.
\]
For general \(d\ge1\), adding a box of content \(c\) changes \(q_d\)
by \((c+1)^d-2c^d+(c-1)^d\). Summation from the empty partition gives
\[
\begin{aligned}
 q_3&=6C_1,& q_4&=12C_2+2C_0,\\
 q_5&=20C_3+10C_1,& q_6&=30C_4+30C_2+2C_0 .
\end{aligned}
\]
The corresponding top commutators are evaluated in
Lemma~\ref{lem:highest-central-commutator}.
\end{example}

\section{Pointwise polynomiality and highest differences}\label{sec:product}

The signed-corner product bound in Theorem~\ref{thm:degree-intro}
follows from the classical iterated Leibniz rule for inner derivations
\cite[Chapter~I, Section~2, (12)]{Jacobson}.
We record its operator form, which preserves the order of the
multiplication factors.

\begin{lemma}[Iterated Leibniz identity; \cite{Jacobson}]\label{lem:noncommutative-multinomial}
Let \(D,M_1,\ldots,M_\ell\) be linear operators on a common invariant
domain. For every nonnegative integer \(r\), with the sum over
nonnegative integers \(j_0,\ldots,j_\ell\),
\[
  D^rM_1\cdots M_\ell
  =\sum_{j_1+\cdots+j_\ell+j_0=r}
    \frac{r!}{j_1!\cdots j_\ell!j_0!}
    \bigl(\ad_D^{j_1}M_1\bigr)
    \cdots
    \bigl(\ad_D^{j_\ell}M_\ell\bigr)
    D^{j_0}.
\]
\end{lemma}

\subsection{The highest central term}

\begin{lemma}[Highest commutator]\label{lem:highest-central-commutator}
Let \(\psi\) have degree \(p\ge1\) and leading coefficient \(a\), and let \(\mathsf M_\psi\)
be multiplication by \(Q_\psi\). With the convention that
\(\binom{p}{p/2}=0\) for odd \(p\),
\begin{equation}\label{eq:highest-central-commutator}
 \operatorname{ad}_D^p(\mathsf M_\psi)
       =a\,p!\binom{p}{p/2}\,I,
 \qquad \operatorname{ad}_D^{p+1}(\mathsf M_\psi)=0.
\end{equation}
In particular, the \(p\)th commutator itself vanishes when \(p\) is odd.
\end{lemma}
\begin{proof}
Let \(\mathcal L_\nabla(T)=-[\nabla,T]\) on one-particle matrices.
By Lemma~\ref{lem:partition-lifts} the one-particle polynomial is
\(h(z)=\psi(z+1)-\psi(z)\), of degree \(p-1\) with leading coefficient \(ap\).
In \eqref{eq:one-particle-difference} each finite difference lowers
 the degree by one and both have the same leading coefficient. After
\(p-1\) applications no lower-degree term survives, and the two
choices of shift at each step give exactly
\[
 \mathcal L_\nabla^{p-1}\bigl(h(\mathsf N)\bigr)
      =ap!\,(S+S^{-1})^{p-1},\qquad
 \mathcal L_\nabla^p\bigl(h(\mathsf N)\bigr)=0.
\]
The assertion also holds when \(p=1\), with zero iterations.

By Lemma~\ref{lem:lifted-commutator}, any scalar from an earlier
step is killed by the next commutator. The scalar at step \(p\)
is therefore \(-c(\nabla,\mathcal L_\nabla^{p-1}(h(\mathsf N)))\).
By \eqref{eq:boundary-cocycle} it suffices to add the coefficients
of \(S^{\pm1}\) in \(T_0^{p-1}\) and multiply by \(ap!\). If \(p\) is odd,
\(T_0^{p-1}\) has no shifts of size \(1\) or \(-1\). If \(p\) is even, their sum is
\[
 2\binom{p-1}{p/2-1}=\binom{p}{p/2}.
\]
This proves the first identity. Since it is scalar, the next
commutator vanishes. The calculations use the row-finite algebra
of Section~\ref{sec:commutator-calculus}.
\end{proof}

The same highest commutator follows from the vertex-operator relation of
Lascoux--Thibon \cite[Section~3.2, Proposition~3.3 and (16)]{LascouxThibon}.
On partition states put \(\alpha_1=\Theta(S^{-1})\) and
\(\alpha_{-1}=\Theta(S)\), the first two modes of the family used in
Section~\ref{clt:gaussian-currents}. Under the identification of \(\ket\lambda\) with the ordinary Schur
function \(s_\lambda\) of Definition~\ref{def:symmetric-schur}, these
operators are \(\partial/\partial p_1\) and multiplication by \(p_1\),
respectively. The one-box Pieri rule \cite[Chapter~I, (5.16)]{MacdonaldSF} is
\(p_1s_\lambda=\sum_{\nu\gtrdot\lambda}s_\nu\), and differentiation
is its transpose in the Schur basis.
With \(s\) a Laurent variable and \(\vartheta=e^u\) for formal \(u\),
their normally ordered vertex operator is
\[
 V_{\mathrm{LT}}(s;\vartheta)=
 \exp\!\left(\sum_{j\ge1}\frac{(\vartheta^j-1)p_js^j}{j}\right)
 \exp\!\left(\sum_{j\ge1}(1-\vartheta^{-j})s^{-j}
                               \partial_{p_j}\right).
\]
It is read as a formal operator series; its zero mode is evaluated on
homogeneous symmetric functions, as justified in
Section~\ref{sec:hook-two-point}. Its zero mode, meaning
its coefficient of \(s^0\), is diagonal with eigenvalue
\[
 1+(\vartheta-1)(1-\vartheta^{-1})\sum_{\square\in\lambda}\vartheta^{c(\square)}
    =Q_{\psi_u}(\lambda),\qquad \psi_u(z)=e^{uz}.
\]
Here \(Q_{\psi_u}\) is defined coefficientwise in \(u\), and the equality
follows from the unsigned content sum on the left gaining
\((\vartheta-1)(1-\vartheta^{-1})\vartheta^c\) on adding a box of content \(c\), the same
increment as \(Q_{\psi_u}\).
Lascoux--Thibon's commutation rule and \(D=\alpha_1-\alpha_{-1}-I\) give
\[
 [D,V_{\mathrm{LT}}(s;\vartheta)]
 =\bigl((\vartheta-1)s+(1-\vartheta^{-1})s^{-1}\bigr)V_{\mathrm{LT}}(s;\vartheta).
\]
After \(j\) commutators, the zero mode is divisible by \(u^j\).
Its coefficient of \(u^j\) is the constant term of \((s+s^{-1})^jI\),
since \(V_{\mathrm{LT}}(s;1)=I\).
Taking \(d![u^d]\) for \(\psi(z)=z^d\), followed by linearity in
\(\psi\), gives \eqref{eq:highest-central-commutator}.

\begin{corollary}[Parity and attainment for signed-corner products]
\label{cor:corner-parity-degree}
Let \(\psi_i\ne0\) have degree \(d_i\) and leading coefficient
\(a_i\), put \(G=\prod_iQ_{\psi_i}\), and set \(w=\sum_i d_i\), as in
Theorem~\ref{thm:degree-intro}.
Then \eqref{eq:corner-product-top-intro} holds and
\[
 D^{1+\sum_i2\lfloor d_i/2\rfloor}G=0.
\]
Thus \(D^wG\ne0\) exactly when every \(d_i\) is even.
\end{corollary}
\begin{proof}
For even \(d_i\), the last possible commutator is the scalar in
Lemma~\ref{lem:highest-central-commutator}; for odd \(d_i\), the
commutator of order \(d_i\) already vanishes. A degree-zero factor
is the scalar \(a_iI\). Apply the multinomial identity to
\(\mathsf M_{\psi_1}\cdots \mathsf M_{\psi_s}\mathbf1\).
At total order \(w\), the only possible tuple is
\((d_1,\ldots,d_s)\), giving \eqref{eq:corner-product-top-intro}.
The individual bounds are \(2\lfloor d_i/2\rfloor\), giving the
stronger vanishing assertion. If every degree is even the displayed
scalar is nonzero, and otherwise it is zero.
\end{proof}

\subsection{A highest-symbol theorem}
For a polynomial \(A\) in \(q_2,q_3,\ldots\) of corner weight at most
\(w\), write \(A_{[w]}\) for its weight-\(w\) component, and set
\[
 \varepsilon(A_{[w]})=A_{[w]}(a_2,a_3,\ldots),
 \qquad a_{2r}=\binom{2r}{r},\quad a_{2r+1}=0.
\]
The statement is valid for any specified polynomial expression, so no
algebraic-independence assertion is needed.

\begin{theorem}[Pointwise polynomiality and highest symbols]
\label{thm:corner-symbol}
Let \(w\in\mathbb Z_{\ge0}\), let \(A\) have corner weight at most \(w\),
and let \(\mathsf M_A\)
denote its multiplication operator. Then
\begin{equation}\label{eq:corner-symbol-operator}
 \operatorname{ad}_D^{w+1}(\mathsf M_A)=0,\qquad
 \operatorname{ad}_D^{w}(\mathsf M_A)
          =w!\,\varepsilon(A_{[w]})I.
\end{equation}
Consequently, on every partition,
\begin{equation}\label{eq:corner-symbol-pointwise}
 D^{w+1}A=0,\qquad D^w A=w!\,\varepsilon(A_{[w]})\mathbf1.
\end{equation}
For a fixed inner partition \(\mu\), \(m=|\mu|\), the polynomial factor
in Theorem~\ref{thm:skew-exponential-polynomiality} satisfies
\begin{equation}\label{eq:corner-symbol-skew-leading}
 [x^{m+w}]P_{\mu,A}(x)=\frac{f^\mu}{m!}\varepsilon(A_{[w]}).
\end{equation}
If this scalar is nonzero, \(A\) has exact \(D\)-degree \(w\), and
\(P_{\mu,A}\) has exact degree \(m+w\) for every fixed \(\mu\).
\end{theorem}
\begin{proof}
For a monomial \(A=q_{d_1}\cdots q_{d_s}\) of weight
\(v=\sum_i d_i\), the iterated Leibniz identity distributes the
commutator order among the multiplication factors.
Lemma~\ref{lem:highest-central-commutator} bounds the \(i\)th order by
\(d_i\), and its last commutator is \(d_i!a_{d_i}I\).
At total order \(v\), only the tuple \((d_1,\ldots,d_s)\) can remain;
the multinomial coefficient cancels these factorials and gives
\(v!\prod_i a_{d_i}I\). Every order above \(v\) vanishes.
For a polynomial of weight at most \(w\), all lower-weight monomials
vanish at order \(w\); linearity proves
\eqref{eq:corner-symbol-operator}. Constants and the empty monomial
give the same assertions at \(w=0\).
Since \(D\mathbf1=0\),
\(D^jA=(\operatorname{ad}_D^j\mathsf M_A)\mathbf1\), proving
\eqref{eq:corner-symbol-pointwise}.
In the skew transfer polynomial, the only possible term of degree
\(m+w\) has \(\eta=\varnothing\) and \(j=w\).
Its coefficient is \(f^\mu D^wA(\varnothing)/(m!w!)\), which is
\eqref{eq:corner-symbol-skew-leading}. Nonvanishing proves exactness.
\end{proof}

For \(A\) of corner weight at most \(w\) and \(B\) of corner weight
at most \(v\), the highest symbols multiply.
\[
 \frac{D^{w+v}(AB)}{(w+v)!}
   =\varepsilon(A_{[w]})\varepsilon(B_{[v]})\mathbf1 .
\]
If two polynomial expressions represent the same function on
partitions, their multiplication operators coincide. Applying
\eqref{eq:corner-symbol-operator} to their difference with a common
weight bound shows that their highest-symbol values coincide.
Thus no algebraic-independence assumption is needed to use this value.

The filtration calculation takes place in the row-finite operator algebra.
For example, \(Dq_2=4|Y(\lambda)|+2\) follows by counting addable and
removable corners.

\begin{proof}[Proof of the corner-product assertions in
Theorem~\ref{thm:degree-intro}]
Expand each \(Q_{\psi_i}\) in the \(q_j\).
Its weight is at most \(\deg\psi_i\), with the stated leading term.
Theorem~\ref{thm:corner-symbol} gives the vanishing bound and the
highest scalar. The parity refinement is
Corollary~\ref{cor:corner-parity-degree}.
\end{proof}

\subsection{Content powers and skew-content products}
\label{sec:quantitative-degree}

\begin{lemma}[Content powers as signed corner sums]\label{lem:content-antiderivative}
For every nonnegative integer \(d\) there exists a polynomial \(\Phi_d\)
satisfying
\begin{equation}\label{eq:content-antiderivative}
 \Phi_d(z+1)-2\Phi_d(z)+\Phi_d(z-1)=z^d,\qquad \Phi_d(0)=0.
\end{equation}
It has degree \(d+2\), leading coefficient \(1/((d+1)(d+2))\), and
\begin{equation}\label{eq:content-corner-antiderivative}
 C_d=Q_{\Phi_d}.
\end{equation}
\end{lemma}
\begin{proof}
The centered second difference takes \(z^{j+2}\) to a polynomial of
degree \(j\) with leading coefficient \((j+1)(j+2)\). Successively
matching the coefficients of \(z^d\) gives \(\Phi_d\), and its constant
term can be chosen zero. Two such choices differ by a multiple of \(z\),
for which \(Q_z=0\). Both sides of \eqref{eq:content-corner-antiderivative} vanish at
the empty partition and increase by \(c^d\) on adding a box of content
\(c\), by Lemma~\ref{lem:corner-increments}.
\end{proof}

\begin{corollary}[Symmetric functions of squared skew contents]\label{cor:skew-content}
For the extension \(F_{\nu,\mu}\) in \eqref{eq:skew-content-extension-intro},
put \(\ell=\ell(\nu)\) and \(d=2|\nu|+2\ell\). For integers \(0\le k\le d\),
\begin{align}\label{eq:skew-content-egf}
 &e^{-x-x^2/2}
 \sum_{n\ge0}\frac{x^n}{n!}
 \sum_{|\lambda/\mu|=n} f^{\lambda/\mu}D^kF_{\nu,\mu}(\lambda)
 \notag\\
 &\qquad=
 \sum_{\eta\subseteq\mu}
 \frac{f^{\mu/\eta}x^{|\mu/\eta|}}{|\mu/\eta|!}
 \sum_{j=0}^{d-k}\frac{D^{j+k}F_{\nu,\mu}(\eta)}{j!}x^j
\end{align}
is a polynomial in \(x\) of degree at most \(|\mu|+2|\nu|+2\ell-k\).
In particular, at \(k=0\) the inner sum on the left is the skew-content
sum
\[
 \sum_{|\lambda/\mu|=n}f^{\lambda/\mu}
                  p_\nu(c(\square)^2:\square\in\lambda/\mu).
\]
The same exponential-polynomial conclusion holds for any symmetric
function \(F\) of bounded degree evaluated on the squared contents of
\(\lambda/\mu\).
\end{corollary}

\begin{proof}
By Lemma~\ref{lem:content-antiderivative},
\(C_r-C_r(\mu)=Q_{\Phi_r-C_r(\mu)}\), because the signed sum of the
constant polynomial \(1\) is \(q_0=1\).
For \(r=2\nu_i\) the polynomial in the subscript has degree \(2\nu_i+2\).
Theorem~\ref{thm:corner-symbol} now gives
\(D^{d+1}F_{\nu,\mu}=0\). Apply Lemma~\ref{lem:skew-transfer} to
obtain the generating function. Linearity extends the assertion to
symmetric functions of bounded degree by the power-sum basis.
\end{proof}

\begin{remark}
For \(\eta\subseteq\mu\), the values \(D^jF_{\nu,\mu}(\eta)\) in
\eqref{eq:skew-content-egf} are taken from the global extension
\eqref{eq:skew-content-extension-intro}, including when boxes of \(\mu\)
have been deleted. They do not refer to a skew diagram \(\eta/\mu\).
\end{remark}

\section{Content polynomials and exact sums}\label{sec:content}

Fix an inner partition \(\mu\) and put \(n=|\lambda/\mu|\).
We compute the leading skew-content moments under \(M_{\mu,n}\)
from the transfer formula and the highest central commutator.
For \(\mu=\varnothing\), these constants follow from M\'eliot's
results \cite{Meliot}.
When applying \(D\) to \(C_d(\lambda/\mu)\), we use the function
\(C_d(\lambda)-C_d(\mu)\) on all partitions. Its signed-corner
representation is given by Lemma~\ref{lem:content-antiderivative}.

\begin{lemma}\label{lem:D-content}
For every integer \(d\ge0\),
\begin{equation}\label{eq:D-Cd}
 DC_d(\lambda)=\sum_{x\in X(\lambda)}x^d+\sum_{y\in Y(\lambda)}y^d.
\end{equation}
\end{lemma}
\begin{proof}
Adding a box of content \(x\) increases \(C_d\) by \(x^d\), whereas deleting a
box of content \(y\) decreases it by \(y^d\). Substitution into the definition
of \(D\) gives the sum of these increments over both kinds of corners.
\end{proof}

\begin{proposition}\label{prop:highest-content-diff}
For every integer \(r\ge0\),
\begin{equation}\label{eq:highest-content-diff}
 D^{2r+2}C_{2r}
 =\frac{(2r)!(2r+2)!}{(r+1)!^2}\,\mathbf1,
\end{equation}
and, on every partition,
\begin{equation}\label{eq:content-vanishing}
 D^{2r+3}C_{2r}=0.
\end{equation}
For odd \(d\ge1\), multiplication by \(C_d\) satisfies
\(\operatorname{ad}_D^{d+2}(\mathsf M_{\Phi_d})=0\), and hence \(D^{d+2}C_d=0\).
\end{proposition}
\begin{proof}
Apply Lemma~\ref{lem:highest-central-commutator} to \(\Phi_d\), with
\(p=d+2\) and \(a=1/((d+1)(d+2))\). Since \(D\mathbf1=0\), the
noncommutative binomial identity gives
\[
 D^j C_d=(\operatorname{ad}_D^j\mathsf M_{\Phi_d})\mathbf1.
\]
For even \(d=2r\) the resulting scalar is
\((2r)!\binom{2r+2}{r+1}\), and for odd \(d\) it is zero.
\end{proof}

\begin{theorem}[First skew-content moment]\label{thm:first-moment}
For every fixed partition \(\mu\) and every integer \(r\ge0\), as
\(n=|\lambda/\mu|\to\infty\),
\[
 \E_{M_{\mu,n}}C_{2r}(\lambda/\mu)\sim c_r n^{r+1},
 \qquad c_r=\frac{(2r)!}{(r+1)!^2}.
\]
\end{theorem}
\begin{proof}
Subtracting \(C_{2r}(\mu)\) leaves all positive differences unchanged.
Use the top difference \eqref{eq:highest-content-diff} in
Corollary~\ref{cor:expectation-asymptotic} and divide by \((2r+2)!\).
\end{proof}

\subsection{Products and concentration}

\begin{lemma}[Top product identity]\label{lem:top-product-rule}
Let \(d_1,\ldots,d_s\) be nonnegative integers with \(s\ge1\),
put \(p_i=d_i+2\) and \(P=\sum_i p_i\), and define
\[
 L_i=d_i!\binom{d_i+2}{(d_i+2)/2}.
\]
The binomial coefficient, and thus \(L_i\), is zero when \(d_i\) is odd.
Then
\begin{equation}\label{eq:top-product-rule}
 D^P\left(\prod_{i=1}^s C_{d_i}\right)
   =\frac{P!}{p_1!\cdots p_s!}\prod_{i=1}^s L_i\,\mathbf1,
\end{equation}
and
\begin{equation}\label{eq:top-product-vanishing}
 D^{P+1}\left(\prod_{i=1}^s C_{d_i}\right)=0.
\end{equation}
The same two identities hold after replacing each \(C_{d_i}\) by
\(C_{d_i}-a_i\), for any fixed constants \(a_i\).
If \(d\) is odd, one has the stronger bound
\(D^{2d+3}((C_d-a)^2)=0\) for every fixed constant \(a\).
\end{lemma}
\begin{proof}
Lemma~\ref{lem:content-antiderivative} gives weight \(p_i=d_i+2\)
and highest coefficient \(\frac{1}{p_i(p_i-1)}\) for \(C_{d_i}\).
Apply Theorem~\ref{thm:corner-symbol} to their product.
The highest-symbol values are \(L_i/p_i!\), which gives both displayed
identities. Subtracting constants does not change them.
For odd \(d\), Corollary~\ref{cor:corner-parity-degree} applied to two
degree-\(d+2\) factors gives \(D^{2d+3}((C_d-a)^2)=0\).
\end{proof}

\begin{example}[A highest difference of a content product]
\label{ex:multinomial-C0-C2}
Let \(\mathsf M_{C_2}\) denote multiplication by \(C_2\) and recall that
\(\mathsf N_0\) multiplies by \(C_0=|\lambda|\).
Their highest commutators are
\[
 \operatorname{ad}_D^2\mathsf N_0=2I,\quad
 \operatorname{ad}_D^3\mathsf N_0=0,\qquad
 \operatorname{ad}_D^4\mathsf M_{C_2}=12I,\quad
 \operatorname{ad}_D^5\mathsf M_{C_2}=0.
\]
In the multinomial expansion of \(D^6(\mathsf N_0\mathsf M_{C_2}\mathbf1)\), the only
surviving indices are \((j_1,j_2,j_0)=(2,4,0)\): the two commutator
orders cannot exceed \(2,4\), and \(D\mathbf1=0\).
Consequently
\[
 D^6(C_0C_2)=\binom62\,2\cdot12\,\mathbf1=360\,\mathbf1,
 \qquad D^7(C_0C_2)=0.
\]
Replacing \(C_j\) by \(C_j-C_j(\mu)\) gives the same value \(360\)
for the shifted product. As a second specialization,
\(D^8(C_2^2)=\binom84\,12^2\,\mathbf1=10080\,\mathbf1\).
\end{example}

\begin{theorem}[Second and mixed leading skew-content moments]
\label{thm:second-moment}
For every fixed partition \(\mu\) and fixed nonnegative integers
\(r_1,\ldots,r_s\), as \(n=|\lambda/\mu|\to\infty\),
\[
 \E_{M_{\mu,n}}\prod_{i=1}^s C_{2r_i}(\lambda/\mu)
 \sim\left(\prod_{i=1}^s c_{r_i}\right)
               n^{\sum_i(r_i+1)}.
\]
In particular,
\[
 \E_{M_{\mu,n}}C_{2r}(\lambda/\mu)^2\sim c_r^2n^{2r+2},
 \qquad
 n^{-r-1}C_{2r}(\lambda/\mu)\longrightarrow c_r
       \quad\hbox{in }L^2(M_{\mu,n}).
\]
\end{theorem}
\begin{proof}
Apply Lemma~\ref{lem:top-product-rule} with
\(a_i=C_{2r_i}(\mu)\) and \(P=\sum_i(2r_i+2)\).
All the highest scalar factors \(L_i\) are nonzero for these even indices.
Corollary~\ref{cor:expectation-asymptotic} now gives the leading
coefficient
\[
 \prod_i\frac{L_i}{(2r_i+2)!}=\prod_i c_{r_i}.
\]
For a single even index, expanding
\(\E_{M_{\mu,n}}(n^{-r-1}C_{2r}(\lambda/\mu)-c_r)^2\)
and using the first and second moments proves the \(L^2\) assertion.
\end{proof}

\begin{proof}[Completion of the proof of Theorem~\ref{thm:degree-intro}]
Apply Lemma~\ref{lem:top-product-rule} to the shifted factors with
\(d_i=2\nu_i\) and \(a_i=C_{2\nu_i}(\mu)\). Its leading coefficient
is positive. For the empty product \(R=0\) and \(D\mathbf1=0\).
\end{proof}

\begin{example}[Highest differences and concentration]\label{ex:skew-C2-top}
For every fixed \(\mu\), put \(F_\mu=C_2-C_2(\mu)\). The scalar subtraction
does not change any positive-order commutator of its multiplication
operator. The highest commutator and top product identities give
\[
 D^4F_\mu=12\mathbf1,\qquad D^5F_\mu=0,\qquad
 D^8(F_\mu^2)=\binom84\,12^2\mathbf1=10080\mathbf1,
 \qquad D^9(F_\mu^2)=0.
\]
Consequently, under \(M_{\mu,n}\),
\[
 \mathbb EF_\mu\sim\frac{12}{4!}n^2=\frac{n^2}{2},\qquad
 \mathbb EF_\mu^2\sim\frac{10080}{8!}n^4=\frac{n^4}{4},\qquad
 \frac{C_2(\lambda/\mu)}{n^2}\longrightarrow\frac12
 \quad\hbox{in }L^2.
\]
For \(\mu=(2)\) this also follows from the exact transfer in Example~\ref{ex:exact-skew-C2}:
the numerator polynomial has leading term \(x^6/4\), whereas the
normalizer polynomial has leading term \(x^2/2\).
\end{example}

\begin{example}[An exact skew-content average]\label{ex:exact-skew-C2}
Take \(\mu=(2)\) and extend the skew statistic to all partitions by
\[
 F(\lambda)=C_2(\lambda)-C_2((2))=C_2(\lambda)-1.
\]
In particular, \(F(\varnothing)=F((1))=-1\) and \(F((2))=0\),
as prescribed by the extension to all partitions.
The local difference operator gives
\[
\begin{array}{c|rrrrr}
 \eta & F(\eta)&DF(\eta)&D^2F(\eta)&D^3F(\eta)&D^4F(\eta)\\\hline
 \varnothing&-1&0&2&8&12\\
 (1)&-1&2&10&20&12\\
 (2)&0&6&16&20&12
\end{array}
\]
and \(D^5F=0\) on every partition. Substitution into the skew transfer
therefore yields the explicit identity
\begin{align*}
 e^{-x-x^2/2}\sum_{n\ge0}\frac{x^n}{n!}
   \sum_{|\lambda/(2)|=n}f^{\lambda/(2)}C_2(\lambda/(2))
  &=5x+\frac{19}{2}x^2+\frac{25}{3}x^3
    +\frac{13}{3}x^4+\frac76x^5+\frac14x^6.
\end{align*}
For \(n=2\), the outer shapes \((4),(3,1),(2,2),(2,1,1)\) have squared
skew-content sums \(13,5,1,5\). Hence their weighted sum is
\(13+2\cdot5+1+5=29\), and the average is \(29/5\).
More generally, with
\(A_\mu(n)=\sum_{|\lambda/\mu|=n}f^{\lambda/\mu}C_2(\lambda/\mu)\),
the first values are
\[
\begin{array}{c|rrr|rrr}
 &\multicolumn{3}{c|}{\mu=(1)}&\multicolumn{3}{c}{\mu=(2)}\\
 n&Z_\mu(n)&A_\mu(n)&A_\mu(n)/Z_\mu(n)
  &Z_\mu(n)&A_\mu(n)&A_\mu(n)/Z_\mu(n)\\\hline
 0&1&0&0&1&0&0\\
 1&2&2&1&2&5&5/2\\
 2&4&14&7/2&5&29&29/5\\
 3&10&68&34/5&13&137&137/13\\
 4&26&300&150/13&38&612&306/19
\end{array}
\]
\end{example}

\subsection{The factorial content polynomial}
\label{sec:factorial-content-proof}

We first compute a single content exponential by the finite-band
calculus of Section~\ref{sec:commutator-calculus}. Throughout this
subsection, \(u\) is a formal variable and
\[
 r=e^{u/2},\qquad z=r-r^{-1}=2\sinh(u/2),\qquad
 w=r+r^{-1}.
\]
In particular \(w^2=z^2+4\). Since \(z=u+O(u^3)\) has a unique formal
inverse, the identities may also be read as power series in \(z\).
All matrix series
are interpreted coefficientwise. At any fixed power of \(u\) and
\(x\), only finitely many bands are involved.

\begin{lemma}[An exact content exponential]
\label{lem:exact-content-exponential}
Put
\[
 \mathcal C_\lambda(u)=\sum_{\square\in\lambda}e^{u c(\square)},\qquad
 \mathcal E(x,u)=\mathbb E_x\mathcal C_\lambda(u),
\]
where \(\mathbb E_x\) is the normalized single-tableau sum in
\eqref{eq:poissonized-motivation}. For integers \(d\ge0\), define
\[
 b_d(x,z)=\sum_{m\ge0}\frac{(xz)^{2m+d}}{m!(m+d)!},\qquad
 b_{-d}(x,z)=b_d(x,z).
\]
Then
\begin{equation}\label{eq:exact-content-exponential}
 \frac{\partial\mathcal E}{\partial x}
 =b_0+\sum_{d\ge1}\bigl(r^d+(-1)^d r^{-d}\bigr)b_d
                          +\frac{2b_1}{z}.
\end{equation}
The quotient \(b_1/z\) is regular at \(z=0\). Equivalently, the
undifferentiated series has the constant-term expression
\begin{equation}\label{eq:factorial-content-ct}
 1+z^2\mathcal E(x,u)=\operatorname{CT}_t
   e^{xz(t+t^{-1})}\frac{1-t^{-2}}{1-zt^{-1}-t^{-2}} .
\end{equation}
Here \(\operatorname{CT}_t=[t^0]\) selects the constant term in \(t\).
We write \(R_{\mathrm c}(t,z)=(1-t^{-2})/(1-zt^{-1}-t^{-2})\) for the
rational factor, expanded in \(t^{-1}\). At each degree in \(x\) the
exponential contributes a Laurent polynomial in \(t\).
\end{lemma}
\begin{proof}
Let \(\psi_u(t)=e^{ut}\). Adding a box of content \(c\) changes
\(Q_{\psi_u}\) by
\(e^{uc}(e^u+e^{-u}-2)=z^2e^{uc}\). Since its value at the empty
partition is \(1\),
\begin{equation}\label{eq:content-exponential-corners}
 Q_{\psi_u}(\lambda)=1+z^2\mathcal C_\lambda(u).
\end{equation}
All operator series in this proof are read coefficientwise in \(x,u\).
Each fixed coefficient is finite-band before lifting.
Lemma~\ref{lem:partition-lifts}, applied coefficientwise in \(u\), gives
\[
 \mathsf M_{\psi_u}=I+\Theta(T_{\mathrm{in}}),\qquad
 T_{\mathrm{in}}=(e^u-1)e^{u\mathsf N},\qquad D=-\Theta(\nabla)-I,
 \quad \nabla=S-S^{-1}.
\]
Using \(S^{\pm1}e^{u\mathsf N}=e^{\mp u}e^{u\mathsf N}S^{\pm1}\),
we obtain
\[
 T_x=e^{-x\operatorname{ad}_\nabla}T_{\mathrm{in}}
   =(e^u-1)e^{u\mathsf N}
           \exp\bigl(xz(rS^{-1}+r^{-1}S)\bigr).
\]
The Laurent expansion
\[
 \exp\bigl(xz(t+t^{-1})\bigr)
       =\sum_{d\in\mathbb Z}b_d(x,z)t^d
\]
shows that the coefficient of \(S^d\) in \(T_x\) is
\((e^u-1)e^{u\mathsf N}r^{-d}b_d\). Thus the scalar in
\eqref{eq:boundary-cocycle} is
\[
 -c(\nabla,T_x)=2zb_1=\frac{\partial b_0}{\partial x}.
\]
Differentiate, use Lemma~\ref{lem:lifted-commutator}, and compare
the initial values at \(x=0\). Uniqueness of the formal coefficient
recursion gives
\begin{equation}\label{eq:content-exponential-lift}
 e^{x\operatorname{ad}_D}\mathsf M_{\psi_u}=b_0 I+\Theta(T_x).
\end{equation}

For a row-finite operator \(T\), let
\(\phi(T)=(T\mathbf1)(\varnothing)\) be the sum of its vacuum row.
Since \(D\mathbf1=0\), the empty-base transfer of
Lemma~\ref{lem:skew-transfer} gives
\(\mathbb E_xQ_{\psi_u}=\phi(e^{x\operatorname{ad}_D}\mathsf M_{\psi_u})\).
This is first an identity for each polynomial coefficient of
\(\psi_u\), and hence an identity of formal series.

We now evaluate the vacuum row directly. For \(d\ge1\), a move of
length \(-d\) can end at the vacuum only by moving a particle from
\(i+d\ge0\) to \(-d\le i<0\). For \(i=-\ell\), its CAR sign is
\((-1)^{\ell-1}\): after deleting the source particle, precisely
\(\ell-1\) occupied sites lie to the right of the target. Its source
has shape \((d-\ell+1,1^{\ell-1})\), and occurs once in the row sum. Therefore
\[
 \phi\bigl(\Theta(e^{u\mathsf N}S^{-d})\bigr)
   =\sum_{\ell=1}^d(-1)^{\ell-1}r^{-2\ell}
   =\frac{1-(-r^{-2})^d}{r^2+1}.
\]
For nonnegative shifts this vacuum row is zero, including the
normal-ordered diagonal term. Substituting in
\eqref{eq:content-exponential-lift} and using
\eqref{eq:content-exponential-corners} yields
\begin{equation}\label{eq:content-exponential-integrated}
 \mathcal E(x,u)=\frac{b_0-1}{z^2}
  +\frac1{zw}\sum_{d\ge1}
                \bigl(r^d-(-1)^dr^{-d}\bigr)b_d.
\end{equation}
The first numerator is divisible by \(z^2\), each \(b_d\) is divisible
by \(z^d\), and \(w\) has constant term \(2\). To differentiate, use
\(\partial_xb_0=2zb_1\) and
\(\partial_xb_d=z(b_{d-1}+b_{d+1})\) for \(d\ge1\).
If \(A_d=r^d-(-1)^dr^{-d}\), then
\(A_0=0\), \(A_1=w\), and
\(A_{d+1}+A_{d-1}=w(r^d+(-1)^dr^{-d})\).
These identities give \eqref{eq:exact-content-exponential}.
For \eqref{eq:factorial-content-ct}, factor the denominator of the
rational function \(R_{\mathrm c}(t,z)\) on the right as
\((1-r/t)(1+r^{-1}/t)\). Its constant term is \(1\), and its
coefficient of \(t^{-d}\), for \(d\ge1\), is
\(z(r^d-(-r^{-1})^d)/w\).
Since \([t^{d}]e^{xz(t+t^{-1})}=b_d\), the constant term is
\(1+z^2\mathcal E\) by \eqref{eq:content-exponential-integrated}.
\end{proof}

The constant-term identity \eqref{eq:factorial-content-ct} can also
be obtained by applying the zero mode of \(V_{\mathrm{LT}}(s;e^u)\)
to Littlewood's sum of Schur functions and specializing to \(p_1=x\)
and \(p_j=0\) for \(j>1\). The variable in
\eqref{eq:factorial-content-ct} is \(t=e^{u/2}s\), as follows from the
vertex commutator in Section~\ref{sec:product}.
See \cite[Section~3.2, Proposition~3.3]{LascouxThibon} for the vertex operator
and Lemma~\ref{lem:external-schur-identities} for the Schur sum.
The coefficient extraction below gives the complete factorial content
polynomial and the integrality assertion.

\begin{lemma}[Central factorial expansion]
\label{lem:content-central-factorial-identities}
With \(z,r,u\) as above and an indeterminate \(c\),
\begin{equation}\label{eq:content-central-factorial-series}
 \cosh(cu)=\sum_{k\ge0}\frac{z^{2k}}{(2k)!}
                      \prod_{j=0}^{k-1}(c^2-j^2).
\end{equation}
\end{lemma}
\begin{proof}
Write \(Y(z)=\cosh(cu(z))\), where \(u(z)\) is the formal inverse
of \(2\sinh(u/2)\). Then
\[
 (1+z^2/4)Y''+(z/4)Y'=c^2Y,\qquad Y(0)=1,\quad Y'(0)=0.
\]
If \(Y=\sum_{k\ge0}a_kz^{2k}\), coefficient comparison gives
\((2k+2)(2k+1)a_{k+1}=(c^2-k^2)a_k\), which proves
\eqref{eq:content-central-factorial-series}.
\end{proof}

\begin{proof}[Proof of Theorem~\ref{thm:factorial-content-intro}]
Define \(Q_k\) initially by dividing the series on the left of
\eqref{eq:factorial-content-egf} by \(I(x)\). Its constant term
vanishes because \(\lambda=\varnothing\) contributes zero.
Conjugating \(\lambda\) preserves its weight and negates every
content. We may therefore replace \(e^{uc}\) by \(\cosh(cu)\) in
\(\mathcal E\). By \eqref{eq:content-central-factorial-series},
\begin{equation}\label{eq:factorial-content-master-series}
 1+z^2\mathcal E(x,u)=1+\sum_{k\ge0}\frac{Q_k(x)}{(2k)!}z^{2k+2}.
\end{equation}
In \eqref{eq:factorial-content-ct} expand the rational factor as
\[
 R_{\mathrm c}(t,z)=\sum_{b\ge0} z^b t^{-b}(1-t^{-2})^{-b}.
\]
Fix \(d\ge0\) and set \(n=d+1\). The coefficient of \(x^n\) in the
exponential is \(z^n t^n(1+t^{-2})^n/n!\).
To obtain \(z^{2k+2}\) we must take \(b=2k+2-n\).
To obtain \(t^0\) we must then take degree
\(h=(n-b)/2=d-k\) in \(y=t^{-2}\). If \(b<0\) or \(h<0\), the
coefficient vanishes. Otherwise \(k\le d\le 2k+1\), and
\begin{equation}\label{eq:factorial-content-integer-extraction}
 \frac{d!}{(2k)!}[x^d]Q'_k
     =[y^h]\frac{(1+y)^n}{(1-y)^b}.
\end{equation}
Indeed, differentiating the term of degree \(n\) multiplies it by
\(n\), leaving \(d!\) in the denominator.

To evaluate the right-hand side, put \(s=y/(1+y)^2\). This change
of formal variable has linear coefficient \(1\) and satisfies
\[
 \frac{\dd y}{y}=\frac{1+y}{1-y}\frac{\dd s}{s},\qquad
 \frac{1-y}{1+y}=\sqrt{1-4s}.
\]
With \(F(y)=(1+y)^n/(1-y)^b\), write \([y^h]F(y)\) as the coefficient
of \(y^{-1}\) in \(F(y)y^{-h-1}\). Apply the formal change of variable
to \(F(y)y^{-h-1}\,\dd y\). Since \(n=2h+b\), this gives
\[
 [y^h]\frac{(1+y)^n}{(1-y)^b}
   =[s^h](1-4s)^{-(b+1)/2}
   =4^h\binom{h+(b-1)/2}{h}
   =4^{d-k}\binom{d/2}{d-k}.
\]
This proves \eqref{eq:factorial-content-derivative}. Integration
with \(Q_k(0)=0\) gives \eqref{eq:factorial-content-complete}.

For integrality return to \eqref{eq:factorial-content-integer-extraction}.
The series \((1-y)^{-b}\) has nonnegative integral coefficients
and constant term \(1\), including \(b=0\). Since \(0\le h\le n\),
the coefficient extracted is at least \(\binom nh>0\). The multiplier
\((2k)!/d!\) is an integer for \(d\le 2k\). At the remaining
degree \(d=2k+1\), we have \(b=0\), and the coefficient of the
derivative is
\[
 \frac{1}{2k+1}\binom{2k+2}{k+1}=2\operatorname{Cat}_k,\qquad
 \operatorname{Cat}_k=\binom{2k}{k}-\binom{2k}{k+1}.
\]
Thus every derivative coefficient on the stated support is a positive
integer, and \(Q_k\) has degree exactly \(2k+2\). Substitution of \(d=2k\)
and \(d=k\) gives the two other values in
\eqref{eq:factorial-content-extreme-coefficients}.

For \(d=2k-2j+1\), \(0\le j\le\lfloor(k+1)/2\rfloor\), and
\(d=2k-2j\), \(0\le j\le\lfloor k/2\rfloor\), respectively,
the two ordinary-binomial factorizations of \([x^d]Q'_k\) are
\begin{align}\label{eq:factorial-content-integral-coefficients}
 \frac{2\,j!(2k)!}{(2j)!(k-2j+1)!(k-j)!}
  &=2(j!)^2\binom{2k}{2j}\binom{k-j+1}{j}
                                     \operatorname{Cat}_{k-j},\notag\\
 \frac{2^{2k-4j}(k-j)!(2k)!}{j!(2k-2j)!(k-2j)!}
  &=4^{k-2j}(2j)!\binom{2k}{2j}\binom{k-j}{j}.
\end{align}
There is also a diagrammatic check of the lowest term. For \(k\ge1\)
the summand at a box vanishes for \(|c|<k\). At size \(k+1\) only
the single row and single column contribute, each with \(H_\lambda=(k+1)!\)
and one nonzero summand \(\prod_{j=0}^{k-1}(k^2-j^2)=(2k)!/2\).
All smaller diagrams give zero. This gives
\([x^{k+1}]Q_k=(2k)!/(k+1)!\) directly. For \(k=0\), the statistic
is the size and \(Q_0=x+x^2\).
\end{proof}

\begin{example}[The first factorial content polynomials]
\label{ex:factorial-content-small}
The first four polynomials are
\begin{align*}
 Q_0&=x+x^2,\\
 Q_1&=x^2+\frac43x^3+\frac12x^4,\\
 Q_2&=4x^3+6x^4+\frac{16}{5}x^5+\frac23x^6,\\
 Q_3&=30x^4+48x^5+30x^6+\frac{64}{7}x^7+\frac54x^8.
\end{align*}
For instance, \(F_1^{\mathrm c}=C_2\), and the two diagrams of size
two each have content-square sum \(1\) and hook product \(2\).
Their total contribution is \(x^2\), explaining the lowest term
of \(Q_1\). Differentiating gives
\begin{align*}
 Q'_0&=1+2x,\\
 Q'_1&=2x+4x^2+2x^3,\\
 Q'_2&=12x^2+24x^3+16x^4+4x^5,\\
 Q'_3&=120x^3+240x^4+180x^5+64x^6+10x^7.
\end{align*}
For example, the two lines of
\eqref{eq:factorial-content-integral-coefficients} give
\([x^3]Q'_3=120\) and \([x^4]Q'_3=240\), respectively.
The coefficients \([x^7]Q_3=64/7\) and \([x^8]Q_3=5/4\)
illustrate the distinction between \(Q_k\) and its derivative.
\end{example}

The binomial version of the statistic has a different normalization.
For a nonnegative integer \(m\), use the same polynomial convention
\(\binom{t}{m}=t(t-1)\cdots(t-m+1)/m!\), with \(\binom t0=1\).
The identity
\[
 \prod_{j=0}^{k-1}(c^2-j^2)
   =\frac{(2k)!}{2}
        \left(\binom{c+k}{2k}+\binom{c+k-1}{2k}\right)
\]
also holds for \(k=0\). Since
\(\binom{-c+k}{2k}=\binom{c+k-1}{2k}\), diagram conjugation gives
\begin{equation}\label{eq:factorial-content-binomial-normalization}
 e^{-x-x^2/2}\sum_\lambda\frac{x^{|\lambda|}}{H_\lambda}
       \sum_{\square\in\lambda}\binom{c(\square)+k}{2k}
                         =\frac{Q_k(x)}{(2k)!}.
\end{equation}
The integrality assertion in Theorem~\ref{thm:factorial-content-intro}
uses the product normalization in
\eqref{eq:factorial-content-egf}.

\section{Exact moment identities from the difference operator}\label{sec:exact-moments}

\subsection{Observables and the exact Poissonized identity}
\label{clt:operator-scope}

We compute exact moments by conjugating multiplication operators with \(D\).
The resulting identities give the hook coefficients in
Section~\ref{sec:hook-polynomials} and the joint skew Gaussian limits in
Section~\ref{sec:limit-shape}. With the corner powers \(q_d=Q_{z^d}\) of
Definition~\ref{def:corner-observables-intro}, let \(\mathsf M_d\) be the operator
on functions of partitions given by
\[
  (\mathsf M_d F)(\lambda)=q_d(\lambda)F(\lambda).
\]
Thus \(\mathsf M_d\) is diagonal in the partition basis, with entry
\(q_d(\lambda)\) at \(\lambda\). Its particle representation will be used in
Section~\ref{clt:finite-band}.

We use the Poissonized expectation \(\mathbb E_x\) defined in
\eqref{eq:poissonized-motivation}. An observable \(F\) has polynomial
growth if \(|F(\lambda)|\le C(1+|\lambda|)^K\) for some constants
\(C,K\ge0\) independent of \(\lambda\). A diagram of size \(N\) has at
most \(2N+1\) addable and removable corners in total, and their contents
have absolute value at most \(N\). It follows that \(D\) preserves
polynomial growth and every fixed product of the \(q_d\) has this
growth bound.

\begin{lemma}[Poissonized differentiation]\label{clt:poisson-derivative}
For every observable of polynomial growth,
\[
 \frac{d}{dx}\mathbb E_xF=\mathbb E_x(DF).
\]
In particular, if \(D^{R+1}F=0\) for an integer \(R\ge0\), then
\begin{equation}\label{clt:poisson-taylor}
 \mathbb E_xF=\sum_{j=0}^{R}\frac{x^j}{j!}D^jF(\varnothing).
\end{equation}
\end{lemma}
\begin{proof}
Set \(b_n(F)=\sum_{|\lambda|=n}f^\lambda F(\lambda)\) and \(b_{-1}(F)=0\).
The recurrence \eqref{eq:transfer-recurrence} for \(b_n(F)\) gives
\(b_n(DF)=b_{n+1}(F)-b_n(F)-n b_{n-1}(F)\). Sum against \(x^n/n!\)
and differentiate the normalization using \(I'(x)=(1+x)I(x)\).
The series for both expectations converge locally uniformly in
complex \(x\), since \(I(x)=e^{x+x^2/2}\) is entire and polynomial factors
in the rank correspond to applying powers of \(x\,d/dx\) to \(I(x)\).
Iterating the identity and using \(\mathbb E_0F=F(\varnothing)\)
proves \eqref{clt:poisson-taylor}.
\end{proof}

\subsection{Finite-band operators and finite matrix elements}
\label{clt:finite-band}

We use the matrix convention of Definition~\ref{def:partition-matrix-element}
and the row-finite lifts of Section~\ref{sec:finite-band-operators}. For the
bilateral shift \(S e_i=e_{i+1}\) and the position operator \(\mathsf N e_i=i e_i\), put
\[
 \nabla=S-S^{-1},\qquad T_0=S+S^{-1},\qquad g_d(i)=(i+1)^d-i^d .
\]
Lemma~\ref{lem:partition-lifts} identifies \(D\) with \(-\Theta(\nabla)-I\) and,
for \(d\ge1\), \(\mathsf M_d\) with \(\Theta(g_d(\mathsf N))\).
On the partition sector, \(\Theta(I)=0\) because the number of particles
relative to the vacuum is zero. The total displacement from the vacuum is
\(\sum_i((\lambda_i-i)-(-i))=|\lambda|\), so
\(\Theta(\mathsf N)|\lambda\rangle=|\lambda|\,|\lambda\rangle\).
In particular \(g_1=1\) and \(g_2(i)=2i+1\) give \(q_1=0\) and
\(q_2(\lambda)=2|\lambda|\).

We extend the vacuum row functional \(\phi\) to an arbitrary row.
For a row-finite matrix \(C\) and a partition \(\eta\), put
\begin{equation}\label{clt:phi}
 \phi_\eta(C)=(C\mathbf1)(\eta)
     =\sum_\lambda\langle\eta|C|\lambda\rangle,\qquad
 \phi=\phi_\varnothing .
\end{equation}
Here \(\mathbf1\) is the constant partition function. The notation
\(\langle\varnothing|C|\mathbf1\rangle\) denotes the finite sum of the
vacuum row. For formal operator series we evaluate each coefficient
separately.

\subsection{Polynomial conjugation and exact moment identities}
\label{clt:conjugation}

For every positive integer \(d\), Proposition~\ref{prop:local-nilpotence} gives
\(\operatorname{ad}_D^{d+1}(\mathsf M_d)=0\). The formal conjugate
\(e^{xD}\mathsf M_de^{-xD}\) is therefore the polynomial
\[
 \mathsf M_d(x)=\sum_{j=0}^{d}\frac{x^j}{j!}
             \operatorname{ad}_D^j(\mathsf M_d),\qquad
 \operatorname{ad}_D(A)=[D,A].
\]

\begin{lemma}[Exact Poissonized operator moments]
\label{clt:operator-moments-lemma}
For positive integers \(d_1,\ldots,d_s\),
\begin{equation}\label{clt:poissonized-operator-moments}
 \mathbb E_x\prod_{i=1}^s q_{d_i}
     =\phi\left(\prod_{i=1}^s\mathsf M_{d_i}(x)\right).
\end{equation}
The operators \(\mathsf M_d(x)\) commute with one another.
\end{lemma}
\begin{proof}
Put \(B=\mathsf M_{d_1}\cdots\mathsf M_{d_s}\), so \(B\mathbf1=\prod_iq_{d_i}\).
The identity \(D\mathbf1=0\) gives
\(D^r(B\mathbf1)=\operatorname{ad}_D^r(B)\mathbf1\) by induction.
By Lemma~\ref{lem:noncommutative-multinomial}, the formal exponential
\(e^{x\operatorname{ad}_D}\) preserves products. Applied to \(B\), it gives the
product of the polynomials \(\mathsf M_{d_i}(x)\), of degree at most
\(R=\sum_i d_i\). Thus \(D^{R+1}(B\mathbf1)=0\), and
Lemma~\ref{clt:poisson-derivative} gives
\[
 \mathbb E_x(B\mathbf1)
 =\sum_{r=0}^{R}\frac{x^r}{r!}\phi\bigl(\operatorname{ad}_D^r(B)\bigr)
 =\phi\left(\prod_{i=1}^{s}\mathsf M_{d_i}(x)\right).
\]
The same product identity preserves the commutativity of the
multiplication operators \(\mathsf M_d\).
\end{proof}

For example, at \(d=2\) the constant, linear and quadratic coefficients of
\(\phi(\mathsf M_2(x))\) are \(0,2,2\). Indeed \(q_2=2|\lambda|\),
\(Dq_2(\varnothing)=2\) and \(Dq_2((1))=6\), so
\(D^2q_2(\varnothing)=4\). This gives \(\mathbb E_x q_2=2x+2x^2\),
in agreement with \(2xI'(x)/I(x)\).

\subsection{The Gaussian law of the limiting currents}
\label{clt:gaussian-currents}

For positive integers \(r\), define the currents (lifts of the shifts
by \(r\) sites) by
\[
 \alpha_r=\Theta(S^{-r}),\qquad \alpha_{-r}=\Theta(S^r),
 \qquad A_r=\alpha_r+\alpha_{-r}.
\]
Lemma~\ref{lem:lifted-commutator} gives the Heisenberg relation
\[
 [\alpha_r,\alpha_s]=r\,\delta_{r,-s}I\qquad(r,s\in\mathbb Z\setminus\{0\}).
\]
Here \(\delta_{r,-s}\) is the Kronecker delta. Opposite shifts of
length \(r>0\) cross the vacuum boundary at
\(r\) sites, and shifts with nonzero net displacement have zero central term.
In particular the \(A_r\) commute. The vacuum also satisfies
\[
 \alpha_r|\varnothing\rangle=0,\qquad
 \langle\varnothing|\alpha_{-r}=0\qquad(r>0).
\]
Every occupied vacuum site \(j<0\) has \(j-r\) occupied as well, so
the left shift has no allowed move. The second relation is its transpose:
reversing a particle move preserves the sign.

\begin{definition}[Symmetric functions and Schur functions]
\label{def:symmetric-schur}
Let \(\Lambda=\mathbb C[p_1,p_2,\ldots]\) be the ring of symmetric
functions over \(\mathbb C\), obtained from \(\Lambda_{\mathbb Q}\) by
extension of scalars. Its grading is given by \(\deg p_r=r\).
Evaluation in \(L\) variables sends \(p_r\) to \(\sum_{i=1}^Lz_i^r\).
These evaluations are symmetric and compatible with setting the last
variable to zero. Write \(\Lambda_k\) for its degree-\(k\) subspace. The completion
by degree is \(\widehat\Lambda=\prod_{k\ge0}\Lambda_k\), with products
computed degree by degree. Define \(h_j\) by
\[
 \sum_{j\ge0}h_jv^j
   =\exp\left(\sum_{r\ge1}\frac{p_rv^r}{r}\right),
 \qquad h_0=1,\qquad h_j=0\ (j<0).
\]
For partitions \(\eta\subseteq\lambda\) and any
\(L\ge\max(\ell(\lambda),\ell(\eta))\), define
\[
 s_{\lambda/\eta}
   =\det(h_{\lambda_i-\eta_j-i+j})_{1\le i,j\le L},
 \qquad s_\lambda=s_{\lambda/\varnothing}.
\]
The determinant is independent of \(L\). These conventions agree with
\cite{MacdonaldSF,StanleyEC2,MacdonaldLittlewood,OkounkovWedge}.
\end{definition}

The first identity below is the transpose, with integer indices, of
\cite[Appendix~A, (A.15)--(A.16)]{OkounkovWedge}. The second is
Littlewood's Schur-sum identity (see \cite[p.~160]{MacdonaldLittlewood}).

\begin{lemma}[Schur matrix elements and the Schur sum; \cite{MacdonaldSF,MacdonaldLittlewood,OkounkovWedge}]
\label{lem:external-schur-identities}
Let \(t_1,t_2,\ldots\) be formal variables of finite support, and use the
partition states and currents \(\alpha_r\) defined above. For
\(\eta\subseteq\lambda\),
\[
 \left\langle\eta\left|
     \exp\left(\sum_{r\ge1}t_r\alpha_r\right)
     \right|\lambda\right\rangle
      =s_{\lambda/\eta}\big|_{p_r=rt_r}.
\]
The matrix element is zero when \(\eta\not\subseteq\lambda\).
In \(\widehat\Lambda\),
\[
 \sum_{\lambda}s_\lambda
 =\prod_i(1-z_i)^{-1}\prod_{i<j}(1-z_iz_j)^{-1}
 =\exp\left(\sum_{r\text{ odd}}\frac{p_r}{r}
                       +\frac12\sum_{r\ge1}\frac{p_r^2}{r}\right).
\]
Infinite products mean their compatible finite-variable evaluations.
The specialized identities are coefficientwise formal: every specified
weighted coefficient involves only finitely many partitions.
\end{lemma}

\begin{lemma}[Gaussian currents]\label{clt:currents-lemma}
For a finitely supported family \((t_r)_{r\ge1}\) of formal parameters, let
\(e_r=\mathbf1_{\{r\ {\rm odd}\}}\) be the indicator that \(r\) is odd. Then
\begin{equation}\label{eq:positive-current-gaussian}
 \phi\!\left(\exp\!\left(\sum_{r\ge1}t_r\alpha_r\right)\right)
 =\exp\!\left(\sum_{r\ge1}e_rt_r+\tfrac12\sum_{r\ge1}r t_r^2\right)
\end{equation}
and
\begin{equation}\label{clt:current-gaussian}
 \phi\!\left(\exp\!\left(\sum_{r\ge1}t_r A_r\right)\right)
 =\exp\!\left(\sum_{r\ge1}(e_rt_r+r t_r^2)\right).
\end{equation}
Thus the polynomial moments of \((A_1,A_2,\ldots)\) under \(\phi\) are exactly those
of independent real Gaussian variables of means \(e_r\) and variances \(2r\).
\end{lemma}
\begin{proof}
Apply \(\langle\varnothing|\) on the left and sum the Schur matrix
elements of Lemma~\ref{lem:external-schur-identities} over the input
partition. Littlewood's sum at \(p_r=rt_r\) gives
\eqref{eq:positive-current-gaussian}. Put \(U=\sum_r t_r\alpha_{-r}\),
\(V=\sum_r t_r\alpha_r\), and \(C=[V,U]=\sum_r rt_r^2I\).
Since \(C\) is central,
\(e^{tU}Ve^{-tU}=V-tC\). For
\(\mathcal H(t)=e^{tU}e^{tV}e^{t^2C/2}\), differentiation gives
\(\mathcal H'(t)=(U+V)\mathcal H(t)\) and \(\mathcal H(0)=I\).
Uniqueness of the formal coefficient recursion therefore gives
\[
 e^{U+V}=e^Ue^V\exp\left(\frac12\sum_r r t_r^2\right).
\]
The vacuum bra kills every positive-degree term of \(e^U\), so applying \(\phi\)
proves \eqref{clt:current-gaussian}. The final assertion follows by coefficient
comparison with the usual Gaussian moment generating function.
\end{proof}

\subsection{The two leading terms of the conjugated observables}
\label{clt:leading-conjugation}

With the coefficients \(a_d,b_{dr}\) from \eqref{eq:intro-clt-coefficients}, set
\begin{equation}\label{clt:transfer-coefficients}
 L_d=\sum_{r\ge1}b_{dr}A_r,\qquad d\ge2.
\end{equation}
Only \(r\le d-1\) with \(r\equiv d-1\pmod2\) occur.

\begin{proposition}[Leading conjugation formula]
\label{clt:leading-conjugation-prop}
For every \(d\ge2\) there are fixed row-finite matrices \(O_{d,j}\) such that
\begin{equation}\label{clt:operator-leading}
 \mathsf M_d(x)=a_d x^d I+x^{d-1}L_d
                    +\sum_{j=0}^{d-2}x^j O_{d,j}.
\end{equation}
The evaluation by \(\phi\) of any finite product of the \(O_{d,j}\)
and \(A_r\) is a finite scalar independent of \(x\) and the diagram size.
\end{proposition}
\begin{proof}
By Lemma~\ref{lem:lifted-commutator} there is a scalar polynomial \(c_d(x)\) such that
\[
 \mathsf M_d(x)=\Theta(\mathsf T_d(x))+c_d(x)I,
 \qquad
 \mathsf T_d(x)=e^{-x\operatorname{ad}_\nabla}g_d(\mathsf N),
 \qquad c_d(0)=0.
\]
The one-particle commutators satisfy \([\nabla,\mathsf N]=-T_0\) and
\([\nabla,T_0]=0\). The conjugation sends \(\mathsf N\) to \(\mathsf N+xT_0\)
and preserves products. Thus
\[
 \mathsf T_d(x)=g_d(\mathsf N+xT_0)
      =d x^{d-1}T_0^{d-1}+\sum_{j=0}^{d-2}x^j \mathsf T_{d,j},
\]
where every \(\mathsf T_{d,j}\) is a finite-band matrix with polynomial
coefficients. The highest power of \(x\) in \(g_d(\mathsf N+xT_0)\)
comes from choosing \(xT_0\) in every factor of its term of degree \(d-1\).
Since \(\Theta(I)=0\), the binomial expansion of \(T_0^{d-1}\) gives
\(d\Theta(T_0^{d-1})=L_d\).

The scalar in the lifted commutator is given by \eqref{eq:boundary-cocycle}.
Comparing the scalar terms in
\(\partial_x\mathsf M_d(x)=[D,\mathsf M_d(x)]\) gives
\(c_d'(x)=-c(\nabla,\mathsf T_d(x))\). Thus the coefficient of \(x^d\)
in \(c_d(x)\) is
\[
 -c(\nabla,T_0^{d-1})=
 \begin{cases}
 2\binom{d-1}{d/2-1}=\binom{d}{d/2},&d\text{ even},\\
 0,&d\text{ odd}.
 \end{cases}
\]

It remains to exclude a scalar term of degree \(d-1\). Every lifted matrix has
zero vacuum diagonal, so
\(c_d(x)=\langle\varnothing|\mathsf M_d(x)|\varnothing\rangle\).
Let \(J|\lambda\rangle=|\lambda'\rangle\) be diagram conjugation. It preserves
\(D\), fixes the vacuum, and sends \(\mathsf M_d\) to \((-1)^d\mathsf M_d\).
Thus \(c_d(x)=0\) for odd \(d\). For even \(d\) let
\(P|\lambda\rangle=(-1)^{|\lambda|}|\lambda\rangle\). Since each off-diagonal
term of \(D\) changes size by one,
\[
 PDP=-D-2I,
 \qquad P\mathsf M_dP=\mathsf M_d.
\]
Conjugating each iterated commutator therefore gives
\(P\mathsf M_d(x)P=\mathsf M_d(-x)\). The vacuum is fixed by \(P\), so \(c_d\) is
an even polynomial. Its degree-\(d-1\) coefficient is again zero. Thus
\(c_d(x)-a_dx^d\) has degree at most \(d-2\). Combining the scalar and
lifted polynomial remainders proves \eqref{clt:operator-leading}.
Row-finiteness proves the final assertion.
\end{proof}

For example, \(a_2=2\), \(L_2=2A_1\), \(L_3=3A_2\), and
\(L_4=4(A_3+3A_1)\). The first current records the Poissonized size fluctuation.

\section{Hook-length polynomials}\label{sec:hook-polynomials}

We express even hook sums in corner and content coordinates, then use
their highest differences to prove Theorems~\ref{thm:hook-polynomiality-intro}
and~\ref{thm:factorial-hook-intro}.

\subsection{A finite boundary identity for hook lengths}

For a nonnegative integer \(r\), recall that
\[
 \mathcal H_r(\lambda)=\sum_{\square\in\lambda}
                           h_\lambda(\square)^{2r},
 \qquad \mathcal H_0(\lambda)=|\lambda|.
\]
Thus \(\mathcal H_r\) is the \(r\)th power sum of the squared hook
lengths. Let
\[
 \mathcal Q_\lambda(z)=\sum_{a\in X(\lambda)}z^a
                         -\sum_{b\in Y(\lambda)}z^b.
\]
This is a Laurent polynomial in \(z\). At \(z=e^t\) it has the formal
expansion \(\sum_{j\ge0}q_j(\lambda)t^j/j!\).

\begin{lemma}[Hook lengths as boundary differences]
\label{lem:hook-corner-identity}
For every partition \(\lambda\),
\begin{equation}\label{eq:hook-corner-identity}
 \sum_{\square\in\lambda}
        \bigl(z^{h_\lambda(\square)}+z^{-h_\lambda(\square)}\bigr)
 =\frac{\mathcal Q_\lambda(z)\mathcal Q_\lambda(z^{-1})-1}
                    {(z-1)(1-z^{-1})}.
\end{equation}
The quotient on the right is a Laurent polynomial. In particular,
as an identity of formal power series,
\begin{equation}\label{eq:hook-cosh-identity}
 \sum_{\square\in\lambda}\cosh\bigl(t h_\lambda(\square)\bigr)
 =\frac{\mathcal Q_\lambda(e^t)\mathcal Q_\lambda(e^{-t})-1}
                          {8\sinh^2(t/2)}.
\end{equation}
The numerator in \eqref{eq:hook-cosh-identity} is divisible by \(t^2\).
\end{lemma}
\begin{proof}
Choose \(L\ge\max\{1,\ell(\lambda)\}\), pad \(\lambda\) by zero
parts, and put
\[
 \beta_i=\lambda_i+L-i\quad(1\le i\le L),\qquad
 \mathcal B=\{\beta_1,\ldots,\beta_L\},\qquad \mathfrak v(z)=\sum_{a\in \mathcal B}z^a.
\]
The \(\beta_i\), called beta numbers, are distinct nonnegative integers.
We use the hook-length multiset description of
\cite[Lemma~4.1]{StanleyHooks}. For a box \((i,j)\) put
\[
 a=\lambda_i+L-i,\qquad b=L+j-\lambda'_j-1.
\]
The integers \(b_j=L+j-\lambda'_j-1\), \(1\le j\le\lambda_1\),
increase strictly and enumerate the complement of \(\mathcal B\) in
\(\{0,\ldots,\lambda_1+L-1\}\). Indeed, if \(i\le\lambda'_j\),
then \(\beta_i>b_j\), whereas if \(i>\lambda'_j\), then
\(\beta_i<b_j\). None of the \(b_j\) belongs to \(\mathcal B\), and their
number is the size of that complement. Moreover,
\(\beta_i>b_j\) holds precisely when \(j\le\lambda_i\).
Thus boxes correspond bijectively to pairs
\[
 a\in \mathcal B,\qquad b\notin \mathcal B,\qquad 0\le b<a,
 \qquad h_\lambda(i,j)=a-b.
\]

First sum over all \(0\le b<a\) and then remove the terms with
\(b\in \mathcal B\). The first sum is
\[
 \sum_{a\in \mathcal B}\sum_{h=1}^{a}(z^h+z^{-h})
       =\frac{z\mathfrak v(z)-\mathfrak v(z^{-1})}{z-1}-L.
\]
The removed sum is \(\mathfrak v(z)\mathfrak v(z^{-1})-L\), since a pair of distinct
elements of \(\mathcal B\) contributes both opposite differences. Hence the
hook sum on the left of \eqref{eq:hook-corner-identity} is
\[
 \frac{z\mathfrak v(z)-\mathfrak v(z^{-1})}{z-1}-\mathfrak v(z)\mathfrak v(z^{-1}).
\]
The occupied boundary positions are \(\beta_i-L\), together with
all integers below \(-L\). Their signed changes give
\[
 \mathcal Q_\lambda(z)=z^{-L}\bigl(1+(z-1)\mathfrak v(z)\bigr).
\]
For example, an occupied position \(a-L\) contributes
\(z^{a-L+1}-z^{a-L}\); the lower occupied tail contributes
\(z^{-L}\). Multiplying the last display by its version at \(z^{-1}\)
now proves \eqref{eq:hook-corner-identity}. Substitution \(z=e^t\)
and division by two prove \eqref{eq:hook-cosh-identity}.
\end{proof}

The polynomiality assertion below also follows from
Han--Xiong \cite[Theorem~7.1]{HanXiong}: their factorial hook sums form
a triangular basis of the even hook sums. The boundary identity gives
the following explicit corner and content formulas.

\begin{corollary}[Hook sums in corner and content coordinates]
\label{cor:hook-corner-polynomial}
For each nonnegative integer \(r\), the hook sum \(\mathcal H_r\) is a
polynomial over \(\mathbb Q\) in \(q_2,\ldots,q_{2r+2}\) of corner
weight at most \(2r+2\).
Its component of corner weight \(2r+2\) is
\begin{equation}\label{eq:hook-highest-weight}
 \frac{(2r)!}{2}
 \sum_{a+b=2r+2}\frac{(-1)^b q_aq_b}{a!b!},
 \qquad q_0=1,\quad q_1=0.
\end{equation}
In content coordinates,
\begin{equation}\label{eq:hook-content-polynomial}
 \mathcal H_r
 =C_{2r}+
   \sum_{\substack{\ell\ge1,\ a,b\ge0\\a+b+2\ell=2r}}
       \frac{(2r)!(-1)^b}{a!b!(2\ell)!}\,C_aC_b.
\end{equation}
An empty sum is zero. The indices \(a,b\) in this formula may both
be odd.
\end{corollary}
\begin{proof}
The numerator of \eqref{eq:hook-cosh-identity} has zero constant
and linear terms. The series \(t^2/(8\sinh^2(t/2))\) is even and
has constant term \(1/2\). Extracting the coefficient of \(t^{2r}\)
proves polynomiality and \eqref{eq:hook-highest-weight}; all other
terms have weights \(2r,2r-2,\ldots\).

For the content formula, put
\(\mathcal C_\lambda(t)=\sum_{\square\in\lambda}e^{t c(\square)}\)
and \(s(t)=e^t+e^{-t}-2\). The one-box corner change gives
\[
 \mathcal Q_\lambda(e^t)=1+s(t)\mathcal C_\lambda(t).
\]
Equation \eqref{eq:hook-cosh-identity} consequently becomes
\[
 \sum_{\square}\cosh(t h_\lambda(\square))
 =\frac12\bigl(\mathcal C_\lambda(t)+\mathcal C_\lambda(-t)
               +s(t)\mathcal C_\lambda(t)\mathcal C_\lambda(-t)\bigr).
\]
Use \(s(t)=2\sum_{\ell\ge1}t^{2\ell}/(2\ell)!\) and extract the
coefficient of \(t^{2r}\) to obtain \eqref{eq:hook-content-polynomial}.
\end{proof}

\begin{example}[The first even hook sums]\label{ex:hook-content-small}
Writing \(n=|\lambda|\), the first formulas are
\[
 \mathcal H_0=n,\qquad
 \mathcal H_1=C_2+n^2,\qquad
 \mathcal H_2=C_4+12nC_2-12C_1^2+n^2.
\]
For \(\lambda=(2)\), with contents \(0,1\) and hook lengths \(2,1\),
the last formula gives \(1+24-12+4=17=2^4+1^4\).
\end{example}

\subsection{Proof of the Gelfand polynomiality theorem}
\label{sec:gelfand-stanley-proof}

We first express the shifted parts in the same corner coordinates as the
contents and hooks. Put \(\mathcal R_j(\lambda)=p_j(\boldsymbol a_\lambda)\)
for \(j\ge1\).

\begin{lemma}[Shifted parts in corner coordinates]
\label{lem:shifted-parts-corners}
For \(j\ge1\), let \(\mathfrak s_j(t)\in\mathbb Q[t]\) be the unique
polynomial satisfying
\[
 \mathfrak s_j(0)=0,\qquad
 \mathfrak s_j(t+1)-\mathfrak s_j(t)=t^j.
\]
It has degree \(j+1\). With \(N=|\lambda|\),
\begin{equation}\label{eq:shifted-parts-content-polynomial}
 \mathcal R_j(\lambda)=\mathfrak s_j(N)
    +\sum_{\square\in\lambda}
       \bigl((N+c(\square))^j-(N+c(\square)-1)^j\bigr).
\end{equation}
Consequently \(\mathcal R_j\) is a rational polynomial in the \(q_d\)
of corner weight at most \(2j+2\).
\end{lemma}
\begin{proof}
The map \(A(t)\mapsto A(t+1)-A(t)\) takes \(t^{m+1}\) to a polynomial
of degree \(m\) with leading coefficient \(m+1\). Successive coefficient
matching gives \(\mathfrak s_j\), and its value at zero fixes the additive
constant uniquely. At each integer \(N\ge0\), telescoping gives
\(\mathfrak s_j(N)=\sum_{i=1}^N(N-i)^j\). In row \(i\), with box
columns \(b=1,\ldots,\lambda_i\), the sum
\[
 \sum_{b=1}^{\lambda_i}\bigl((N+b-i)^j-(N+b-i-1)^j\bigr)
       =(N+\lambda_i-i)^j-(N-i)^j
\]
telescopes. Sum over \(1\le i\le N\), padding by zero rows, to prove
\eqref{eq:shifted-parts-content-polynomial}, including \(N=0\).
Now \(N=q_2/2\) has corner weight two, and
Lemma~\ref{lem:content-antiderivative} expresses \(C_b\) with weight
at most \(b+2\). Every term in the difference
\((N+c)^j-(N+c-1)^j\) has the form \(N^ac^b\) with \(a+b\le j-1\).
After summation its weight is at most \(2a+b+2\le2j\). The remaining
term \(\mathfrak s_j(N)\) has weight at most \(2j+2\).
\end{proof}

\begin{proof}[Proof of Theorem~\ref{thm:gelfand-stanley-polynomiality}]
The evaluations of \(p_j\) on the three lists are respectively
\(\mathcal H_j\), \(C_j\) and \(\mathcal R_j\). They have corner weights
at most \(2j+2\), \(j+2\) and \(2j+2\), by
Corollary~\ref{cor:hook-corner-polynomial},
Lemma~\ref{lem:content-antiderivative} and
Lemma~\ref{lem:shifted-parts-corners}. Hence each power-sum monomial
in \(g_F\) has weight at most the corresponding term in \(\mathfrak w(F)\).
The corner-product bound of Theorem~\ref{thm:corner-symbol}, applied to
each monomial in the \(q_d\), gives \(D^{\mathfrak w(F)+1}g_F=0\).
All coefficients are rational. The fixed-base transfer in
Lemma~\ref{lem:skew-transfer} therefore gives the asserted generating
functions, the degree bound and the stated finite formula for
\(P_{\mu,F}\). Setting \(\mu=\varnothing\) gives the ordinary identity.
\end{proof}

\begin{corollary}[Highest shifted-part differences]
\label{cor:shifted-parts-highest-differences}
Fix a partition \(\mu\) and integers \(j_1,\ldots,j_s\ge1\). Put
\[
 F^{\mathrm a}_{\boldsymbol j,\mu}(\lambda)
   =\prod_{i=1}^s
       \bigl(\mathcal R_{j_i}(\lambda)-\mathcal R_{j_i}(\mu)\bigr),
 \qquad W=2\sum_{i=1}^s(j_i+1).
\]
As a function on all partitions, this product has exact \(D\)-degree
\(W\), and
\begin{equation}\label{eq:shifted-parts-highest-difference}
 D^W F^{\mathrm a}_{\boldsymbol j,\mu}
   =W!\prod_{i=1}^s\frac{1}{j_i+1}\,\mathbf1,
 \qquad D^{W+1}F^{\mathrm a}_{\boldsymbol j,\mu}=0.
\end{equation}
In particular, for every \(j\ge1\),
\begin{equation}\label{eq:shifted-part-single-highest-difference}
 D^{2j+2}\mathcal R_j
       =\frac{(2j+2)!}{j+1}\,\mathbf1,
 \qquad D^{2j+3}\mathcal R_j=0.
\end{equation}
The empty product has \(W=0\) and satisfies the same convention.
\end{corollary}
\begin{proof}
The polynomial \(\mathfrak s_j\) in Lemma~\ref{lem:shifted-parts-corners}
has leading coefficient \(1/(j+1)\). That lemma therefore gives
\[
 \mathcal R_j=\frac{N^{j+1}}{j+1}+B_j,
 \qquad N=|\lambda|=q_2/2,
 \qquad \operatorname{wt}B_j\le2j.
\]
Assume \(s\ge1\), and put \(K=\sum_i(j_i+1)\), so \(W=2K\).
Subtracting the fixed values at \(\mu\) and multiplying yields
\[
 F^{\mathrm a}_{\boldsymbol j,\mu}
    =\left(\prod_i\frac1{j_i+1}\right)N^K+B,
 \qquad \operatorname{wt}B\le W-2.
\]
The corner-product bound of Theorem~\ref{thm:corner-symbol} gives
\(D^{W-1}B=0\). Applying Lemma~\ref{lem:top-product-rule} to
\(K\) copies of \(C_0=N\) gives
\(D^{2K}N^K=(2K)!\mathbf1\) and \(D^{2K+1}N^K=0\).
This proves \eqref{eq:shifted-parts-highest-difference}; its scalar
is positive, so the degree is exact. The single-factor formula follows,
and the empty product is immediate from \(D\mathbf1=0\).
\end{proof}

\begin{example}[Hook squares and shifted parts]
\label{ex:gelfand-stanley-small}
Write \(N(\lambda)=|\lambda|\). The first shifted-part sum is
\(\mathcal R_1=N(N+1)/2\). The empty-base transfer gives
\[
 \mathbb E_xN=x+x^2,\qquad
 \mathbb E_xN^2=x+3x^2+2x^3+x^4,
\]
since applying \((x\,d/dx)^2\) to \(I(x)=e^{x+x^2/2}\) counts the
squared size. Thus
\[
 \mathbb E_x\mathcal R_1=x+2x^2+x^3+\frac12x^4.
\]
Also \(\mathcal H_1=C_2+N^2\), and
Example~\ref{ex:factorial-content-small} gives
\(\mathbb E_xC_2=x^2+4x^3/3+x^4/2\). Hence
\[
 \mathbb E_x\mathcal H_1=x+4x^2+\frac{10}{3}x^3+\frac32x^4.
\]
For the base \(\mu=(1)\), the transfer polynomial is
\(P_{(1),F}=(1+x)P_{\varnothing,F}+P'_{\varnothing,F}\): indeed,
\(\mathsf U=D+\mathsf L+I\) and \((\mathsf Lh)(\varnothing)=0\) give
\(D^jg((1))=D^{j+1}g(\varnothing)+D^jg(\varnothing)\).
In particular, for \(F=p_1(\boldsymbol h)\),
\[
 P_{(1),F}=1+9x+15x^2+\frac{40}{3}x^3+\frac{29}{6}x^4+\frac32x^5.
\]
Its constant term is \(\mathcal H_1((1))=1\); its coefficient of \(x\),
after multiplication by \(I(x)\), is \(10\), the sum of the squared-hook
sums of \((2)\) and \((1,1)\).
\end{example}

\begin{example}[A shape-dependent mixed product]
\label{ex:shifted-parts-mixed-sharp}
The second shifted-part sum is
\[
 \mathcal R_2=\frac{N^3}{3}+\frac{3N^2}{2}-\frac{5N}{6}+2C_1.
\]
For example, \(\mathcal R_2((2))=9\) and
\(\mathcal R_2((1,1))=5\). Corollary~\ref{cor:shifted-parts-highest-differences}
gives \(D^6\mathcal R_2=240\mathbf1\).
For the product \(g=\mathcal H_1C_2\mathcal R_2\), the three
highest-weight parts evaluate to \(3/2\), \(1/2\) and \(1/3\),
respectively, under \(q_{2r}\mapsto\binom{2r}{r}\),
\(q_{2r+1}\mapsto0\).
The scalar-top identity therefore gives
\[
 D^{14}g=\frac{14!}{4}\mathbf1,\qquad D^{15}g=0.
\]
For every fixed \(\mu\), with \(m=|\mu|\), the polynomial factor
\(P_{\mu,g}\) in Theorem~\ref{thm:skew-exponential-polynomiality}
thus has degree exactly \(m+14\), with leading coefficient
\(f^\mu/(4m!)\). The unique top-degree term in the transfer formula
has \(\eta=\varnothing\) and \(j=14\).
\end{example}

\subsection{General hook polynomials with a fixed inner partition}

\begin{proof}[Proof of Theorem~\ref{thm:hook-polynomiality-intro}]
By Corollary~\ref{cor:hook-corner-polynomial}, \(\mathcal H_r\) has
corner weight at most \(2r+2\). The corner product bound in
Theorem~\ref{thm:corner-symbol} implies that \(D^{w+1}\) annihilates
every polynomial of weight at most \(w\).

To compute the highest difference, evaluate the component
\eqref{eq:hook-highest-weight} at \(q_d=a_d\) as in
Theorem~\ref{thm:corner-symbol}, using \(q_0=1\). Only even indices
remain, and
Vandermonde's identity gives
\[
 \frac{(2r)!}{2}\sum_{j=0}^{r+1}
                 \frac1{j!^2(r+1-j)!^2}
 =\frac{(2r)!}{2(r+1)!^2}\binom{2r+2}{r+1}
 =\beta_r^{\mathrm h}.
\]
The weight-\((2r+2)\) term is the same for \(\mathcal H_r\) and
\(\mathcal H_r-\mathcal H_r(\mu)\). Multiplicativity of the highest
symbol gives \eqref{eq:hook-highest-difference}. Its coefficient
is positive, so the degree is exact. Expansion in power sums gives
\eqref{eq:hook-symmetric-degree}, and Lemma~\ref{lem:skew-transfer}
gives the generating function.
\end{proof}

\begin{corollary}[Leading hook moments]
\label{cor:hook-leading-moments}
For every fixed \(\mu\) and fixed list \(r_1,\ldots,r_s\ge0\),
\[
 \mathbb E_{M_{\mu,n}}F^{\mathrm h}_{\boldsymbol r,\mu}
 \sim\left(\prod_i\beta_{r_i}^{\mathrm h}\right)
                                  n^{\sum_i(r_i+1)}.
\]
For each \(r\ge0\) and finite \(p\ge1\),
\[
 n^{-r-1}\bigl(\mathcal H_r(\lambda)-\mathcal H_r(\mu)\bigr)
       \longrightarrow\beta_r^{\mathrm h}
                   \quad\text{in }L^p(M_{\mu,n}).
\]
\end{corollary}
\begin{proof}
The first assertion follows from
Corollary~\ref{cor:expectation-asymptotic} and
\eqref{eq:hook-highest-difference}. Repeating a single index gives
convergence of every integer moment of the normalized statistic
to the corresponding power of \(\beta_r^{\mathrm h}\). Expanding
its centered \(2m\)th power shows convergence to zero in \(L^{2m}\).
Choose \(2m\ge p\) and use the probability-space \(L^p\) inequality.
\end{proof}

\begin{example}[Hook differences and changing old hooks]
\label{ex:hook-fixed-base-distinction}
For \(\mu=(1)\) and \(\lambda=(2)\),
\[
 \mathcal H_1(\lambda)-\mathcal H_1(\mu)=5-1=4.
\]
The old box's squared hook increases from \(1\) to \(4\), and the new
box has squared hook \(1\). The single box of \(\lambda/\mu\) has outer
hook length \(1\). Hence the fixed-base difference also includes the
changes in the hooks of boxes belonging to the inner partition.
For example \(\beta_0^{\mathrm h}=1\), \(\beta_1^{\mathrm h}=3/2\)
and \(\beta_2^{\mathrm h}=20/3\). Under any fixed base,
\(\mathcal H_1/n^2\to3/2\) and \(\mathcal H_2/n^3\to20/3\) in
every finite \(L^p\).
\end{example}

\begin{proposition}[The first hook power is not a \(D\)-polynomial]
\label{prop:odd-hook-counterexample}
Let \(T(\lambda)=\sum_{\square\in\lambda}h_\lambda(\square)\).
Then
\[
 A(x)=e^{-x-x^2/2}\sum_\lambda
                       \frac{x^{|\lambda|}}{H_\lambda}T(\lambda)
\]
is not a polynomial. In particular, \(T\) is not a \(D\)-polynomial.
\end{proposition}
\begin{proof}
For real \(x>0\), \(A(x)=\mathbb E_xT\). Cauchy's inequality on the
boxes of a diagram gives \(T(\lambda)^2\le|\lambda|\mathcal H_1(\lambda)\).
Applying it once more to the positive Poissonized measure yields
\[
 0\le A(x)\le
 \bigl(\mathbb E_x|\lambda|\,
                    \mathbb E_x\mathcal H_1\bigr)^{1/2}=O(x^3)
                    \qquad(x\longrightarrow+\infty).
\]
Here \(\mathbb E_x|\lambda|=x+x^2\), and
Theorem~\ref{thm:hook-polynomiality-intro} bounds the degree of
\(\mathbb E_x\mathcal H_1\) by four. These expectations are finite.
Every hook is at most \(|\lambda|\), so \(T\le|\lambda|^2\);
the size distribution has all moments. The same bound makes the
defining series for \(A\) absolutely convergent for every complex \(x\).

Put \(a_n=\sum_{\lambda\vdash n}T(\lambda)/H_\lambda\).
Direct evaluation of the diagrams of sizes zero through four gives
\[
 (a_0,a_1,a_2,a_3,a_4)
                  =\left(0,1,3,\frac{11}{3},\frac72\right).
\]
For example, at size four the two row or column shapes contribute
\(2(10/24)\), the shapes \((3,1)\) and \((2,1,1)\) contribute
\(2(8/8)\), and \((2,2)\) contributes \(8/12\).
Since \(e^{-x-x^2/2}=1-x+\frac13x^3+O(x^4)\),
\[
 [x^4]A(x)=\frac72-\frac{11}{3}+\frac13=\frac16.
\]
A polynomial with the growth bound \(O(x^3)\) has degree at most
three, contradicting this nonzero fourth coefficient.
If \(T\) were a \(D\)-polynomial, the empty-base transfer would make
\(A\) a polynomial. This proves the last assertion as well.
\end{proof}

\subsection{The factorial hook polynomial}
\label{sec:hook-degree-extremes}

\begin{proof}[Degree and extreme coefficients]
Expanding the factorial product in \(\mathcal J_k\) gives
\(\mathcal H_k\) plus hook sums of lower order.
Theorem~\ref{thm:hook-polynomiality-intro} yields
\(D^{2k+3}\mathcal J_k=0\), so the empty-base transfer gives
\[
 P_k(x)=\mathbb E_x\mathcal J_k
       =\sum_{j=0}^{2k+2}D^j\mathcal J_k(\varnothing)\frac{x^j}{j!}.
\]
Differentiation gives \eqref{eq:factorial-hook-all-coefficients}.
The lower-order hook sums have corner weights at most \(2k\). Hence
\eqref{eq:hook-highest-weight} and \eqref{eq:hook-highest-difference}
give \([x^{2k+2}]P_k=\beta_k^{\mathrm h}>0\). This agrees with the
Poisson transform of the Okada--Panova formula for the same statistic
under Plancherel measure \cite[Theorem~1]{Panova}. Indeed, under a
Poisson law of mean \(x^2\), the mean of \(N^{\downarrow(k+1)}\) is
\(x^{2k+2}\). Multiplying the leading coefficient by \(2k+2\) gives
the leading coefficient of \(P'_k\).

For \(k\ge1\), the next coefficient uses the term \(x^{d-1}L_d\)
in \eqref{clt:operator-leading}. The current means in
Lemma~\ref{clt:currents-lemma} give, for even \(d\ge2\),
\[
 \phi(L_d)=d\sum_{j=0}^{d/2-1}\binom{d-1}{j}
                         =d\,2^{d-2};
\]
for odd \(d\), \(\phi(L_d)=0\). In the moment of \(q_aq_b\), the term
of degree \(a+b-1\) comes from choosing one leading scalar and one
coefficient \(L_d\) by Lemma~\ref{clt:operator-moments-lemma}.
Here \(q_0=1\) contributes only its constant term, and \(q_1=0\).
Inserting these means in \eqref{eq:hook-highest-weight} gives
\[
 [x^{2k+1}]P_k=(2k)!\sum_{j=0}^{k}
                \frac{4^{k-j}}{j!^2(2k+1-2j)!}.
\]
After differentiation it becomes
\begin{align*}
 \sum_{j=0}^{k}\binom{2k+1}{2j}\binom{2j}{j}4^{k-j}
 &=\frac12[z^0](2+z+z^{-1})^{2k+1}\\
 &=\frac12\binom{4k+2}{2k+1}
  =\binom{4k+1}{2k}.
\end{align*}
Here \([z^0]\) extracts the constant Laurent coefficient, and
\(2+z+z^{-1}=z^{-1}(1+z)^2\) proves the middle equality.

Finally, for \(|\lambda|\le k\), every hook length is in
\(\{1,\ldots,k\}\), so \(\mathcal J_k(\lambda)=0\). At size \(k+1\), a
nonzero contribution requires a hook of length \(k+1\). These
are exactly the hook shapes \((k+1-a,1^a)\), \(0\le a\le k\),
with the unique contributing box \((1,1)\). Its weight is
\[
 \prod_{j=1}^k((k+1)^2-j^2)
       =\frac{k!(2k+1)!}{(k+1)!}.
\]
There are \(\binom ka\) standard tableaux of this hook shape:
choose which \(a\) entries among \(2,\ldots,k+1\) lie below the
first box, and sort the row and column. Their total count is
\(2^k\). By \(1/H_\lambda=f^\lambda/(k+1)!\),
\[
 [x^{k+1}]P_k(x)
 =\frac{2^k}{(k+1)!}\frac{k!(2k+1)!}{(k+1)!}.
\]
The normalization \(e^{-x-x^2/2}=1+O(x)\) leaves this first
nonzero coefficient unchanged. Differentiation proves the last
coefficient and the divisibility assertions. For \(k=0\) the
empty-base size identity gives \(P_0=x+x^2\) directly.
The complete finite formula and integrality are proved in
Theorem~\ref{thm:factorial-hook-complete} below.
\end{proof}

\subsection{A two-point constant term for hook statistics}
\label{sec:hook-two-point}

Under the specialization \(p_1=x,\ p_{j>1}=0\), the Poissonized
single-tableau weights are proportional to \(s_\lambda\). This is the
one-partition case of a Schur process with one free boundary, in which
the terminal partition is summed over all shapes. Its correlation
formulas are given in \cite[Theorem~3.3]{BR} and
\cite[Theorem~2.2]{BBNV}. For contour formulas for products of Macdonald difference operators
on the Pfaffian Schur partition function, see
\cite[Theorem~4.3]{Ghosal}.
We use a double constant term for the hook series. Its coefficient
extraction gives the finite binomial formula and derivative integrality
in Theorem~\ref{thm:factorial-hook-complete}.

\begin{lemma}[A two-point hook series]\label{lem:hook-two-point}
Let \(u\) be formal and put \(z=2\sinh(u/2)\). Since \(z=u+O(u^3)\),
the substitution has a unique formal inverse. Define
\[
 \mathcal G(x,z)=\mathbb E_x\bigl[\mathcal Q_\lambda(e^u)\mathcal Q_\lambda(e^{-u})\bigr],
\]
with \(\mathcal Q_\lambda\) as above. The series \(\mathcal G\) is even in \(z\).
Writing \(\operatorname{CT}_{a,b}=[a^0b^0]\) and setting
\[
 R_z(a)=\frac{1-a^2}{1-za-a^2},\qquad
 L=a+a^{-1}+b+b^{-1},
\]
we have
\begin{equation}\label{eq:hook-two-point-ct}
 \mathcal G(x,z)=\operatorname{CT}_{a,b} e^{xzL}R_z(a)R_z(b)
       \left(1+\frac{z^2ab}{(1+ab)^2}\right)
       \left(1+\frac{z^2(b/a)}{(1+b/a)^2}\right).
\end{equation}
Expand \(R_z(a)\) and \(R_z(b)\) first in \(z\), their remaining
denominators in nonnegative powers of \(a^2,b^2\), and the last two
denominators in nonnegative powers of \(ab,b/a\), respectively.
Every coefficient \(x^nz^d\) of this constant term is a finite sum.
\end{lemma}

\begin{proof}
We use the symmetric-function conventions of
Definition~\ref{def:symmetric-schur}. All vertex calculations can first
be made with finitely many modes. On a homogeneous input of degree
\(N\), a zero mode has weighted derivative degree at most \(N\), and
its creation degree must agree with it. Removing modes of index greater
than \(N\) therefore leaves that zero mode unchanged. Normal ordering
at a finite cutoff is a finite-variable formal operator identity; the
bounded-index argument below gives stabilization of every final
coefficient as the cutoff increases. Put \(r=e^{u/2}\) and
\(\zeta_j=r^j-r^{-j}\). In the vertex operator of Lascoux--Thibon
\cite[Section~3.2, Proposition~3.3]{LascouxThibon}, substitute
\(\vartheta_{\mathrm{LT}}=e^u\) and \(z_{\mathrm{LT}}=e^{-u/2}t\).
Their zero mode then acts on \(s_\lambda\) by
\[
 1+(\vartheta_{\mathrm{LT}}-1)(1-\vartheta_{\mathrm{LT}}^{-1})
       \sum_{\square\in\lambda}\vartheta_{\mathrm{LT}}^{c(\square)}
 =\mathcal Q_\lambda(e^u),
\]
where the equality follows from the one-box corner change.
In this normalization the vertex is
\[
 [t^0]V_u(t)s_\lambda=\mathcal Q_\lambda(e^u)s_\lambda,\qquad
 V_u(t)=
 \exp\left(\sum_{j\ge1}\frac{\zeta_jp_jt^j}{j}\right)
 \exp\left(\sum_{j\ge1}\zeta_jt^{-j}\partial_{p_j}\right).
\]
Here \(\partial_{p_j}\) denotes differentiation in \(p_j\).
By Littlewood's identity, the Schur sum is
\[
 G(p)=\sum_\lambda s_\lambda(p)
 =\exp\left(\sum_{\substack{j\ge1\\j\ {\rm odd}}}\frac{p_j}{j}
                     +\frac12\sum_{j\ge1}\frac{p_j^2}{j}\right).
\]
Each zero mode preserves the grading \(\deg p_j=j\). Before
normalization, the coefficient of \(x^N\) after specialization can
therefore be computed on the degree-\(N\) component of \(G(p)\),
which is a finite sum. By Lemma~\ref{lem:external-schur-identities} and
tableau branching, \(s_\lambda(x,0,\ldots)=f^\lambda x^N/N!\) for
\(N=|\lambda|\), so it equals \(x^{|\lambda|}/H_\lambda\) and
\(G(x,0,\ldots)=I(x)\). Consequently \(\mathcal G\) is obtained by taking the
double zero mode of \(V_u(t)V_{-u}(s)G(p)\), then specializing and
dividing by \(I(x)\).

Moving the derivative exponential of the first vertex past the
multiplication exponential of the second gives
\[
 K(s/t),\qquad
 K(y)=\exp\left(-\sum_{j\ge1}\frac{\zeta_j^2y^j}{j}\right)
     =\frac{(1-r^2y)(1-r^{-2}y)}{(1-y)^2}
     =1-\frac{z^2y}{(1-y)^2}.
\]
The derivatives then translate \(p_j\) by
\(\zeta_j(t^{-j}-s^{-j})\). The factors in the translated \(G\) involving
only \(t\) or only \(s\) are
\[
 R(t,z)=
 \exp\left(\sum_{\substack{j\ge1\\j\ {\rm odd}}}
                 \frac{\zeta_jt^{-j}}{j}
          +\frac12\sum_{j\ge1}\frac{\zeta_j^2t^{-2j}}{j}\right)
 =\frac{1-t^{-2}}{1-zt^{-1}-t^{-2}},
 \qquad R(s,-z).
\]
The cross term is \(K((ts)^{-1})\). Thus, after specialization
and division by \(I(x)\), the expression whose double constant term
is sought equals
\[
 e^{xz(t+t^{-1}-s-s^{-1})}
       R(t,z)R(s,-z)K(s/t)K((ts)^{-1}).
\]
Replace \(s\) by \(-s\), then set \(a=s^{-1}\), \(b=t^{-1}\).
Thus \(K(-s/t)\) becomes \(K(-b/a)\) and \(K(-(ts)^{-1})\)
becomes \(K(-ab)\).
Since \(K(-y)=1+z^2y/(1+y)^2\), this is
\eqref{eq:hook-two-point-ct}. The displayed expansion convention is
the one inherited from the normal-ordered product. For coefficientwise finiteness, the two kernel expansions have
nonnegative powers of \(b\), so the exponent supplied by the exponential
bounds their indices at fixed \(x^nz^d\); the remaining geometric-series
indices are then bounded by the constant-term conditions. For \(n\ge2\),
the cancellation in \eqref{eq:hook-kernel-cancellation} gives a shorter
finite extraction. Finally the definition of \(\mathcal G\) is unchanged
under \(u\mapsto-u\), proving evenness in \(z\).
\end{proof}

\subsection{All factorial-hook coefficients and their integrality}
\label{sec:hook-all-coefficients}

The two kernel factors in \eqref{eq:hook-two-point-ct} cancel most
of the summation indices. We first evaluate the one-variable
coefficient that remains. For nonnegative integers \(p,m,h\),
the coefficient has the explicit binomial evaluation
\begin{equation}\label{eq:hook-one-variable-explicit}
 \mathfrak f_{p,m}(h)
 :=[y^h]\frac{(1+y)^p}{(1-y)^m}
 =\sum_{i=0}^{\min(p,h)}\binom pi
                         \binom{m+h-i-1}{h-i}.
\end{equation}
The binomial coefficient on the right uses the polynomial convention;
in particular \(\binom{-1}{0}=1\). Thus the same formula includes
\(m=0\), when it reduces to \(\binom ph\). Put
\(\mathfrak f_{p,m}(h)=0\) if \(h\notin\mathbb Z_{\ge0}\).
The central specialization admits the product evaluation
\begin{equation}\label{eq:hook-central-binomial}
 \mathfrak f_{m+2h,m}(h)
   =4^h\binom{h+(m-1)/2}{h}
   =\frac{2^h}{h!}\prod_{i=1}^h(m+2i-1).
\end{equation}
The empty product is one. For example,
\(\mathfrak f_{4,0}(2)=6\),
\(\mathfrak f_{3,1}(1)=4\) and
\(\mathfrak f_{6,2}(2)=30\).

For all integers \(p,h\ge0\) define
\[
 t_{p,h}=\begin{cases}
  4^h\displaystyle\binom{(p-1)/2}{h},&2h\le p,\\
  0,&2h>p,
 \end{cases}
 \qquad
 \operatorname{Cat}_i=\frac1{i+1}\binom{2i}{i}.
\]
Also put \(t_{p,h}=0\) for negative \(h\).
These are nonnegative integers, by
\eqref{eq:hook-central-binomial} with \(m=p-2h\).
For integers \(n\ge0\) and \(0\le2h\le n\), set
\begin{align}\label{eq:hook-short-binomial-arrays}
 \mathfrak a_{n,h}
  &=\sum_{p=0}^n\binom np\sum_{i=0}^h
                            t_{p,i}t_{n-p,h-i},\notag\\
 \mathfrak e_{n,h}
  &=\sum_{i=0}^h\binom n{2i}\operatorname{Cat}_i
                            t_{n-2i,h-i}.
\end{align}
Both arrays are zero if \(h<0\) or \(2h>n\).
For \(n\ge2\) and an even integer \(d\ge0\), define
\begin{equation}\label{eq:hook-finite-binomial-sum}
 \mathfrak b_{n,d}
  =\mathfrak a_{n,n-d/2}
    +\frac12\mathfrak a_{n,n-d/2+1}
    +2\mathfrak e_{n-2,n-d/2}.
\end{equation}
For odd \(d\), put \(\mathfrak b_{n,d}=0\).
Although the middle term has a factor \(1/2\),
\(\mathfrak a_{n,h}\) is even for \(n\ge1\), so every
\(\mathfrak b_{n,d}\) is an integer. Formula
\eqref{eq:hook-finite-binomial-sum} uses only the two sums in
\eqref{eq:hook-short-binomial-arrays}.

\begin{theorem}[Complete factorial-hook coefficients]
\label{thm:factorial-hook-complete}
One has \(P_0=x+x^2\). For \(k\ge1\),
\begin{equation}\label{eq:factorial-hook-complete-coefficients}
 P'_k(x)=\frac{(k!)^2}{2}
 \sum_{n=k+1}^{2k+2}\frac{x^{n-1}}{(n-1)!}
 \sum_{j=\lfloor(n-1)/2\rfloor}^k
       (-1)^{k-j}\binom{2j}{j}\mathfrak b_{n,2j+2}.
\end{equation}
Here \(\mathfrak b\) is the explicit binomial expression
\eqref{eq:hook-finite-binomial-sum}. Since \(P_k(0)=0\), the formula
determines every coefficient of \(P_k\). Moreover,
\begin{equation}\label{eq:factorial-hook-derivative-integral}
 P'_k(x)\in x^k\mathbb Z[x],\qquad \deg P'_k=2k+1.
\end{equation}
\end{theorem}

\begin{proof}
Expansion of the numerator and denominator proves
\eqref{eq:hook-one-variable-explicit}. For the central value, use
the formal change of variable \(w=y/(1+y)^2\). It is invertible
at zero and satisfies
\[
 \frac{\dd y}{y}=\frac{1+y}{1-y}\frac{\dd w}{w},\qquad
 \frac{1-y}{1+y}=\sqrt{1-4w}.
\]
The residue giving the coefficient consequently becomes
\[
 [y^h]\frac{(1+y)^{m+2h}}{(1-y)^m}
   =[w^h](1-4w)^{-(m+1)/2}
   =4^h\binom{h+(m-1)/2}{h}.
\]
This also proves the product form in
\eqref{eq:hook-central-binomial}.

Write \(A=a+a^{-1}\), \(B=b+b^{-1}\), and \(L=A+B\).
An elementary rational identity gives
\begin{equation}\label{eq:hook-kernel-cancellation}
 \left(1+\frac{z^2ab}{(1+ab)^2}\right)
 \left(1+\frac{z^2(b/a)}{(1+b/a)^2}\right)
       =1+\frac{z^2(AB+4)+z^4}{L^2}.
\end{equation}
For instance, multiplying the sum of the two fractions by \(L^2\)
gives \(AB+4\), and multiplying their product by \(L^2\) gives one.
After multiplication by \(L^n\), \(n\ge2\), the right side is a
Laurent polynomial. Thus the expansion convention in
Lemma~\ref{lem:hook-two-point} is preserved, and the infinite
kernel expansions disappear before any constant term is taken.

Put
\[
 T_p(z)=\operatorname{CT}_a A^pR_z(a)
       =\sum_{h=0}^{\lfloor p/2\rfloor}t_{p,h}z^{p-2h}.
\]
Indeed, the coefficient of \(z^m\) in this constant term is
\(\mathfrak f_{p,m}((p-m)/2)\), evaluated above. Define
\[
 \mathfrak A_n(z)=\sum_{p=0}^n\binom np T_p(z)T_{n-p}(z),\qquad
 \mathfrak E_n(z)=\sum_{i=0}^{\lfloor n/2\rfloor}
                 \binom n{2i}\operatorname{Cat}_iT_{n-2i}(z).
\]
Their coefficients of \(z^{n-2h}\) are respectively
\(\mathfrak a_{n,h}\) and \(\mathfrak e_{n,h}\).
Pascal's identity in the explicit formula for \(t_{p,h}\) gives
\[
 T_{p+2}(z)=(z^2+4)T_p(z)
             -\begin{cases}
 2\operatorname{Cat}_{p/2},&p\text{ even},\\
 0,&p\text{ odd}.
 \end{cases}
\]
For the constant coefficient when \(p=2r\), the correction is
\(-2\operatorname{Cat}_r\); all other coefficients satisfy Pascal's
identity without a correction. If
\(\mathfrak D_n=\sum_{p=0}^n\binom np T_{p+1}T_{n-p+1}\),
expanding \(\mathfrak A_{n+2}\) by the two Pascal identities for
\(\binom{n+2}{p}\) yields
\[
 \mathfrak A_{n+2}
    =2(z^2+4)\mathfrak A_n+2\mathfrak D_n-4\mathfrak E_n.
\]
It follows from \eqref{eq:hook-kernel-cancellation} that, for \(n\ge2\),
\begin{align*}
 n![x^n]\mathcal G(x,z)
 &=z^n\bigl(\mathfrak A_n
       +z^2(\mathfrak D_{n-2}+4\mathfrak A_{n-2})
       +z^4\mathfrak A_{n-2}\bigr)\\
 &=z^n\bigl((1+z^2/2)\mathfrak A_n
                                      +2z^2\mathfrak E_{n-2}\bigr).
\end{align*}
Consequently
\(n![x^nz^d]\mathcal G=\mathfrak b_{n,d}\).
This proves the short binomial extraction. The symmetry
\(p\leftrightarrow n-p\) in \(\mathfrak A_n\) pairs all terms
unless \(n\) is even and \(p=n/2\); in that case
\(\binom n{n/2}\) is even. Thus \(\mathfrak a_{n,h}\) is even
for every \(n\ge1\), as asserted.

For an indeterminate \(h\), put
\[
 F_j(h)=\prod_{i=0}^{j-1}(h^2-i^2),\qquad
 J_k(h)=\prod_{i=1}^k(h^2-i^2).
\]
The polynomial identity
\begin{equation}\label{eq:hook-central-factorial-inversion}
 J_k(h)=(k!)^2\sum_{j=0}^k\frac{(-1)^{k-j}}{(j!)^2}F_j(h)
\end{equation}
follows on multiplication by \(h^2\) from
\(h^2F_j=F_{j+1}+j^2F_j\), by telescoping.
Lemma~\ref{lem:content-central-factorial-identities} and the hook
boundary identity give
\[
 \mathbb E_x\sum_{\square\in\lambda}
        \cosh(uh_\lambda(\square))
   =\frac{\mathcal G(x,z)-1}{2z^2},\qquad
 \cosh(uh)=\sum_{j\ge0}\frac{F_j(h)}{(2j)!}z^{2j},
 \quad z=2\sinh(u/2).
\]
Therefore, for \(n\ge2\),
\[
 [x^{n-1}]P'_k(x)
   =\frac{(k!)^2}{2(n-1)!}
      \sum_{j=0}^k(-1)^{k-j}\binom{2j}{j}
                                         \mathfrak b_{n,2j+2}.
\]
The definition of \(\mathfrak b\) makes it zero for \(n>2j+2\).
Thus \(j\) starts at \(\lfloor(n-1)/2\rfloor\).
The initial vanishing and degree established in
Section~\ref{sec:hook-degree-extremes} restrict \(n\) to
\(k+1,\ldots,2k+2\), giving
\eqref{eq:factorial-hook-complete-coefficients}.

For integrality it is convenient to write a summand before factoring
out \((k!)^2\), as
\[
 \left(\frac{k!}{j!}\right)^2
       \frac{(2j)!}{2(n-1)!}\mathfrak b_{n,2j+2}.
\]
The first factor is an integer. If \(n\le2j\), the factorial ratio
\((2j)!/(n-1)!\) is an even integer. If \(n=2j+1\), the formula
\eqref{eq:hook-finite-binomial-sum} reduces to
\(\mathfrak b_{n,n+1}=\mathfrak a_{n,(n-1)/2}\), which is even.
Finally, if \(n=2j+2\), only the constant terms of the two
\(R_z\) factors contribute, and
\[
 \mathfrak b_{2j+2,2j+2}
     =\operatorname{CT}_{a,b}L^{2j+2}
     =\binom{2j+2}{j+1}^{\!2}.
\]
The last equality follows by writing
\(L=(a+b)(1+(ab)^{-1})\). The required integer is then
\[
 \frac1{2(2j+1)}\binom{2j+2}{j+1}^{\!2}
       =2\operatorname{Cat}_j\binom{2j+1}{j}.
\]
These cases exhaust all nonzero terms. They prove
\eqref{eq:factorial-hook-derivative-integral}; the nonzero leading
coefficient was computed in Section~\ref{sec:hook-degree-extremes}.
\end{proof}

\begin{example}[The first factorial hook polynomials]
\label{ex:factorial-hook-small}
The first five polynomials are
\[
 P_0=x+x^2,\qquad
 P_1=3x^2+\frac{10}{3}x^3+\frac32x^4,
\]
\[
 P_2=\frac{80}{3}x^3+40x^4+\frac{126}{5}x^5+\frac{20}{3}x^6.
\]
The next two are
\begin{align*}
 P_3={}&420x^4+\frac{3906}{5}x^5+616x^6
                +\frac{1716}{7}x^7+\frac{175}{4}x^8,\\
 P_4={}&\frac{48384}{5}x^5+20832x^6+\frac{135168}{7}x^7
                  +9696x^8+\frac{24310}{9}x^9+\frac{1764}{5}x^{10}.
\end{align*}
All coefficients are also given by the finite binomial formula
\eqref{eq:factorial-hook-complete-coefficients}; in particular,
Theorem~\ref{thm:factorial-hook-complete} proves that every coefficient
of \(P'_k\) is an integer.
For \(k=1\), \(\mathcal J_1=\mathcal H_1-\mathcal H_0=C_2+|\lambda|^2-|\lambda|\),
so this example also follows from the content and size identities.
For a fixed nonempty \(\mu\), Lemma~\ref{lem:skew-transfer}
applies to \(\mathcal J_k\) with degree bound \(2k+2\) and gives a polynomial
of degree at most \(|\mu|+2k+2\) multiplying the same involution
exponential. The coefficients in \eqref{eq:factorial-hook-three-coefficients}
are those for the empty base.
\end{example}

\section{Limit shape and Gaussian fluctuations from the difference operator}
\label{sec:limit-shape}

We prove Theorems~\ref{thm:skew-shape-intro}
and~\ref{thm:skew-asymptotics-intro} from the local operator \(D\).
The highest-difference formulas and the exact skew transfer give the
content laws and uniform limit shape. The polynomial conjugation of
Section~\ref{sec:exact-moments}, followed by coefficient extraction,
gives the joint Gaussian limits directly from any fixed inner shape.

Fix \(\mu\), put \(m=|\mu|\), \(n=|\lambda/\mu|\ge1\), and
\(N=n+m\). For a boundary profile, its charge density is half its
excess over the absolute-value profile. Write
\[
 \begin{gathered}
 \sigma_\lambda(x)=\frac{\lambda(x)-|x|}{2},\qquad
 \sigma_\Omega(x)=\frac{\Omega(x)-|x|}{2},\\
 \bar\lambda_N(x)=N^{-1/2}\lambda(\sqrt N x),\qquad
 \sigma_N(x)=\frac{\bar\lambda_N(x)-|x|}{2}.
 \end{gathered}
\]
For \(t>0\), define the outer fluctuation by
\begin{equation}\label{eq:outer-profile-fluctuation}
 \Delta_{\lambda,t}^{Y,\mathrm{out}}(x)
 =\frac{\sqrt t}{2}
   \left(t^{-1/2}\lambda(\sqrt t x)-\Omega(x)\right).
\end{equation}
The genuine skew profile and its fluctuation are
\begin{equation}\label{def:skew-profile}
 \omega_{\lambda/\mu}(x)=\lambda(x)-\mu(x)+|x|,\quad
 \overline\omega_{\lambda/\mu,n}(x)
       =n^{-1/2}\omega_{\lambda/\mu}(\sqrt n x),\quad
 \Delta_{\lambda/\mu,n}^{Y}
       =\frac{\sqrt n}{2}(\overline\omega_{\lambda/\mu,n}-\Omega).
\end{equation}
These are the conventions of Section~\ref{sec:introduction}.
The symbol \(\delta_a\) denotes a unit point mass at \(a\);
distributional derivatives are defined by integration by parts against
smooth compactly supported tests.

\subsection{Content moments and the uniform limit shape}

\begin{lemma}[Content laws from highest differences]
\label{lem:operator-content-lln}
For every fixed \(\mu\), integer \(j\ge0\), and \(1\le p<\infty\),
\[
 n^{-j/2-1}C_j(\lambda/\mu)\longrightarrow\gamma_j
                    \quad\text{in }L^p(M_{\mu,n}).
\]
The assertion is unchanged with \(N\) in place of \(n\).
\end{lemma}
\begin{proof}
For \(j=2r\), Theorem~\ref{thm:second-moment}, whose proof uses the
highest \(D\)-differences and Corollary~\ref{cor:expectation-asymptotic},
gives, for every fixed positive integer \(b\),
\[
 \mathbb E_{M_{\mu,n}}(n^{-r-1}C_{2r}(\lambda/\mu))^b
       \longrightarrow c_r^b,\qquad
 c_r=\frac{(2r)!}{(r+1)!^2}=\gamma_{2r}.
\]
Expanding a centered even power proves convergence in every finite
\(L^p\). For odd \(j\), Lemma~\ref{lem:content-antiderivative}
expresses \(C_j-C_j(\mu)\) as a signed corner sum with defining
polynomial of odd degree \(j+2\). For every positive integer \(a\),
Corollary~\ref{cor:corner-parity-degree} therefore gives
\[
 D^{2a(j+1)+1}(C_j-C_j(\mu))^{2a}=0.
\]
Corollary~\ref{cor:expectation-asymptotic} bounds its expectation by
\(O(n^{a(j+1)})\); after normalization the \(2a\)-th absolute moment
is \(O(n^{-a})\). Taking \(2a\ge p\) proves the assertion with
\(\gamma_j=0\). Finally, \(N/n\to1\).
\end{proof}

\begin{lemma}[Tent decomposition]\label{lem:tent}
Let \(T(u)=(1-|u|)_+\), where \(r_+=\max\{r,0\}\). Then
\begin{equation}\label{eq:tent-decomp}
 \sigma_\lambda(x)=\sum_{\square\in\lambda}T(x-c(\square)).
\end{equation}
For every integer \(k\ge0\),
\begin{equation}\label{eq:moment-charge-content}
 \int_{\mathbb R}x^k\sigma_\lambda(x)\,dx
 =\sum_{j=0}^{\lfloor k/2\rfloor}
   \binom{k}{2j}\frac{2}{(2j+1)(2j+2)}\,C_{k-2j}(\lambda).
\end{equation}
\end{lemma}
\begin{proof}
Adding a box of content \(c\) increases the Russian boundary by \(2T(x-c)\).
Induction from the empty diagram proves the first identity. For the second,
expand \((c+u)^k\) and use
\[
 \int_{-1}^1u^{2j}T(u)\,du=\frac{2}{(2j+1)(2j+2)},\qquad
 \int_{-1}^1u^{2j+1}T(u)\,du=0.
\]
In particular, the coefficient of \(C_k\) is \(1\).
\end{proof}

\begin{lemma}[Moments of the limit profile]\label{lem:omega-moments}
For every integer \(r\ge0\),
\begin{equation}\label{eq:Omega-charge-moments}
 \int_{\mathbb R}x^{2r}\sigma_\Omega(x)\,dx
 =\frac{(2r)!}{(r+1)!^2},\qquad
 \int_{\mathbb R}x^{2r+1}\sigma_\Omega(x)\,dx=0.
\end{equation}
\end{lemma}
\begin{proof}
The function \(\sigma_\Omega\) is even and supported on \([-2,2]\).
On \((-2,0)\cup(0,2)\), its second derivative is
\(1/(\pi\sqrt{4-x^2})\), and its first derivative has a jump of \(-1\)
at \(0\). Integrating by parts twice against \(x^{k+2}\),
whose value at \(0\) is zero, gives
\[
 \int x^k\sigma_\Omega(x)\,dx
 =\frac{1}{(k+1)(k+2)}
   \int_{-2}^{2}\frac{x^{k+2}}{\pi\sqrt{4-x^2}}\,dx.
\]
There are no endpoint boundary terms, since \(\sigma_\Omega\) and its first
derivative vanish at \(\pm2\). Substitution \(x=2\cos\theta\) gives
\eqref{eq:Omega-charge-moments}.
\end{proof}

\begin{lemma}[Uniform convergence from charge moments]
\label{lem:charge-compactness}
Let \(s_j\) be nonnegative \(1\)-Lipschitz functions on \(\R\), with
\(\int_\R s_j=1\). Suppose that, for every nonnegative integer \(k\),
\[
 \int_\R x^k s_j(x)\,\dd x
       \longrightarrow\int_\R x^k\sigma_\Omega(x)\,\dd x.
\]
Then \(\|s_j-\sigma_\Omega\|_\infty\to0\).
\end{lemma}
\begin{proof}
If $a=s_j(x)$, the Lipschitz bound gives
$s_j(y)\ge(a-|y-x|)_+$, hence $a^2\le1$. Integrating the outward
half of this triangle when $|x|>R$ gives
\begin{equation}\label{eq:tail-height}
 \sup_{|x|>R}s_j(x)
 \le\left(2\int_{|y|>R}s_j(y)\,\dd y\right)^{1/2}
 \le\frac{\sqrt{2\int_\R y^2s_j(y)\,\dd y}}R.
\end{equation}
The uniformly bounded second moments and Arzel\`a--Ascoli therefore
give, from every subsequence, a further subsequence converging uniformly
on $\R$ to a continuous nonnegative function $s$. For $2q>k$,
\[
 \int_{|x|>R}|x|^k s_j(x)\,\dd x
       \le R^{k-2q}\int_\R x^{2q}s_j(x)\,\dd x.
\]
The right side tends uniformly to zero as $R\to\infty$; Fatou gives
the same tail bound for $s$. Every moment of $s$, including its mass,
thus equals that of $\sigma_\Omega$. Its even moments are at most
$4^q$, so
$\int_{|x|>2+\varepsilon}s(x)\,\dd x
 \le4^q/(2+\varepsilon)^{2q}\to0$.
Both measures are supported on $[-2,2]$. Weierstrass approximation
gives their equality, and continuity gives $s=\sigma_\Omega$.
Uniqueness of the subsequential limit proves the assertion.
\end{proof}

\begin{theorem}[Uniform shape from a fixed skew base]
\label{thm:skew-limit-shape}\label{cor:genuine-skew-shape}
For every fixed \(\mu\), under \(M_{\mu,n}\),
\[
 \|\bar\lambda_N-\Omega\|_\infty\longrightarrow0,\qquad
 \|\overline\omega_{\lambda/\mu,n}-\Omega\|_\infty\longrightarrow0
                  \quad\text{in probability}.
\]
The outer assertion also holds with \(n\) in place of \(N\).
\end{theorem}
\begin{proof}
The tent identity and the boundary slopes show that each \(\sigma_N\) is nonnegative,
\(1\)-Lipschitz, and has mass one. By
\eqref{eq:moment-charge-content} and Lemma~\ref{lem:operator-content-lln},
for every integer \(k\ge0\),
\[
 \int x^k\sigma_N(x)\,\dd x
 =N^{-k/2-1}\sum_{j=0}^{\lfloor k/2\rfloor}
       \binom{k}{2j}\frac{2C_{k-2j}(\lambda)}{(2j+1)(2j+2)}
 \longrightarrow\int x^k\sigma_\Omega(x)\,\dd x
\]
in probability. Indeed, the \(j=0\) term has limit \(\gamma_k\),
the fixed inner contribution vanishes, and every \(j\ge1\) term
has smaller order. Lemma~\ref{lem:omega-moments} identifies the limit.
Realize the partitions on a common product probability space. From
any subsequence, choose a further subsequence on which all these
countably many moment convergences hold almost surely.
Lemma~\ref{lem:charge-compactness} then gives uniform convergence
on each such sample path. The subsequence criterion proves
\(\|\sigma_N-\sigma_\Omega\|_\infty\to0\) in probability, hence the
outer assertion because \(\bar\lambda_N-\Omega=2(\sigma_N-\sigma_\Omega)\).

For \(t_N=N/n\), the charge of the outer profile at the increment
scale is \(\sqrt{t_N}\sigma_N(x/\sqrt{t_N})\). Its uniform difference
from \(\sigma_\Omega\) is bounded by
\[
 \sqrt{t_N}\|\sigma_N-\sigma_\Omega\|_\infty+
 \|\sqrt{t_N}\sigma_\Omega(\cdot/\sqrt{t_N})-\sigma_\Omega\|_\infty.
\]
The second term tends to zero since \(t_N\to1\) and
\(\sigma_\Omega\) is continuous with compact support. Finally,
\[
 \overline\omega_{\lambda/\mu,n}(x)
     -n^{-1/2}\lambda(\sqrt n x)=-2n^{-1/2}\sigma_\mu(\sqrt n x),
\]
whose uniform norm tends to zero. This proves the skew assertion.
Together with Lemma~\ref{lem:operator-content-lln}, it proves
Theorem~\ref{thm:skew-shape-intro}.
\end{proof}

\subsection{Coefficient extraction and its local limit}
\label{clt:transfer-subsection}

The functionals \(\phi_\eta\) and \(\phi\) in \eqref{clt:phi} are unital.
For \(n\ge1\), write \(\varrho=\sqrt n\) and define
\begin{equation}\label{clt:coefficient-functional}
 \mathcal T_n(P)=\frac{n![z^n]P(z)I(z)}{I_n}\quad(P\in\mathbb C[z]),
 \qquad J_n(c)=n![z^n]e^{cz+z^2/2}.
\end{equation}
Here \([z^n]\) extracts the coefficient of \(z^n\).
The translated variable is \(h=z-\varrho\). On \(\mathbb C[h]\) define
\[
 \tau_n(P(h))=\mathcal T_n(P(z-\varrho)),\qquad
 \tau(h^a)=a![t^a]e^{-t/2-t^2/4}\quad(a\in\mathbb Z_{\ge0}),
\]
extending the latter definition linearly from the monomials. The value
\(\tau(h^2)-\tau(h)^2=-1/2\) shows that \(\tau\) is a formal moment
functional. When \(C\) is a row-finite operator, the tensor-product
notation means \((\tau\otimes\phi)(h^aC)=\tau(h^a)\phi(C)\), extended
linearly to operator-valued polynomials in \(h\).

The coefficient estimate below is related to the involution
de-Poissonization coefficients in \cite[Section~5.2]{Bornemann}.
For \(0\le a\le n\), the functional in
\eqref{clt:coefficient-functional} satisfies
\[
 \mathcal T_n(z^a)=\fall{n}{a}\frac{I_{n-a}}{I_n},
 \qquad
 \mathcal T_n(e^{tz})=\frac{J_n(1+t)}{J_n(1)}.
\]
For \(a>n\), \(\mathcal T_n(z^a)=0\). The exponential identity can be
read as an identity of Taylor series in \(t\), with \(\mathcal T_n\)
applied to each coefficient.

\begin{lemma}[Uniform analytic estimate]\label{clt:analytic-coefficient-estimate}
Uniformly for \(c\) in each compact set \(\mathcal K\subset\{c\in\C:\operatorname{Re}c>0\}\),
\begin{equation}\label{clt:Jn-asymptotic}
 J_n(c)=\frac{\varrho^n e^{-n/2}}{\sqrt2}
       e^{c\varrho-c^2/4}\bigl(1+O(\varrho^{-1})\bigr).
\end{equation}
The relative error is holomorphic in \(c\). For every fixed nonnegative integer \(a\),
\begin{equation}\label{clt:eta-moment-limit}
 \tau_n(h^a)=\tau(h^a)+O_a(n^{-1/2}),
 \qquad
 \sum_{a\ge0}\tau(h^a)\frac{t^a}{a!}
       =e^{-t/2-t^2/4}.
\end{equation}
Each fixed moment on the left is uniformly bounded in \(n\).
Explicitly, its limiting formal moments are
\[
 \tau(h^a)=a!\sum_{j=0}^{\lfloor a/2\rfloor}
      \frac{(-1/2)^{a-2j}(-1/4)^j}{(a-2j)!\,j!}.
\]
\end{lemma}

\begin{proof}
For a standard real Gaussian variable \(Z\), comparison of exponential
moments gives \(J_n(c)=\E(c+Z)^n\).  A change of variables for real \(c\),
followed by analytic continuation, yields
\begin{equation}\label{clt:Jn-gaussian-integral}
 J_n(c)=\frac{e^{-c^2/2}}{\sqrt{2\pi}}
 \int_0^\infty t^n e^{-t^2/2}
       \bigl(e^{ct}+(-1)^n e^{-ct}\bigr)\,dt.
\end{equation}
Both sides are entire in \(c\). Write
\[
 J_n^{+}(c)=\int_0^\infty t^n e^{-t^2/2+ct}\,dt.
\]
With \(t=\varrho+u\), this becomes
\[
 J_n^{+}(c)=\varrho^n e^{-n/2+c\varrho}
  \int_{-\varrho}^{\infty}e^{h_\varrho(u)+cu}\,du,
 \qquad
 h_\varrho(u)=\varrho^2\log(1+u/\varrho)-\varrho u-u^2/2.
\]
Since \(h_\varrho(0)=h_\varrho'(0)=0\) and
\[
 h_\varrho''(u)=-\frac{\varrho^2}{(\varrho+u)^2}-1\le-1,
\]
we have \(h_\varrho(u)\le-u^2/2\) on \((-\varrho,\infty)\).  On
\(|u|\le\varrho^{1/6}\), Taylor's formula also gives
\[
 h_\varrho(u)=-u^2+O(|u|^3/\varrho).
\]
For \(c\) in a fixed compact set \(\mathcal K\subset\C\), the integral error on this
central interval after replacing \(h_\varrho(u)\) by \(-u^2\) is bounded by
\[
 \frac C\varrho\int_\R |u|^3e^{-u^2+C|u|}\,du=O(\varrho^{-1}).
\]
The omitted tails are \(O_{\mathcal K}(e^{-c_0\varrho^{1/3}})\) for some
\(c_0>0\), by \(h_\varrho(u)\le-u^2/2\) and boundedness of \(c\).
This stretched-exponential bound is \(o(\varrho^{-1})\).
Therefore, uniformly on such a compact set,
\[
 J_n^{+}(c)=\varrho^n e^{-n/2+c\varrho}\sqrt\pi\,e^{c^2/4}
          \bigl(1+O(\varrho^{-1})\bigr).
\]
If additionally \(\operatorname{Re}c\ge c_{\min}>0\), the contribution of
\(J_n^{+}(-c)\) relative to \(J_n^{+}(c)\) is \(O(e^{-2c_{\min}\varrho})\).
Equation~\eqref{clt:Jn-gaussian-integral} now proves
\eqref{clt:Jn-asymptotic}.  Holomorphy of its normalized error follows
from its explicit definition as an entire function divided by a nonvanishing
exponential.

Fix a closed disk \(|t|\le\epsilon<1\), and apply
\eqref{clt:Jn-asymptotic} with \(c=1+t\) and \(c=1\).  The exact identity
for \(\mathcal T_n(e^{tz})\) gives
\[
 \tau_n(e^{th})
 =e^{-t\varrho}\frac{J_n(1+t)}{J_n(1)}
 =e^{-t/2-t^2/4}\bigl(1+O(\varrho^{-1})\bigr)
\]
uniformly on this disk.  Cauchy's coefficient formula implies convergence
of every Taylor coefficient, which is precisely
\eqref{clt:eta-moment-limit}.  It also proves boundedness of each fixed
moment. Expanding the limiting exponential gives the displayed formula
for \(\tau(h^a)\); no positivity assertion is needed.
\end{proof}

\subsection{Fixed-size moments from polynomial conjugation}
\label{clt:fixed-size-subsection}

Let \(e_r=\mathbf1_{\{r\text{ odd}\}}\), and let
\(\Xi_r\sim\mathcal N(e_r,2r)\), \(r\ge2\), be independent. Put
\(V_d=\sum_{r\ge2}b_{dr}\Xi_r\), as in
Theorem~\ref{thm:skew-asymptotics-intro}. We shall identify them by the Gaussian-current identity.

\begin{proposition}[A local fixed-size expansion]
\label{prop:operator-fixed-moment-rate}
For a fixed finite family \(\mathcal D\subseteq\{2,3,\ldots\}\), set
\[
 \mathsf Z_{d,N}(z)=\frac{\mathsf M_d(z)-a_dN^{d/2}I}{N^{(d-1)/2}},\qquad
 Y_{d,N}=\frac{q_d-a_dN^{d/2}}{N^{(d-1)/2}}.
\]
Write \(\mathsf Z_N=(\mathsf Z_{d,N})_{d\in\mathcal D}\) and
\(Y_N=(Y_{d,N})_{d\in\mathcal D}\). For every fixed polynomial \(P\),
\[
 \mathbb E_{0,N}P((Y_{d,N})_{d\in\mathcal D})
       =\mathbb EP((V_d)_{d\in\mathcal D})+O_P(N^{-1/2}).
\]
For bounded degree, the error is uniform with a constant times the
sum of the absolute values of the coefficients of \(P\).
\end{proposition}
\begin{proof}
Put \(\varrho=\sqrt N\). The exact moment identity of
Lemma~\ref{clt:operator-moments-lemma} gives
\[
 \mathbb E_{0,N}P(Y_N)
 =\mathcal T_N\bigl[\phi(P(\mathsf Z_N(z)))\bigr].
\]
The components of \(\mathsf Z_N\) commute, being simultaneous conjugates of
multiplication operators. The leading conjugation formula gives
\begin{equation}\label{clt:operator-local-expansion}
 \mathsf Z_{d,N}(\varrho+h)=d a_dhI+L_d+R_{d,N}(h),
\end{equation}
where the remainder has bounded degree in \(h\), and its coefficients
are \(O(\varrho^{-1})\) in a fixed finite operator span. Indeed,
expand \((\varrho+h)^d-\varrho^d\) and use that the remaining terms
of \(\mathsf M_d(z)\) have degree at most \(d-2\).
Lemma~\ref{clt:analytic-coefficient-estimate} now yields
\[
 \mathbb E_{0,N}P(Y_N)
 =(\tau\otimes\phi)P((d a_dhI+L_d)_{d\in\mathcal D})
                           +O_P(\varrho^{-1}).
\]
All products involve finitely many fixed operator words and scalar
moments. The binomial definitions give \(b_{d1}=d a_d/2\).
To identify the leading functional, put
\(s_r=\sum_d b_{dr}t_d\) and \(v=\sum_d d a_dt_d=2s_1\).
The current identity gives the formal exponential series
\begin{align}\label{clt:mode-one-cancellation}
 (\tau\otimes\phi)e^{\sum_dt_d(d a_dhI+L_d)}
 &=e^{-v/2-v^2/4}
       e^{\sum_{r\ge1}e_rs_r+\sum_{r\ge1}rs_r^2}\notag\\
 &=e^{\sum_{r\ge2}e_rs_r+\sum_{r\ge2}rs_r^2}.
\end{align}
Both terms of the size mode cancel. The last expression is the
Gaussian moment generating function of \((V_d)_d\). Comparing
coefficients proves the assertion; the finite number of monomials
of bounded degree gives its uniform form.
\end{proof}

\begin{corollary}[Quantitative corner moments from the skew \(D\)-transfer]
\label{thm:skew-corner-clt}
Fix \(\mu\) and a finite set \(\mathcal D\subseteq\{2,3,\ldots\}\).
Under \(M_{\mu,n}\), the vector
\[
 V_N^\mu=\left(\frac{q_d-a_dN^{d/2}}{N^{(d-1)/2}}\right)_{d\in\mathcal D}
\]
converges in distribution and in all mixed moments to
\(V=(V_d)_{d\in\mathcal D}\). For every fixed polynomial \(P\),
\[
 \mathbb E_{M_{\mu,n}}P(V_N^\mu)=\mathbb EP(V)+O_{\mu,P}(n^{-1/2}).
\]
For bounded degree the error is bounded by a constant, depending only
on \(\mu,\mathcal D\) and that degree, times the sum of absolute
coefficients of \(P\). The same rate holds after separate exact
expectation centering, after replacing \(N\) by \(n\), or after
subtracting \(q_d(\mu)\) from \(q_d(\lambda)\).
\end{corollary}
\begin{proof}
The product-commutator argument in
Lemma~\ref{clt:operator-moments-lemma} works in every row. Thus, for
a polynomial \(g\) in a fixed finite family of corner observables,
\[
 (e^{zD}g)(\eta)=\phi_\eta\bigl(g(\mathsf M_2(z),\mathsf M_3(z),\ldots)\bigr).
\]
The \(D\)-exponential truncates by pointwise nilpotence. Its coefficients
may depend on \(n\), which is fixed while taking the differences.
Put \(\varrho=\sqrt n\) and, for \(\eta\subseteq\mu\), set
\[
 r_\eta=m-|\eta|,\quad c_\eta=\frac{f^{\mu/\eta}}{r_\eta!},\quad
 P_\mu(z)=\sum_{\eta\subseteq\mu}c_\eta z^{r_\eta},\quad
 B_{d,n}(z)=\frac{\mathsf M_d(z)-a_dn^{d/2}I}{n^{(d-1)/2}}.
\]
Let \(B_n=(B_{d,n})_{d\in\mathcal D}\) and
\(\widetilde V_n=((q_d-a_dn^{d/2})/n^{(d-1)/2})_{d\in\mathcal D}\).
Lemma~\ref{lem:skew-transfer}, applied to the centered polynomial
observable, gives the exact identity
\begin{equation}\label{eq:skew-asy-exact-centered}
 \mathbb E_{M_{\mu,n}}P(\widetilde V_n)
 =\frac{\mathcal T_n\!\left[
     \sum_{\eta\subseteq\mu}c_\eta z^{r_\eta}
                 \phi_\eta(P(B_n(z)))\right]}
        {\mathcal T_n(P_\mu)}.
\end{equation}
In particular, \(\mathcal T_n(P_\mu)=Z_\mu(n)/I_n\).
Formula~\eqref{clt:operator-local-expansion}, with \(N\) replaced by
\(n\), shows that each
\(p_{\eta,n}(h)=\phi_\eta(P(B_n(\varrho+h)))\) has bounded degree
and bounded coefficients; its coefficients differ by
\(O(\varrho^{-1})\) from those of
\(\phi_\eta(P((d a_dhI+L_d)_{d\in\mathcal D}))\).
Only finitely many fixed operator words occur, and \(\mu\) has
finitely many subpartitions, so these bounds hold for every required row.
By Lemma~\ref{clt:analytic-coefficient-estimate},
\begin{equation}\label{eq:skew-asy-denominator}
 \mathcal T_n(P_\mu)=c_\varnothing\varrho^m
                         (1+O_\mu(\varrho^{-1})),\qquad
 c_\varnothing=f^\mu/m!>0.
\end{equation}
After division by this denominator, every nonempty \(\eta\) contributes
\(O_{\mu,P}(\varrho^{-|\eta|})\). The empty term equals
\[
 (\tau\otimes\phi)P((d a_dhI+L_d)_{d\in\mathcal D})
       +O_{\mu,P}(\varrho^{-1}),
\]
because \((\varrho+h)^m/\varrho^m=1+O(\varrho^{-1})\)
coefficientwise in bounded degree. The size-mode cancellation
\eqref{clt:mode-one-cancellation} identifies this functional with
\(\mathbb EP(V)\). This proves the quantitative assertion at the
increment scale. For each \(d\),
\[
 V_{d,N}^\mu
 =\left(\frac nN\right)^{(d-1)/2}\widetilde V_{d,n}
    +a_d\frac{n^{d/2}-N^{d/2}}{N^{(d-1)/2}}.
\]
The coefficient is \(1+O_\mu(n^{-1})\) and the translation is
\(O_\mu(n^{-1/2})\). The same bound applies to subtraction of the
fixed \(q_d(\mu)\). Bounded moments justify these affine changes;
the first-moment estimate also permits separate exact centering.
The finite number of monomials of bounded degree proves the uniform
coefficient bound. Gaussian moment determinacy, applied to every real
linear combination, gives joint convergence in distribution.
\end{proof}

\subsection{Contents and polynomial profile tests}

For the fixed inner shape, abbreviate
\(V_{d,N}^\mu=(q_d-a_dN^{d/2})/N^{(d-1)/2}\) to \(V_{d,N}\).
The law of these statistics remains \(M_{\mu,n}\).

\begin{lemma}[Finite geometric identities]
\label{lem:clt-geometric-moments}
For \(k\ge0\), put
\(\Theta_{k,N}=\int x^k\Delta_{\lambda,N}^{Y,\mathrm{out}}(x)\,\dd x\).
Then
\begin{equation}\label{eq:clt-geometric-moments}
 \Theta_{k,N}=\frac{V_{k+2,N}}{(k+1)(k+2)},\qquad \Theta_{0,N}=0.
\end{equation}
For \(j\ge0\) there are constants \(\beta_{jb}\) such that
\begin{equation}\label{eq:skew-asy-content-leading}
 C_j=\frac{q_{j+2}}{(j+1)(j+2)}+\sum_{b=0}^{j}\beta_{jb}q_b.
\end{equation}
For \(j=0\), this reads \(C_0=q_2/2\).
Finally, if \(t_N=N/n\), then
\begin{align*}
 \int x^k\Delta_{\lambda,n}^{Y,\mathrm{out}}(x)\,\dd x
 &=t_N^{(k+1)/2}\Theta_{k,N}
       +\sqrt n\bigl(t_N^{(k+2)/2}-1\bigr)\gamma_k,\\
 \int p(x)(\Delta_{\lambda/\mu,n}^{Y}
             -\Delta_{\lambda,n}^{Y,\mathrm{out}})(x)\,\dd x
 &=-n^{-1/2}\int p(u/\sqrt n)\sigma_\mu(u)\,\dd u
                    =O_{\mu,p}(n^{-1/2}).
\end{align*}
\end{lemma}
\begin{proof}
The slope jumps give
\(\sigma_\lambda''=\sum_{x\in X}\delta_x-\sum_{y\in Y}\delta_y-\delta_0\).
Twice integrating against \(x^{k+2}\) and rescaling proves the first
identity, since
\((k+1)(k+2)\int x^k\sigma_\Omega(x)\,\dd x=a_{k+2}\).
The latter integral follows directly by the substitution
\(x=2\cos\theta\) in the formula for \(\Omega\).
For \eqref{eq:skew-asy-content-leading} with \(j\ge1\), use
Lemma~\ref{lem:content-antiderivative}: the centered second difference
of \(z^{j+2}\) has no \(z^{j-1}\) term, so matching coefficients
forces the \(z^{j+1}\) coefficient of \(\Phi_j\) to vanish. Its
leading coefficient is \(1/((j+1)(j+2))\). For \(j=0\), use
\(C_0=q_2/2\) directly; the free linear term of \(\Phi_0\) contributes
only a multiple of \(q_1=0\). The last two formulas are changes of variables
and subtraction of the fixed inner profile.
\end{proof}

\begin{corollary}[Content and profile fluctuations]
\label{cor:skew-content-clt}\label{thm:diagram-clt}
The content and polynomial-profile assertions of
Theorem~\ref{thm:skew-asymptotics-intro} hold jointly in distribution
and in all mixed moments, with the same Gaussian family
\((\Xi_r)_{r\ge2}\). Exact expectation centering gives the corresponding
centered Gaussian limits.
\end{corollary}
\begin{proof}
By \eqref{eq:skew-asy-content-leading}, for each \(j\ge0\),
\[
 \frac{C_j(\lambda/\mu)-\gamma_jN^{j/2+1}}{N^{(j+1)/2}}
 =\frac{V_{j+2,N}^\mu}{(j+1)(j+2)}+O_{L^p}(N^{-1/2})
\]
for every fixed finite \(p\). Indeed the remaining corner terms have
indices at most \(j\), and all their normalized moments are bounded
by Corollary~\ref{thm:skew-corner-clt}; the fixed inner subtraction
has the same or smaller order.

For the profile modes use the finite identity
\begin{equation}\label{eq:finite-chebyshev-transform}
 x^k=\sum_{a=0}^{\lfloor k/2\rfloor}
       \left(\binom{k}{a}-\binom{k}{a-1}\right)u_{k-2a}(x),\qquad
 \binom{k}{-1}=0.
\end{equation}
This follows by induction from \(xu_j=u_{j+1}+u_{j-1}\), with
\(u_{-1}=0\). By \eqref{eq:clt-geometric-moments}, the limit of
\(\Theta_{k,N}\) is \(V_{k+2}/((k+1)(k+2))\). Its coefficient of
\(\Xi_{k+1-2a}\), for \(k+1-2a\ge2\), is
\[
 \frac{1}{k+1}\binom{k+1}{a}
 =\frac{\binom{k}{a}-\binom{k}{a-1}}{k-2a+1}.
\]
The triangular identity \eqref{eq:finite-chebyshev-transform}, starting
with \(\Theta_{0,N}=0\), therefore shows that
\[
 \int u_k(x)\Delta_{\lambda,N}^{Y,\mathrm{out}}(x)\,\dd x
       \longrightarrow\frac{\Xi_{k+1}}{k+1}\qquad(k\ge1).
\]
Every fixed family of these tests is a finite linear transformation
of one corner vector, so the convergence holds jointly with all mixed
moments. The independent Gaussian coordinates have the means and
variances stated in Theorem~\ref{thm:skew-asymptotics-intro}.
The exact scale comparison in Lemma~\ref{lem:clt-geometric-moments}
changes each fixed polynomial pairing by \(O(n^{-1/2})\) in every
fixed \(L^p\), since \(N/n=1+O_\mu(n^{-1})\). The fixed-inner
subtraction has the same deterministic order. Both skew and limiting
charge densities have mass one, hence the zeroth skew mode is exactly
zero. Mixed-moment convergence permits exact expectation centering.
\end{proof}

\begin{proof}[Proof of Theorem~\ref{thm:skew-asymptotics-intro}]
Combine Corollaries~\ref{thm:skew-corner-clt}
and~\ref{cor:skew-content-clt}. Their proofs use the same finite
corner vector and the same current limits, so the three classes of
statistics converge jointly with the common family \((\Xi_r)_{r\ge2}\).
\end{proof}

Write \(C_c^\infty(\mathbb R)\) for the smooth compactly supported
test functions, where
\(\supp\varphi=\overline{\{x:\varphi(x)\ne0\}}\).
A distribution is a linear functional \(T\) on this
space such that, for each compact set \(K\), there are
\(C_K\ge0\) and an integer \(r_K\ge0\) with
\[
 |\langle T,\varphi\rangle|
 \le C_K\max_{0\le j\le r_K}\sup_{x\in K}|\varphi^{(j)}(x)|
 \qquad(\supp\varphi\subseteq K).
\]
It has order at most \(r\) if one can always take \(r_K\le r\).
It is supported in a closed set \(F\) if its pairing vanishes for every
test supported in \(\mathbb R\setminus F\)
\cite[Chapter~II, Sections~2.1--2.3]{HormanderDistributions}.
An almost surely defined family of distributions is a Gaussian random
distribution if its test pairings are measurable and every finite family
of real test pairings is jointly Gaussian.
Pairing a compactly supported distribution with a polynomial means
inserting a smooth cutoff equal to one near its support; the support
condition makes the value independent of that choice.

\begin{corollary}[Gaussian field representing the operator limits]
\label{thm:polynomial-gaussian-field}
The limiting polynomial pairings are represented by the Gaussian
random distribution supported on \([-2,2]\) with angular series
\begin{equation}\label{eq:clt-field-series}
 \Delta^Y(2\cos\theta)=\frac1\pi\sum_{r\ge2}
                 \frac{\Xi_r}{r}\sin(r\theta).
\end{equation}
The series is interpreted by testing against
\(2\sin\theta\,\varphi(2\cos\theta)\), where
\(\varphi\in C_c^\infty(\mathbb R)\). Its mean equals
\(m_Y(x)=\frac14-\sqrt{4-x^2}/(2\pi)\) for \(|x|<2\), and is zero
outside \([-2,2]\). In particular
\(\langle\Delta^Y,u_0\rangle=0\) and
\(\langle\Delta^Y,u_k\rangle=\Xi_{k+1}/(k+1)\) for \(k\ge1\).
\end{corollary}
\begin{proof}
Two integrations by parts give the uniform bound
\[
 \left|\int_0^\pi\varphi(2\cos\theta)\sin\theta\sin(r\theta)
               \,\dd\theta\right|
 \le Cr^{-2}\|\varphi\|_{C^2([-2,2])},
\]
where the norm is the maximum of the suprema of \(\varphi\) and its
first two derivatives. Writing
\(\Xi_r=e_r+\sqrt{2r}\,\zeta_r\) with independent standard Gaussian
\(\zeta_r\), one has
\(\mathbb E\sum_{r\ge2}|\zeta_r|r^{-5/2}<\infty\).
Thus on one event of probability one the centered tested series is
bounded by a finite random constant times this norm, simultaneously
for all tests. It defines a distribution of order at most two. Orthogonality gives the displayed polynomial modes.
For the mean, use
\(\sum_{j\ge0}\sin((2j+1)\theta)/(2j+1)=\pi/4\) on \(0<\theta<\pi\),
and remove the \(r=1\) term. This gives the displayed \(m_Y\).
The polynomial modes identify the field with the limits just proved.
The convergence of finite diagrams concerns fixed polynomial tests,
with all mixed moments.
\end{proof}

\paragraph{The ordinary specialization.}
For \(\mu=\varnothing\), one has \(n=N\),
\(M_{\mu,n}=M_N\), and \(\omega_{\lambda/\mu}=\lambda\).
Thus Theorems~\ref{thm:skew-shape-intro}
and~\ref{thm:skew-asymptotics-intro} recover the classical Gelfand
limit shape and Gaussian fluctuations of M\'eliot
\cite{Meliot,MeliotKerov}. In our half-profile convention,
\(\Delta^Y\) is half the field in \cite[Theorem~13]{Meliot}.
The general fixed-inner-shape statements have been obtained directly
from \(D\)-identities and their coefficient asymptotics.

\section{Quantitative inner-shape corrections}\label{sec:clt}

The difference-operator argument of Section~\ref{sec:limit-shape}
already gives the fixed-inner-shape limit laws with quantitative corner
moments. We first record the corresponding error estimate for contents
and profile tests, then use inserted operator moments to compute the
corrections at a common outer size.

\subsection{An error bound for mixed moments}

\begin{corollary}[An error bound for fixed mixed moments]
\label{cor:fixed-mixed-moment-rate}
Fix a partition \(\mu\), put \(N=n+|\mu|\), and sample \(\lambda\)
from \(M_{\mu,n}\). Let \(X_n\) be any fixed finite vector drawn from
\[
 \begin{gathered}
 \frac{q_d(\lambda)-a_dN^{d/2}}{N^{(d-1)/2}}\quad(d\ge2),\qquad
 \frac{C_j(\lambda/\mu)-\gamma_jN^{j/2+1}}{N^{(j+1)/2}}\quad(j\ge0),\\
 \int_\R p(x)\Delta_{\lambda/\mu,n}^{Y}(x)\,\dd x
 \quad(p\in\R[x]).
 \end{gathered}
\]
Let \(X\) be the corresponding joint Gaussian vector, using the same
family \((\Xi_r)_{r\ge2}\), whose respective coordinates are
\[
 V_d,\qquad \frac{V_{j+2}}{(j+1)(j+2)},\qquad
 \langle\Delta^Y,p\rangle.
\]
For every fixed polynomial \(P\) in the coordinates of this vector,
\begin{equation}\label{eq:fixed-mixed-moment-rate}
 \E_{M_{\mu,n}}P(X_n)=\E P(X)+O(n^{-1/2}).
\end{equation}
The estimate also holds with both vectors centered by their respective
exact means. One may replace \(N\) by \(n\) in the content and corner
normalizations, replace \(q_d(\lambda)\) by
\(q_d(\lambda)-q_d(\mu)\), or use the outer-profile polynomial tests
in place of the genuine skew tests. The implied constants may depend
on \(\mu\), the chosen finite family and \(P\).
\end{corollary}
\begin{proof}
Let \(V_N\) be the vector of the finitely many normalized corner
statistics required by the chosen coordinates, and let \(V\) be its
limit from Corollary~\ref{thm:skew-corner-clt}. The content identity
\eqref{eq:skew-asy-content-leading}, the geometric identity
\eqref{eq:clt-geometric-moments} and the scale comparisons in
Lemma~\ref{lem:clt-geometric-moments} give
\[
 X_n=A_nV_N+b_n,\qquad X=AV,\qquad
 A_n=A+O(n^{-1/2}),\qquad b_n=O(n^{-1/2}),
\]
where \(A_n,A\) are deterministic matrices, \(b_n\) is a deterministic
vector, and the bounds hold entrywise. The same bounds apply to the
normalization variants in the statement. The polynomials
\(P(A_ny+b_n)\) and \(P(Ay)\) differ coefficientwise by \(O(n^{-1/2})\).
Corollary~\ref{thm:skew-corner-clt} bounds every fixed moment of \(V_N\)
and gives the rate for the fixed test \(P(Ay)\). This proves
\eqref{eq:fixed-mixed-moment-rate}. Applying the estimate to the
coordinate functions gives \(\E X_n-\E X=O(n^{-1/2})\), so subtracting
the exact means preserves the same coefficient estimates.
\end{proof}

\subsection{The first correction caused by a fixed inner partition}
\label{subsec:adv-inner-correction}

Fix a partition \(\mu\) of size \(m\) and let \(N\) be the common outer
size. Write \(\mathbb E_{\mu,N}\) and \(\operatorname{Var}_{\mu,N}\)
for expectation and variance under \(M_{\mu,N-m}\). Use subscripts
\(0,N\) for the ordinary law \(M_N\). Both laws are on partitions of
size \(N\). We use the same observable under either law and write
\[
 Y_{d,N}=\frac{q_d-a_dN^{d/2}}{N^{(d-1)/2}},\qquad
 Y_d=\sum_{r\ge2}b_{dr}\Xi_r,\qquad d\ge2.
\]
Thus \(Y_{d,N}\) is \(V_{d,N}^{\mu}\) from
Corollary~\ref{thm:skew-corner-clt} under the skew law, and \(Y_d=V_d\).
The centering \(a_dN^{d/2}\) is common and deterministic. We specify
separate exact expectation centering when it is used.

For a finite index set \(\mathcal D\subseteq\{2,3,\ldots\}\), put
\(Y_N=(Y_{d,N})_{d\in\mathcal D}\), \(Y=(Y_d)_{d\in\mathcal D}\),
and use \(x=(x_d)_{d\in\mathcal D}\) for the polynomial test variables.
For a polynomial \(P\) let \(\|P\|_1\) be the sum of the absolute
values of its monomial coefficients, and define the directional
derivative
\[
 \mathfrak D_r=\sum_{d\in\mathcal D}b_{dr}\partial_{x_d}
 \qquad(r\ge2).
\]
Expectations with no subscripts use the Gaussian variables \(\Xi_r\)
in limiting formulas. Tests may have complex coefficients. The
covariance \(\operatorname{Cov}(U,V)=\mathbb E(UV)-\mathbb EU\,\mathbb EV\)
is bilinear for such tests.

\begin{theorem}[The first inner-shape correction]
\label{thm:adv-first-inner-correction}
For every fixed \(\mu\), finite \(\mathcal D\) and integer \(K\ge0\),
uniformly over polynomials \(P\) of degree at most \(K\),
\begin{align}
 \mathbb E_{\mu,N}P(Y_N)-\mathbb E_{0,N}P(Y_N)
 &=\frac{C_1(\mu)}N\,
            \mathbb E\!\left[\Xi_2P(Y)\right]
       +O_{\mu,\mathcal D,K}\!\left(\|P\|_1N^{-3/2}\right) \label{eq:adv-first-inner}\\
 &=\frac{4C_1(\mu)}N\,
            \mathbb E[(\mathfrak D_2P)(Y)]
       +O_{\mu,\mathcal D,K}\!\left(\|P\|_1N^{-3/2}\right).\notag
\end{align}
The $N^{-1}$ correction vanishes for every fixed corner-polynomial
test if and only if $C_1(\mu)=0$. In particular,
\[
 \mathbb E_{\mu,N}Y_{3,N}-\mathbb E_{0,N}Y_{3,N}
       =\frac{12C_1(\mu)}N+O(N^{-3/2}).
\]
If $C_1(\mu)=0$, set
$\kappa_3(\mu)=C_2(\mu)-\binom m2$. Then
\begin{align}
 \mathbb E_{\mu,N}P(Y_N)-\mathbb E_{0,N}P(Y_N)
 &=\frac{\kappa_3(\mu)}{N^{3/2}}
       \mathbb E\!\left[(\Xi_3-1)P(Y)\right]
            +O_{\mu,\mathcal D,K}\!\left(\|P\|_1N^{-2}\right) \label{eq:adv-second-inner}\\
 &=\frac{6\kappa_3(\mu)}{N^{3/2}}
       \mathbb E[(\mathfrak D_3P)(Y)]
            +O_{\mu,\mathcal D,K}\!\left(\|P\|_1N^{-2}\right).\notag
\end{align}
Subject to $C_1(\mu)=0$, this correction vanishes for every fixed
corner-polynomial test if and only if $\kappa_3(\mu)=0$.
\end{theorem}

The necessity assertions concern all fixed corner families. For a
particular family the corresponding directional derivative may vanish;
for example, $b_{d2}=0$ for every even $d$. The condition $C_1(\mu)=0$
can hold for non-self-conjugate shapes, as illustrated after
Corollary~\ref{cor:response-three-levels}.

We use Gaussian integration by parts in the form
\[
 \mathbb E[(\Xi_r-e_r)P(Y)]
       =2r\,\mathbb E[(\mathfrak D_rP)(Y)]\qquad(r\ge2)
\]
(see \cite{SteinNormal}). For a one-dimensional Gaussian density of
mean \(a\) and variance \(v>0\), its derivative equals
\(-(t-a)/v\) times that density. Integration by parts therefore gives
\(\mathbb E[(Z-a)f(Z)]=v\mathbb Ef'(Z)\) for every polynomial \(f\).
Conditioning on the other independent modes gives the displayed
identity, also for complex coefficients; polynomial growth makes the
boundary terms zero.

\paragraph{A finite current insertion identity.}
Expand the Schur polynomial of the fixed partition as
\begin{equation}\label{eq:adv-schur-expansion}
 s_\mu(p)=c_\mu p_1^m+
  \sum_{\substack{\rho:\ m_1(\rho)=0\\2\le |\rho|\le m}}
       a_{\mu,\rho}\,p_\rho p_1^{m-|\rho|},
 \qquad c_\mu=\frac{f^\mu}{m!}.
\end{equation}
Here \(\rho\) ranges over partitions with no parts equal to \(1\).
Call \(w=|\rho|\) the insertion weight. We take \(a_{\mu,\rho}=0\)
when \(|\rho|>m\), and put \(u_{\mu,\rho}=a_{\mu,\rho}/c_\mu\).
Unless marked otherwise, subsequent sums over \(\rho\) use the range
in \eqref{eq:adv-schur-expansion}.

For \(\rho\) with no parts equal to \(1\), including the empty partition,
and \(N\ge\max\{|\rho|,1\}\), put
\[
 \widehat\chi_{\rho,N}(\lambda)
   =\frac{\chi^\lambda(\rho\cup 1^{N-|\rho|})}{f^\lambda}
   =\frac{\Sigma_\rho(\lambda)}{\fall{N}{|\rho|}}
   \qquad(|\lambda|=N).
\]
The notation is that of \eqref{eq:central-character}. For the empty
partition the normalized character equals \(1\). Define
\(\alpha_\rho=\prod_{r\ge2}\alpha_r^{m_r(\rho)}\), with
\(\alpha_{\varnothing}=I\), and introduce the coefficient functional
\begin{equation}\label{eq:adv-insertion-functional}
 \mathcal I_{\rho,N}(P)
 =\frac1{\fall{N}{w}}
   \mathcal T_N\!\left[
        z^w\phi\!\left(\alpha_\rho P(\mathsf Z_N(z))\right)\right],
 \qquad
 \mathsf Z_{d,N}(z)=\frac{\mathsf M_d(z)-a_dN^{d/2}I}{N^{(d-1)/2}}.
\end{equation}
Here \(\mathsf Z_N(z)=(\mathsf Z_{d,N}(z))_{d\in\mathcal D}\). The components commute
because they are simultaneous conjugates of multiplication operators,
as in Lemma~\ref{clt:operator-moments-lemma}. The character interpretation is
\begin{equation}\label{eq:adv-insertion-character}
 \mathcal I_{\rho,N}(P)
   =\mathbb E_{0,N}\!\left[\widehat\chi_{\rho,N}P(Y_N)\right].
\end{equation}
In particular, Stanley's exact skew-character formula
\cite[Theorem~3.1, equation~(12)]{StanleySkew} yields
\begin{equation}\label{eq:adv-exact-insertion-ratio}
 \mathbb E_{\mu,N}P(Y_N)=
 \frac{\mathcal I_{\varnothing,N}(P)+
       \sum_\rho u_{\mu,\rho}\mathcal I_{\rho,N}(P)}
      {1+\sum_\rho u_{\mu,\rho}\mathcal I_{\rho,N}(1)},
 \qquad
 \mathcal I_{\varnothing,N}(P)=\mathbb E_{0,N}P(Y_N).
\end{equation}
This identity holds for \(N\ge\max\{m,1\}\).

\begin{proof}[Proof of the character interpretation and the ratio]
Since \(D=\alpha_1-\alpha_{-1}-I\), the Heisenberg relation gives
\[
 e^{zD}=e^{-z-z^2/2}e^{-z\alpha_{-1}}e^{z\alpha_1},
\]
by the same factorization as in Lemma~\ref{clt:currents-lemma}.
The operator \(\alpha_\rho\) commutes with \(\alpha_{\pm1}\), and
\(\langle\varnothing|\alpha_{-1}=0\). Hence
\(I(z)\langle\varnothing|\alpha_\rho e^{zD}
=\langle\varnothing|\alpha_\rho e^{z\alpha_1}\).
Conjugation of the multiplication operator \(P(Y_N)\), together with
\(D\mathbf1=0\), gives
\[
 I(z)\phi\!\left(\alpha_\rho P(\mathsf Z_N(z))\right)
 =\sum_{l\ge0}\frac{z^l}{l!}
      \langle\varnothing|\alpha_\rho\alpha_1^l
                         P(Y_N)|\mathbf1\rangle .
\]
By Lemma~\ref{lem:external-schur-identities} and the Frobenius expansion
\(s_\lambda=\sum_{\nu\vdash N}\chi^\lambda(\nu)p_\nu/z_\nu\)
\cite[Chapter~I, (7.2), (7.5)]{MacdonaldSF}, for \(|\lambda|=N\) one has
\[
 \langle\varnothing|\alpha_\rho\alpha_1^{N-w}|\lambda\rangle
     =\chi^\lambda(\rho\cup 1^{N-w}).
\]
Indeed, the specialization \(p_r=rt_r\) turns \(p_\nu/z_\nu\) into
\(\prod_r t_r^{m_r(\nu)}/m_r(\nu)!\), so the matching exponential
coefficients give this matrix element. The operator lowers rank by
\(w+(N-w)=N\), so its vacuum row is supported on rank \(N\).
Multiplying the coefficient of \(z^{N-w}\) by \((N-w)!/I_N\) now
proves \eqref{eq:adv-insertion-character}. Every coefficient uses
finitely many row-finite matrix elements.

For the ratio, let \(R_{\mu,N}(\lambda)=f^{\lambda/\mu}/(c_\mu f^\lambda)\),
where \(f^{\lambda/\mu}=0\) if \(\mu\not\subseteq\lambda\).
In \eqref{eq:adv-schur-expansion}, the Frobenius formula gives
\(a_{\mu,\rho}=\chi^\mu(\rho\cup1^{m-w})/(z_\rho(m-w)!)\) for
\(w=|\rho|\le m\), because \(m_1(\rho)=0\). Stanley's formula gives
\begin{equation}\label{eq:adv-density}
 R_{\mu,N}=1+\sum_\rho u_{\mu,\rho}\widehat\chi_{\rho,N},
 \qquad
 \frac{d M_{\mu,N-m}}{dM_N}
   =\frac{R_{\mu,N}}{\mathbb E_{0,N}R_{\mu,N}}.
\end{equation}
The last quotient is the ratio of the masses of the two laws at each
partition of size \(N\). Its denominator is \(Z_\mu(N-m)/(c_\mu I_N)>0\).
Taking the expectation of \(P(Y_N)\) with this density gives
\eqref{eq:adv-exact-insertion-ratio}.
\end{proof}

\paragraph{Inserted Gaussian moments.}
For \(r\ge2\), define the Hermite polynomials \(H_j^{[r]}\) with
variance parameter \(r\) by
\[
 \exp(uw-ru^2/2)=\sum_{j\ge0}H_j^{[r]}(w)\frac{u^j}{j!}.
\]
For variables \(\boldsymbol w=(w_r)_{r\ge2}\) indexed by the Gaussian
modes, put
\(\mathfrak H_\rho(\boldsymbol w)=\prod_{r\ge2}H_{m_r(\rho)}^{[r]}(w_r)\).
Only finitely many factors differ from \(H_0^{[r]}=1\). In particular,
\(H_1^{[r]}(w)=w\) and \(H_2^{[r]}(w)=w^2-r\). The parameter \(r\)
in the generating series comes from the current commutator below.
The limiting mode \(\Xi_r\) has variance \(2r\).
For every fixed $\rho,\mathcal D,K$, uniformly for $\deg P\le K$,
\begin{equation}\label{eq:adv-inserted-limit}
 \mathcal I_{\rho,N}(P)
 =N^{-w/2}\left\{
       \mathbb E\!\left[\mathfrak H_\rho(\Xi)P(Y)\right]
          +O\!\left(\|P\|_1N^{-1/2}\right)\right\}.
\end{equation}

\begin{proof}
Put $\varrho=\sqrt N$ and $z=\varrho+h$. The leading conjugation
formula gives, coefficientwise in a fixed finite operator span,
\[
 \mathsf Z_{d,N}(\varrho+h)=d a_dhI+L_d+R_{d,N}(h),
\]
where \(R_{d,N}\) has bounded degree in \(h\) and coefficients of order
\(O(\varrho^{-1})\) in the sense of \eqref{clt:operator-local-expansion}.
Lemma~\ref{clt:analytic-coefficient-estimate} gives
\(\tau_N(h^a)=\tau(h^a)+O_a(\varrho^{-1})\), and
\(z^w/\fall{N}{w}=N^{-w/2}(1+h/\varrho)^w(1+O(N^{-1}))\).
These observations reduce the leading coefficient to
\[
 (\tau\otimes\phi)\!\left[
      \alpha_\rho P\bigl((d a_dhI+L_d)_{d\in\mathcal D}\bigr)\right].
\]
To identify all these finite moments simultaneously, introduce finitely
supported formal parameters $u_r$ ($r\ge2$) and $t_r$ ($r\ge1$).
To move the terms with negative subscripts to the left, first apply
the factorization of Lemma~\ref{clt:currents-lemma} to \(e^{\sum_r t_r A_r}\).
Commuting \(e^{\sum_{r\ge2}u_r\alpha_r}\) past
\(e^{\sum_{r\ge1}t_r\alpha_{-r}}\) produces the scalar
\(\exp(\sum_{r\ge2}r u_rt_r)\). The negative-subscript exponentials
then act as the identity on the vacuum bra. Evaluating the remaining
positive-subscript exponentials by \eqref{eq:positive-current-gaussian} gives
\begin{align*}
 \phi\!\left(e^{\sum_{r\ge2}u_r\alpha_r}
                   e^{\sum_{r\ge1}t_rA_r}\right)
 =\exp\!\left(
   \sum_{r\ge1}(e_rt_r+rt_r^2)
   +\sum_{r\ge2}(e_ru_r+2rt_ru_r+ru_r^2/2)\right).
\end{align*}
For this display only, also take an independent $\Xi_1\sim \mathcal N(1,2)$.
The expression equals
\[
 \mathbb E\left[e^{\sum t_r\Xi_r}
             \prod_{r\ge2}e^{u_r\Xi_r-ru_r^2/2}\right]
\]
before the size-mode cancellation. Substituting
\(t_r=\sum_{d\in\mathcal D}b_{dr}v_d\), with formal parameters \(v_d\),
and applying \(\tau\) cancels the \(r=1\) terms by
\eqref{clt:mode-one-cancellation}. Coefficient comparison gives the
asserted inserted Gaussian moments. The error bound follows from the
finite number of operator words and scalar moments. For fixed
$\mathcal D,K$ there are finitely many monomials, so the bound is
uniform with the factor $\|P\|_1$. This proves the estimate directly from the exact transfer, leading
conjugation and current identities.
\end{proof}

\begin{lemma}[The first Schur coefficients]
\label{lem:adv-first-schur-coefficients}
The coefficients in \eqref{eq:adv-schur-expansion} satisfy
\[
 u_{\mu,(2)}=C_1(\mu),\qquad
 u_{\mu,(3)}=C_2(\mu)-\binom m2.
\]
\end{lemma}
\begin{proof}
The first formula follows immediately by comparing the coefficient of
\(t^{m-1}\) for \(m\ge1\) in the hook-content specialization
\cite[Chapter~I, Section~3, Example~4]{MacdonaldSF}
\[
 s_\mu(1^t)=\frac1{H_\mu}\prod_{\square\in\mu}(t+c(\square)).
\]
Here \(1^t\) is the alphabet of \(t\) ones for integers \(t\ge0\),
and the identity extends polynomially in \(t\). For \(m=0\) the
coefficient and \(C_1(\mu)\) are both zero by convention.

For use in the next lemma, we compute the one-part coefficients from
Jacobi--Trudi \cite[Chapter~I, (3.4)]{MacdonaldSF}. Take \(L\ge\max\{1,\ell(\mu)\}\), put
\(a_i=\mu_i+L-i\) for \(1\le i\le L\), and define
\(\Pi(x)=\prod_{i=1}^L(x-a_i)\). At \(p_1=1\) and \(p_r=0\) for
\(r\ge2\), the Jacobi--Trudi matrix has entries \(1/(a_i-L+j)!\).
Its derivative with respect to \(p_r\) has entries
\(1/(r(a_i-L+j-r)!)\), with \(1/b!=0\) when \(b\) is a negative
integer. Multiplying row \(i\) by \(a_i!\) and reversing columns gives
the evaluation matrix \((\fall{a_i}{j})_{1\le i\le L,\,0\le j<L}\).
In the derivative, row \(i\) is evaluated at \(a_i-r\) and multiplied
by \(\fall{a_i}{r}/r\). The ratio of its determinant to the original
determinant follows from the Vandermonde formula. Summing over
differentiated rows gives
\begin{align*}
 u_{\mu,(r)}
 &=\frac1r\sum_i \fall{a_i}{r}
          \prod_{j\ne i}\frac{a_i-r-a_j}{a_i-a_j}\\
 &=-\frac1{r^2}[x^{-1}]
          \frac{\fall{x}{r}\Pi(x-r)}{\Pi(x)}.
\end{align*}
The second equality sums the simple-pole residues at the \(a_i\), which
gives the coefficient of \(x^{-1}\) in the Laurent expansion at infinity.
For \(r=2,3\), expand this finite rational expression using
\begin{align*}
 \sum_i(a_i-(L-i))&=m,\\
 \sum_i(a_i^2-(L-i)^2)&=2C_1(\mu)+(2L-1)m,\\
 \sum_i(a_i^3-(L-i)^3)&=3C_2(\mu)+(6L-3)C_1(\mu)
                                  +(3L^2-3L+1)m.
\end{align*}
The resulting values are $C_1(\mu)$ and
$C_2(\mu)-m(m-1)/2$, respectively.
\end{proof}

\begin{lemma}[The weight-four Schur coefficients]
\label{lem:response-weight-four}
With the notation of \eqref{eq:adv-schur-expansion}, for every partition
\(\mu\) of size \(m\),
\begin{align}
 u_{\mu,(2,2)}
  &=\frac12C_1(\mu)^2-\frac32C_2(\mu)+\binom m2,
                       \label{eq:response-u22}\\
 u_{\mu,(4)}
  &=C_3(\mu)-(2m-3)C_1(\mu).
                       \label{eq:response-u4}
\end{align}
Coefficients of weight larger than \(m\) are understood to be zero.
\end{lemma}
\begin{proof}
For \(m\ge2\), compare the coefficient of \(t^{m-2}\) in
\[
 \frac{s_\mu(1^t)}{c_\mu}
  =\prod_{\square\in\mu}(t+c(\square)).
\]
On the left, \(p_r(1^t)=t\), so a term with type \(\rho\)
contributes a power \(t^{m-|\rho|+\ell(\rho)}\). The exponent is
\(m-2\) for the types \((3)\) and \((2,2)\). Thus
\[
 u_{\mu,(3)}+u_{\mu,(2,2)}
       =\frac12(C_1(\mu)^2-C_2(\mu)).
\]
Lemma~\ref{lem:adv-first-schur-coefficients} proves
\eqref{eq:response-u22}. The cases \(m=0,1\) are immediate.

For the second identity use the Laurent-coefficient formula in the
proof of Lemma~\ref{lem:adv-first-schur-coefficients}. Choose
\(L\ge\max\{1,\ell(\mu)\}\), put
\(a_i=\mu_i+L-i\), \(\Pi_\mu(x)=\prod_{i=1}^L(x-a_i)\), and
\(\Pi_0(x)=\fall{x}{L}\). With \(y=x-L\), row telescoping gives
\[
 \frac{\Pi_\mu(x)}{\Pi_0(x)}
     =\prod_{\square\in\mu}\frac{y-c(\square)}{y-c(\square)+1}.
\]
Consequently
\begin{equation}\label{eq:response-content-residue}
 \fall{x}{r}\frac{\Pi_\mu(x-r)}{\Pi_\mu(x)}
   =\fall{y}{r}W_r(y),\qquad
 W_r(y)=\prod_{\square\in\mu}
 \frac{(y-c(\square)-r)(y-c(\square)+1)}
      {(y-c(\square)-r+1)(y-c(\square))}.
\end{equation}
The coefficient at infinity of degree \(-1\) is unchanged by
\(x=y+L\). For \(r=4\), write
\(\log W_4(y)=\sum_{j\ge2}\ell_jy^{-j}\), using the formal logarithm
at infinity. Write \(C_i=C_i(\mu)\) in this calculation. Expansion gives
\begin{align*}
 \ell_2&=-4m,&
 \ell_3&=-8C_1-12m,\\
 \ell_4&=-12C_2-36C_1-44m,&
 \ell_5&=-16C_3-72C_2-176C_1-156m,
\end{align*}
Since
\(\fall{y}{4}=y^4-6y^3+11y^2-6y\),
\begin{align*}
 [y^{-1}]\fall{y}{4}W_4(y)
  &=\ell_5+\ell_2\ell_3
      -6(\ell_4+\ell_2^2/2)+11\ell_3-6\ell_2\\
  &=-16C_3+16(2m-3)C_1.
\end{align*}
Multiplication by \(-1/4^2\), as in the cited Laurent-coefficient
formula, proves \eqref{eq:response-u4}.
\end{proof}

\begin{proof}[Proof of Theorem~\ref{thm:adv-first-inner-correction}]
In \eqref{eq:adv-exact-insertion-ratio}, the only term of insertion
weight two is $C_1(\mu)\mathcal I_{(2),N}$. Equation
\eqref{eq:adv-inserted-limit} gives its leading term
$C_1(\mu)N^{-1}\mathbb E[\Xi_2P(Y)]$. Since
$\phi(\alpha_2)=0$, $\mathcal I_{(2),N}(1)=0$ exactly. All remaining
terms are $O(\|P\|_1N^{-3/2})$. This proves
\eqref{eq:adv-first-inner}. Gaussian integration by parts gives the
directional-derivative form. The test $P=x_3$, with $Y_3=3\Xi_2$,
gives $12C_1(\mu)$ and proves necessity of the vanishing condition.

If $C_1(\mu)=0$, the first possible insertion has weight three and
coefficient $\kappa_3(\mu)$. Its numerator leading term is
$\kappa_3N^{-3/2}\mathbb E[\Xi_3P(Y)]$, and its denominator leading
term is $\kappa_3N^{-3/2}$ because $\mathbb E\Xi_3=1$. Subtracting
the normalizer gives \eqref{eq:adv-second-inner}; the errors are
$O(\|P\|_1N^{-2})$. Gaussian integration by parts and the test
$P=x_4$, with $Y_4=4\Xi_3$, give the derivative form and necessity.
\end{proof}

\begin{theorem}[The first nonzero inner-shape correction]
\label{thm:adv-first-nonzero-inner}
In the finite expansion \eqref{eq:adv-schur-expansion}, let
\[
 s_*(\mu)=\min\{\,|\rho|:2\le|\rho|\le m,\ m_1(\rho)=0,
                         \ a_{\mu,\rho}\ne0\,\}.
\]
Set $s_*(\mu)=\infty$ if the set is empty. For \(m\ge2\) it belongs
to \(\{2,3,4\}\); Corollary~\ref{cor:response-three-levels} gives the
complete classification. If $s=s_*(\mu)<\infty$,
put
\[
 A_{\mu,s}(\boldsymbol w)=\sum_{\substack{|\rho|=s\\m_1(\rho)=0}}
                      u_{\mu,\rho}\mathfrak H_\rho(\boldsymbol w).
\]
For every fixed finite $\mathcal D$ and $K$, uniformly for
$\deg P\le K$,
\begin{equation}\label{eq:adv-general-inner-correction}
 \mathbb E_{\mu,N}P(Y_N)-\mathbb E_{0,N}P(Y_N)
 =N^{-s/2}\operatorname{Cov}\!\left(A_{\mu,s}(\Xi),P(Y)\right)
        +O_{\mu,\mathcal D,K}\!\left(\|P\|_1N^{-(s+1)/2}\right).
\end{equation}
The coefficient in \eqref{eq:adv-general-inner-correction} is nonzero
for at least one fixed corner-polynomial test. More generally, if all
Schur coefficients of weights $2,\ldots,s-1$ vanish, then
\eqref{eq:adv-general-inner-correction} holds with the displayed
$A_{\mu,s}$, and its coefficient vanishes for every fixed
corner-polynomial test if and only if all $a_{\mu,\rho}$ of weight
$s$ vanish. The only bases with $s_*(\mu)=\infty$ are
$\varnothing$ and $(1)$; for these bases the expectations at the same
outer rank agree exactly for every observable.
\end{theorem}

\begin{proof}
Write $F_N(P)=\mathcal I_{\varnothing,N}(P)$ and
$S_N(P)=\sum_\rho u_{\mu,\rho}\mathcal I_{\rho,N}(P)$. The exact
identity gives
\[
 \mathbb E_{\mu,N}P(Y_N)-F_N(P)
       =\frac{S_N(P)-F_N(P)S_N(1)}{1+S_N(1)}.
\]
The comparison uses the exact ordinary moment \(F_N(P)\) at rank
\(N\). By \eqref{eq:adv-inserted-limit}, the first possible insertion
weight \(s\) gives
\begin{align*}
 S_N(P)&=N^{-s/2}\mathbb E[A_{\mu,s}(\Xi)P(Y)]
            +O(\|P\|_1N^{-(s+1)/2}),\\
 S_N(1)&=N^{-s/2}\mathbb E A_{\mu,s}(\Xi)
            +O(N^{-(s+1)/2}),\\
 F_N(P)&=\mathbb EP(Y)+O(\|P\|_1N^{-1/2}).
\end{align*}
All sums over \(\rho\) are finite. The reciprocal of the denominator is
\(1+O(N^{-s/2})\). Its contribution to the correction is
\(O(\|P\|_1 N^{-s})\), which fits the stated remainder since \(s\ge2\). These identities prove
\eqref{eq:adv-general-inner-correction}, with the stated uniformity.
More explicitly, for each monomial of degree at most $K$, only a finite
set of fixed current/operator words and scalar $\tau_N$ moments enters
\eqref{eq:adv-inserted-limit}; the maximum of the corresponding finitely
many constants bounds every polynomial by its coefficient norm
$\|P\|_1$. No index, degree or inner partition is allowed to grow with
$N$ in this bound.

Giving \(w_r\) weight \(r\), the leading monomial of
\(\mathfrak H_\rho(\boldsymbol w)\) is
\(\prod_{r\ge2}w_r^{m_r(\rho)}\). These monomials are distinct for
distinct \(\rho\) of weight \(s\), so a nonzero \(A_{\mu,s}\) is
nonconstant. The
triangular identities
\[
 Y_{r+1}=(r+1)\Xi_r+\text{a linear combination of lower modes}
 \qquad(r\ge2)
\]
express every polynomial in $\Xi_2,\ldots,\Xi_s$ as a polynomial in
$Y_3,\ldots,Y_{s+1}$. Taking $P(Y)=A_{\mu,s}(\Xi)$ makes the displayed
covariance its variance, which is positive because the Gaussian vector
has full support in its Euclidean space. Conversely,
if all coefficients vanish, the covariance vanishes. This proves the
criterion and sharp detectability.

Finally, for $m\ge2$, $s_\mu$ cannot be a scalar multiple of $p_1^m$:
Pieri's identity \cite[Chapter~I, (5.16)]{MacdonaldSF} gives
$ p_1^m=\sum_{\nu\vdash m}f^\nu s_\nu$ with at least two nonzero
terms. The Schur expansion is therefore nontrivial. For
$\mu=\varnothing,(1)$ it consists only of its $p_1$ term, and the exact
ratio reduces to \(F_N(P)\). Equation~\eqref{eq:adv-density} then
gives equality of the two laws on the entire rank.
\end{proof}

\begin{corollary}[The correction after separate exact centering]
\label{cor:adv-centered-inner-correction}
Under the assumptions of Theorem~\ref{thm:adv-first-nonzero-inner},
with \(s=s_*(\mu)<\infty\), let \(M=\mathbb EY\). Define the part of
\(A_{\mu,s}(\Xi)\) orthogonal to constants and linear functions of the
Gaussian modes by
\[
 A_{\mu,s}^{\ge2}(\Xi)=A_{\mu,s}(\Xi)-\mathbb EA_{\mu,s}(\Xi)
   -\sum_{r=2}^s\frac{\operatorname{Cov}(A_{\mu,s}(\Xi),\Xi_r)}{2r}
                           (\Xi_r-e_r).
\]
Then, on each fixed polynomial class,
\begin{align*}
 &\mathbb E_{\mu,N}P(Y_N-\mathbb E_{\mu,N}Y_N)
       -\mathbb E_{0,N}P(Y_N-\mathbb E_{0,N}Y_N)\\
 &\hspace{12mm}=N^{-s/2}
       \mathbb E\!\left[A_{\mu,s}^{\ge2}(\Xi)P(Y-M)\right]
       +O_{\mu,\mathcal D,K}(\|P\|_1N^{-(s+1)/2}).
\end{align*}
Here all tests range over all fixed finite index families and
all fixed degrees. The coefficient vanishes for every fixed centered
corner-polynomial test if and only if $A_{\mu,s}$ is affine in the Gaussian modes.
\end{corollary}
\begin{proof}
Apply \eqref{eq:adv-general-inner-correction} to each coordinate:
\[
 \delta_{d,N}:=\mathbb E_{\mu,N}Y_{d,N}
                -\mathbb E_{0,N}Y_{d,N}
 =N^{-s/2}\operatorname{Cov}(A_{\mu,s}(\Xi),Y_d)
                     +O(N^{-(s+1)/2}).
\]
Apply the theorem to the translated polynomial
\(P(y-\mathbb E_{0,N}Y_N)\), whose coefficients converge at rate
\(O(N^{-1/2})\) to those of \(P(y-M)\). The uniform coefficient bound
permits this \(N\)-dependent test. Write
\(\delta_N=(\delta_{d,N})_{d\in\mathcal D}\). The skew-centered test
uses the additional translation by \(-\delta_N\). Taylor expansion
gives a quadratic remainder \(O(\|P\|_1N^{-s})\). The first derivative
terms, evaluated using Corollary~\ref{cor:fixed-mixed-moment-rate},
subtract
\[
 N^{-s/2}\sum_{d\in\mathcal D}
        \operatorname{Cov}(A_{\mu,s}(\Xi),Y_d)
                         \mathbb E[(\partial_{x_d}P)(Y-M)].
\]
The error from replacing the derivative moments by their limiting
values fits \(O(\|P\|_1N^{-(s+1)/2})\). For \(r\ge2\), Gaussian
integration by parts gives
\[
 \mathbb E[(\Xi_r-e_r)P(Y-M)]
   =2r\,\mathbb E[(\mathfrak D_rP)(Y-M)].
\]
Together with \(Y_d=\sum_{r\ge2}b_{dr}\Xi_r\), this identifies the
linear contribution with the projection subtracted in the definition
of \(A_{\mu,s}^{\ge2}\). The residual is orthogonal to
constants and linear Gaussian polynomials. It is zero exactly when
$A_{\mu,s}$ is affine; otherwise the expectation of its square is positive,
and the triangular mode identities furnish a polynomial test that
detects it.
\end{proof}

\begin{corollary}[The three possible first insertion weights]
\label{cor:response-three-levels}
For \(m\ge2\), the first nonzero insertion weight of
Theorem~\ref{thm:adv-first-nonzero-inner} is
\[
 s_*(\mu)=
 \begin{cases}
  2,&C_1(\mu)\ne0,\\
  3,&C_1(\mu)=0,\quad C_2(\mu)\ne\binom m2,\\
  4,&C_1(\mu)=0,\quad C_2(\mu)=\binom m2.
 \end{cases}
\]
In the last case, necessarily \(m\ge5\), and
\begin{equation}\label{eq:response-fourth-polynomial}
 A_{\mu,4}(\Xi)
   =C_3(\mu)\Xi_4-\frac{m(m-1)}4(\Xi_2^2-2).
\end{equation}
Its residual after projection on constants and linear Gaussian
polynomials is
\begin{equation}\label{eq:response-fourth-residual}
 A_{\mu,4}^{\ge2}(\Xi)
       =-\frac{m(m-1)}4(\Xi_2^2-4).
\end{equation}
Hence for all fixed polynomial tests, under the two content conditions
in the last case,
\begin{align}
 &\mathbb E_{\mu,N}P(Y_N-\mathbb E_{\mu,N}Y_N)
       -\mathbb E_{0,N}P(Y_N-\mathbb E_{0,N}Y_N)\notag\\
 &\qquad=-\frac{4m(m-1)}{N^2}
       \mathbb E[(\mathfrak D_2^2P)(Y-\mathbb EY)]
       +O_{\mu,\mathcal D,K}(\|P\|_1N^{-5/2}).
                    \label{eq:response-fourth-centered}
\end{align}
\end{corollary}
\begin{proof}
The first two cases follow from
Lemma~\ref{lem:adv-first-schur-coefficients}. If both coefficients
vanish, \eqref{eq:response-u22} gives
\(u_{\mu,(2,2)}=-m(m-1)/4\ne0\). Thus \(s_*=4\).
This already excludes \(m=2,3\), since their weight-four coefficients
are zero. At size four, the only partition with \(C_1=0\) is
\((2,2)\), for which \(C_2=2\ne6\); hence \(m\ge5\).
The two insertion types of weight four are \((4)\) and \((2,2)\),
so \eqref{eq:response-u4} gives \eqref{eq:response-fourth-polynomial}.
Since \(\Xi_2\sim \mathcal N(0,4)\) and \(\Xi_4\sim \mathcal N(0,8)\) are independent,
the linear projection removes \(C_3(\mu)\Xi_4\), and removal of the
constant replaces \(\Xi_2^2-2\) by \(\Xi_2^2-4\). This proves
\eqref{eq:response-fourth-residual}. Corollary~\ref{cor:adv-centered-inner-correction}
and Gaussian integration by parts twice give
\eqref{eq:response-fourth-centered}.
\end{proof}

Under the common deterministic centering, the three cases can all be
detected using the statistics \(Y_{3,N},Y_{4,N}\). The respective
tests are \(P=x_3\), \(P=x_4\), and \(P=x_3^2\). The third coefficient is
\(288u_{\mu,(2,2)}=-72m(m-1)\ne0\).
The third case also includes non-self-conjugate shapes. For example,
\(\mu=(9,5,3,3,3,3,2)\) has \(m=28\),
\(C_1=0\), \(C_2=378=\binom{28}{2}\), and \(C_3=630\).
Its polynomial \(A_{\mu,4}\) is \(630\Xi_4-189(\Xi_2^2-2)\).
After separate exact centering the residual is \(-189(\Xi_2^2-4)\).

\begin{lemma}[Three exact inserted second moments]
\label{lem:response-exact-second-moments}
Put
\[
 \Gamma_0(z)=\phi(\mathsf M_3(z)^2),\quad
 \Gamma_3(z)=\phi(\alpha_3\mathsf M_3(z)^2),\quad
 \Gamma_{22}(z)=\phi(\alpha_2^2\mathsf M_3(z)^2).
\]
Then
\begin{align}
 \Gamma_0(z)&=36z^2+72z^3+36z^4,\notag\\
 \Gamma_3(z)&=648+1944z+1764z^2+504z^3+36z^4,\notag\\
 \Gamma_{22}(z)&=2304+5184z+4104z^2+1296z^3+360z^4.
                       \label{eq:response-three-exact-polynomials}
\end{align}
Also
\begin{equation}\label{eq:response-linear-insertion-exact}
 \phi(\alpha_2\mathsf M_3(z))=12(1+z)^2.
\end{equation}
\end{lemma}
\begin{proof}
Lemma~\ref{lem:highest-central-commutator} gives
\(\operatorname{ad}_D^3\mathsf M_3=0\), so the first three
polynomials have degree at most four and the last has degree at most two.
For \(g=q_3^2\), the needed local values are
\[
\begin{array}{c|rrrrr}
 \eta&g(\eta)&Dg(\eta)&D^2g(\eta)&D^3g(\eta)&D^4g(\eta)\\ \hline
 \varnothing&0&0&72&432&864\\
 (3)&324&1080&2088&2160&864\\
 (2,1)&0&216&648&1296&864\\
 (4)&1296&2880&3960&3024&1728\\
 (3,1)&144&576&1224&1296&0\\
 (2,2)&0&288&1368&2160&2592
\end{array}
\]
They follow directly from
\(q_3(\lambda)=3\sum_i\lambda_i(\lambda_i+1-2i)\) and the definition
of \(D\). Diagram conjugation preserves every row value for \(g\).
The particle-move signs give
\begin{align*}
 \langle\varnothing|\alpha_3
   &=\langle(3)|-\langle(2,1)|+\langle(1,1,1)|,\\
 \langle\varnothing|\alpha_2^2
   &=\langle(4)|-\langle(3,1)|+2\langle(2,2)|
                  -\langle(2,1,1)|+\langle(1,1,1,1)|.
\end{align*}
The coefficient of \(z^j\) is obtained by applying the appropriate
bra to \(D^jg\) and dividing by \(j!\), proving
\eqref{eq:response-three-exact-polynomials}.
Finally \(\langle\varnothing|\alpha_2=\langle(2)|-\langle(1,1)|\),
where \((q_3,Dq_3,D^2q_3)((2))=(6,12,12)\); diagram conjugation negates
these values. This proves \eqref{eq:response-linear-insertion-exact}.
\end{proof}

\begin{lemma}[The transposition-square identity]
\label{lem:transposition-square}
Let \(\mathbb C[S_m]\) be the complex vector space with basis \(S_m\),
with the group multiplication extended bilinearly. Write \(I\) for
the identity permutation and let
\(\mathsf T_m=\sum_{1\le i<j\le m}(ij)\). Let
\(\mathsf C_{3,m}\) be the sum of all three-cycles and
\(\mathsf C_{22,m}\) the sum of all products of two disjoint
transpositions. A sum over an empty class is zero. Then
\[
 \mathsf T_m^2=\binom m2 I+3\mathsf C_{3,m}+2\mathsf C_{22,m}.
\]
Consequently the coefficients in \eqref{eq:adv-schur-expansion} satisfy
\[
 u_{\mu,(2)}^2=\binom m2+3u_{\mu,(3)}+2u_{\mu,(2,2)}.
\]
\end{lemma}
\begin{proof}
In the square of \(\mathsf T_m\), two equal transpositions give
the identity, accounting for \(\binom m2\). Two distinct
transpositions either share one letter or are disjoint. Each
three-cycle has exactly three ordered factorizations into two
transpositions, whereas a product of two disjoint transpositions has
exactly two. This proves the group-algebra identity.

For a partition \(\rho\) with no parts equal to one and
\(s=|\rho|\le m\), the Frobenius formula gives
\[
 u_{\mu,\rho}=\frac{m!}{z_\rho(m-s)!}
                  \frac{\chi^\mu(\rho\cup1^{m-s})}{f^\mu}.
\]
The first factor is the size of the corresponding conjugacy class.
Its class sum therefore acts on the irreducible representation
indexed by \(\mu\) as the scalar \(u_{\mu,\rho}\), by centrality
and taking traces. Applying the group-algebra identity to this
representation proves the scalar relation. The same formula holds
when a class is empty, by the stated zero convention.
\end{proof}

\paragraph{Why the constant correction depends only on the inner size.}
The calculation below gives the constant variance correction in the
form
\[
 12u_{\mu,(3)}+8u_{\mu,(2,2)}-4u_{\mu,(2)}^2
   =4\bigl(3u_{\mu,(3)}+2u_{\mu,(2,2)}-u_{\mu,(2)}^2\bigr).
\]
Lemma~\ref{lem:transposition-square} identifies this as
\(-4\binom m2=-2m(m-1)\). Thus the contributions of overlapping
and disjoint transpositions cancel after subtraction of the squared
mean. Only the identity term, which counts the pairs of equal
transpositions on the \(m\) prefix labels, remains. This explains
the independence from the shape of the fixed prefix. The next term,
\(6\kappa_3(\mu)/\sqrt N\), retains its three-cycle central character
and can distinguish inner shapes of the same size.

\begin{theorem}[The first two variance corrections]
\label{thm:response-variance-refined}\label{thm:response-universal-variance}
For every fixed partition \(\mu\) of size \(m\), at a common outer
rank \(N\), put \(\kappa_3(\mu)=C_2(\mu)-\binom m2\). Write
\(C_1=C_1(\lambda)\) for the content sum of the outer shape. Then
\begin{equation}\label{eq:response-variance-refined-content}
 \operatorname{Var}_{\mu,N}(C_1)
       -\operatorname{Var}_{0,N}(C_1)
 =-2m(m-1)+\frac{6\kappa_3(\mu)}{\sqrt N}+O_\mu(N^{-1}).
\end{equation}
Equivalently, for \(Y_{3,N}=q_3/N=6C_1/N\),
\begin{equation}\label{eq:response-variance-refined-corner}
 \operatorname{Var}_{\mu,N}(Y_{3,N})
       -\operatorname{Var}_{0,N}(Y_{3,N})
 =-\frac{72m(m-1)}{N^2}
       +\frac{216\kappa_3(\mu)}{N^{5/2}}+O_\mu(N^{-3}).
\end{equation}
Both variances use their respective exact expectations.
\end{theorem}
\begin{proof}
It suffices to prove \eqref{eq:response-variance-refined-corner}; take
\(N\ge\max\{m,8\}\) in the asymptotic calculation. For the test
\(P=x_3^2\), define
\[
 G_{\rho,N}=\mathcal I_{\rho,N}(P)
                -\mathcal I_{\varnothing,N}(P)\mathcal I_{\rho,N}(1).
\]
The exact insertion ratio gives
\begin{equation}\label{eq:response-refined-exact-subtraction}
 \mathbb E_{\mu,N}Y_{3,N}^2-\mathbb E_{0,N}Y_{3,N}^2
 =\frac{\sum_\rho u_{\mu,\rho}G_{\rho,N}}
        {1+\sum_\rho u_{\mu,\rho}\mathcal I_{\rho,N}(1)}.
\end{equation}
The vacuum-row particle moves used in the proof of
Lemma~\ref{lem:exact-content-exponential} give
\[
 \langle\varnothing|\alpha_r
      =\sum_{a=0}^{r-1}(-1)^a\langle(r-a,1^a)|.
\]
For even \(r\), diagram conjugation pairs the hook shapes with opposite signs.
Both \(D\) and \(q_3^2\) are invariant under diagram conjugation. Hence
\(\phi(\alpha_r\mathsf M_3(z)^2)=\phi(\alpha_r)=0\) for \(r=2,4\),
which proves \(G_{(2),N}=G_{(4),N}=0\) exactly.

Write \(r_N=I_{N-1}/I_N\) and \(t_j=\mathcal T_N(z^j)\).
The involution recurrence gives
\[
 t_0=1,\qquad t_1=Nr_N,\qquad
 t_j=(N-j+2)t_{j-2}-t_{j-1}\quad(2\le j\le8).
\]
Indeed, substitute
\(t_j=\fall{N}{j}I_{N-j}/I_N\) in
\(I_{N-j+2}=I_{N-j+1}+(N-j+1)I_{N-j}\).
Using these finite recurrences in the three exact polynomials of
Lemma~\ref{lem:response-exact-second-moments} gives
\begin{align*}
 G_{(3),N}
 &=\frac{432}{N^2}
   +\frac{216r_N}{(N-1)(N-2)}\\
 &\quad-\frac{36\bigl(N^2r_N^2-8Nr_N+35N+58r_N-80\bigr)}
                 {N^2(N-1)(N-2)},\\
 G_{(2,2),N}
 &=\frac{288}{N^2}
   +\frac{144\bigl((N^2-N)r_N^2+(16N-17)r_N-9N+10\bigr)}
                 {N^2(N-1)(N-2)(N-3)}.
\end{align*}
These are exact rational identities. For clarity, they follow from
\[
 \mathcal I_{\varnothing,N}(x_3^2)
       =\frac{36(t_2+2t_3+t_4)}{N^2},\quad
 \mathcal I_{(3),N}(1)=\frac{t_3}{\fall{N}{3}},\quad
 \mathcal I_{(2,2),N}(1)=\frac{2t_4}{\fall{N}{4}},
\]
and the corresponding inserted second-moment polynomials, followed by
the subtraction defining \(G_{\rho,N}\).
The analytic coefficient estimate gives
\[
 Nr_N=\mathcal T_N(z)=\sqrt N+\tau_N(h)
       =\sqrt N-\frac12+O(N^{-1/2}),
\]
so
\(r_N=N^{-1/2}-(2N)^{-1}+O(N^{-3/2})\).
The remainder numerator in the first exact formula is \(O(N)\),
and the second remainder numerator is also \(O(N)\). Consequently
\begin{align}
 G_{(3),N}&=432N^{-2}+216N^{-5/2}+O(N^{-3}),
                  \label{eq:response-G3-refined}\\
 G_{(2,2),N}&=288N^{-2}+O(N^{-4}).
                  \label{eq:response-G22-refined}
\end{align}

The weight-five insertion types are \((5)\) and \((3,2)\).
Their leading Gaussian covariances are
\[
 \operatorname{Cov}(\Xi_5,9\Xi_2^2)=0,\qquad
 \operatorname{Cov}(\Xi_3\Xi_2,9\Xi_2^2)=0,
\]
by independence and the vanishing odd moments of \(\Xi_2\sim \mathcal N(0,4)\).
The inserted estimate \eqref{eq:adv-inserted-limit}, including the
ordinary estimate at \(\rho=\varnothing\), consequently gives
\(G_{\rho,N}=O(N^{-3})\) for these two types. All insertions of
weight at least six satisfy the same bound. The finite denominator
in \eqref{eq:response-refined-exact-subtraction} is
\(1+O_\mu(N^{-1})\); its reciprocal changes the order-\(N^{-2}\)
numerator only by \(O_\mu(N^{-3})\). Equations
\eqref{eq:response-G3-refined}--\eqref{eq:response-G22-refined} yield
\begin{align*}
 \mathbb E_{\mu,N}Y_{3,N}^2-\mathbb E_{0,N}Y_{3,N}^2
  ={}&(432u_{\mu,(3)}+288u_{\mu,(2,2)})N^{-2}\\
     &+216u_{\mu,(3)}N^{-5/2}+O_\mu(N^{-3}).
\end{align*}

We also need the mean one order more precisely than in the general
correction theorem. Equation~\eqref{eq:response-linear-insertion-exact}
and \(I''(z)-I(z)=(1+z)^2I(z)\) give the exact formula
\[
 \mathcal I_{(2),N}(x_3)
   =\frac{12}{N \fall{N}{2}}
                         \mathcal T_N\bigl(z^2(1+z)^2\bigr)
   =\frac{12}{N}\left(1-\frac{I_{N-2}}{I_N}\right).
\]
The ratio \(I_{N-2}/I_N\) is \(O(N^{-1})\).
Diagram conjugation makes the weight-three insertion of the odd test
\(x_3\) zero exactly, while every insertion of weight at least four
is \(O(N^{-2})\). The ordinary mean is zero. The exact insertion
ratio therefore gives
\[
 \mathbb E_{\mu,N}Y_{3,N}=12u_{\mu,(2)}N^{-1}+O_\mu(N^{-2}).
\]
Subtracting the square of this mean changes the order-\(N^{-2}\)
coefficient by \(-144u_{\mu,(2)}^2\), and contributes no
order-\(N^{-5/2}\) term. Finally,
\[
 432u_{\mu,(3)}+288u_{\mu,(2,2)}-144u_{\mu,(2)}^2
       =-72m(m-1),\qquad u_{\mu,(3)}=\kappa_3(\mu),
\]
by Lemmas~\ref{lem:adv-first-schur-coefficients} and
\ref{lem:response-weight-four}. This proves the result.
\end{proof}

\begin{example}[The bases $(2)$ and $(2,1)$]
For $\mu=(2)$, $s_\mu=(p_1^2+p_2)/2$, so $s_*=2$ and
$A_{\mu,2}=\Xi_2$. In fact,
\[
 \mathbb E_{(2),N}P(Y_N)-\mathbb E_{0,N}P(Y_N)
                  =\mathcal I_{(2),N}(P)
\]
exactly, because $\mathcal I_{(2),N}(1)=0$. Thus the coefficient of
$N^{-1}$ is $12$ for $P=x_3$ and $1296$ for $P=x_3^3$.
For $\mu=(2,1)$,
\[
 s_\mu=(p_1^3-p_3)/3,\qquad s_*=3,\qquad A_{\mu,3}=-\Xi_3.
\]
The correction is
$-N^{-3/2}\operatorname{Cov}(\Xi_3,P(Y))+O(N^{-2})$;
for $P=x_4$ its coefficient is $-24$. These two leading corrections
are pure infinitesimal mean translations and disappear under separate
exact expectation centering.
\end{example}

\begin{example}[Vanishing of the first two response orders]
For $\mu=(3,1,1)$,
\[
 s_\mu=\frac{p_1^5}{20}-\frac{p_1p_2^2}{4}+\frac{p_5}{5}.
\]
Hence $C_1(\mu)=\kappa_3(\mu)=0$, and the first nonzero insertion has
weight four, with
$A_{\mu,4}(\Xi)=-5(\Xi_2^2-2)$. Theorem~\ref{thm:response-variance-refined} gives
\[
 \operatorname{Var}_{\mu,N}(Y_{3,N})
       -\operatorname{Var}_{0,N}(Y_{3,N})
       =-\frac{1440}{N^2}+O(N^{-3}).
\]
Both means vanish exactly by diagram conjugation symmetry. This correction
therefore remains after exact expectation centering.
\end{example}

\paragraph{Centering and scope.}
The classification in Theorem~\ref{thm:adv-first-nonzero-inner}
uses the same deterministic centering \(a_dN^{d/2}\) under both laws.
Separate exact expectation centering removes the order-\(N^{-1}\)
term in \eqref{eq:adv-first-inner}. When \(C_1(\mu)=0\), it also
removes the order-\(N^{-3/2}\) term in \eqref{eq:adv-second-inner}.
The order-\(N^{-2}\) term in the third case of
Corollary~\ref{cor:response-three-levels} has the nonzero quadratic
residual \eqref{eq:response-fourth-residual}.

For transformations of the observables at a common outer size \(N\),
the exact affine relation is
\[
 \frac{q_d(\lambda)-q_d(\nu)-a_dt^{d/2}}{t^{(d-1)/2}}
 =
 \left(\frac Nt\right)^{(d-1)/2}Y_{d,N}
 +\frac{a_d(N^{d/2}-t^{d/2})-q_d(\nu)}{t^{(d-1)/2}},
\]
where \(t>0\) is the chosen scale and \(\nu=\varnothing\) or \(\mu\)
chooses whether to subtract the inner corner sum. If the two laws use
different affine transformations, the terms in this identity must
be included before comparing high-order coefficients. Separate exact
expectation centering cancels any deterministic translation in the
identity, including the subtraction of \(q_d(\mu)\). A change of
scale still multiplies the centered statistic by
\((N/t)^{(d-1)/2}\), and its variance by \((N/t)^{d-1}\).

\section{Conclusion and future directions}\label{sec:conclusion}

We developed a local difference-operator calculus for polynomial
observables on Young's lattice. Pointwise nilpotence and highest
symbols give exponential polynomiality, degree bounds and explicit
highest differences for sums with a fixed inner partition. We obtain
complete factorial content and hook coefficients and prove integrality
of their derivatives. The same \(D\)-identities yield the limit shape,
joint Gaussian fluctuations and quantitative inner-shape corrections.
Taking \(\mu=\varnothing\) recovers the classical Gelfand limit laws.

A natural direction is to extend this approach to general
\(\beta\)-Plancherel measures \cite[(47)]{MeliotKerov}:
\[
 \mathbb P_{N,\beta}(\lambda)
   =\frac{(f^\lambda)^\beta}{\mathcal Z_{N,\beta}},\qquad
 \mathcal Z_{N,\beta}=\sum_{\nu\vdash N}(f^\nu)^\beta,
 \qquad \lambda\vdash N,\quad \beta>0.
\]
The cases \(\beta=1\) and \(\beta=2\) are the Gelfand and Plancherel
measures, respectively.

We expect suitably weighted difference operators to provide a useful
framework for function averages under this family. One objective is to
identify which content, hook and shifted-part averages retain exact
polynomiality and to develop asymptotic expansions more generally.
Leading moment identities could then be used to study limit shapes,
while higher operator terms could describe Gaussian fluctuations,
\(\beta\)-dependent centering and covariance, and finer corrections.
Further questions concern perturbations of the measure and fixed inner
boundary conditions.

\begin{samepage}
\section*{Declarations}
%\subsection*{Data availability} All results are proved in the text. No external datasets were used.

\subsection*{Declaration of generative AI and AI-assisted technologies
in the manuscript preparation process}
During the preparation of this work, the authors used the AI-assisted
tools Codex and Claude for language editing, manuscript organization,
and assistance with proof strategies and symbolic calculations. The
authors reviewed and edited the resulting material and take full
responsibility for the contents of the manuscript.

\end{samepage}
\raggedbottom

\end{document}